\documentclass{amsart}

\usepackage[applemac]{inputenc}
\usepackage[T1]{fontenc}
\usepackage{amssymb}
\usepackage{amsmath}
\usepackage{amsthm}
\usepackage{varioref} 
\usepackage{paralist} 
\setdefaultenum{\upshape(i)}{}{}{}
\usepackage[pdftex]{graphicx}
\usepackage[pdftex,colorlinks,linkcolor=black,filecolor=black,citecolor=black,urlcolor=black]{hyperref}  

\allowdisplaybreaks

\numberwithin{equation}{section}

\newcommand\drawingpath{}

\newtheoremstyle{mythmstyle}%
{1.5\baselineskip}
{\baselineskip}
{\itshape}
{}
{\bf}
{}
{0pt}
{} 

\newtheoremstyle{mydefstyle}%
{1.5\baselineskip}
{\baselineskip}
{}
{}
{\bf}
{}
{0pt}
{} 

\newtheoremstyle{mypreuvestyle}%
{\baselineskip}
{\baselineskip}
{}
{}
{\em}
{}
{0pt}
{} 

\newif\ifmynonumberenvi\mynonumberenvitrue
\theoremstyle{mythmstyle}
\newtheorem{proclaimmythm}[equation]{} 
\newtheorem*{proclaimmythm*}{}
\newenvironment{proclaim}[2][*]{\ifx*#1\mynonumberenvitrue\begin{proclaimmythm*}{\bf#2.} \ignorespaces\else\mynonumberenvifalse\begin{proclaimmythm}{.\kern0.5em\bf#2.}\label{#1} \ignorespaces\fi}{\ifmynonumberenvi\end{proclaimmythm*}\else\end{proclaimmythm}\fi}

\theoremstyle{mydefstyle}
\newtheorem{proclaimmydef}[equation]{} 
\newtheorem*{proclaimmydef*}{}
\newenvironment{definition}[2][*]{\ifx*#1\mynonumberenvitrue\begin{proclaimmydef*}{\bf#2.}\else\mynonumberenvifalse\begin{proclaimmydef}{.\kern0.5em\bf#2.}\label{#1} \ignorespaces\fi{}}{\ifmynonumberenvi\end{proclaimmydef*}\else\end{proclaimmydef}\fi}

\def\QEDbox{\hbox{\lower2.3pt\vbox{\hrule\hbox
   {\vrule\kern1pt\vbox{\kern1.7pt\hbox{$\scriptstyle
   QED$}\kern.6pt}\kern1pt\vrule}\hrule}}}
\def\QED{\hskip0.01em plus 40pt\null{} \null\nobreak\hfill
   \kern3pt\QEDbox} 
\newcommand\QEDici{\\\noalign{\vskip-\baselineskip\smash{\hbox to\linewidth{\vrule width0pt \hfill\global\QEDdejaplacetrue\QEDbox}}\vskip-\baselineskip}}
\newif\ifQEDdejaplace\QEDdejaplacefalse

\theoremstyle{mypreuvestyle}
\newtheorem*{proclaimmypreuve}{}
\newenvironment{preuve}[1][*]{\begin{proclaimmypreuve}{\ifx*#1{\em Proof.}\else{\em#1.}\fi} \ignorespaces}{\ifQEDdejaplace\global\QEDdejaplacefalse\else\QED\fi\end{proclaimmypreuve}}

\def\eindelangbewijs{}
\long\def\langbewijs#1\eindelangbewijs{}
\def\langbewijs{}

\def\eindekortbewijs{}
\long\def\kortbewijs#1\eindekortbewijs{}

\newcommand{\recalf}[1]{(\ref{#1})}
\newcommand{\recalt}[1]{{[\ref{#1}]}}
\newcommand{\recaltt}[2]{{[\ref{#2}.\ref{#1}]}}

\def\dcap_#1{\mathchoice{%
          {\textstyle\bigcap\limits_{#1}}}%
          {\underset{#1}\cap}%
          {\underset{#1}\cap}%
          {\underset{#1}\cap}}
\def\dcup_#1{\mathchoice{%
          {\textstyle\bigcup\limits_{#1}}}%
          {\underset{#1}\cup}%
          {\underset{#1}\cup}%
          {\underset{#1}\cup}}
\def\ddcap_#1^#2{\mathchoice{%
          {\textstyle\bigcap\limits_{#1}^{#2}}}%
          {\overset{#2}{\underset{#1}\cap}}%
          {\overset{#2}{\underset{#1}\cap}}%
          {\overset{#2}{\underset{#1}\cap}}}
\def\ddcup_#1^#2{\mathchoice{%
          {\textstyle\bigcup\limits_{#1}^{#2}}}%
          {\overset{#2}{\underset{#1}\cup}}%
          {\overset{#2}{\underset{#1}\cup}}%
          {\overset{#2}{\underset{#1}\cup}}}

\newcommand\bigrestricted{{\kern1pt\vrule height3.3ex depth1.7ex width0.6pt\kern1pt}}
\newcommand\caprestricted{{\kern1pt\vrule height1.7ex depth0.5ex width0.6pt\kern1pt}}
\newcommand\midrestricted{{\kern1pt\vrule height1.7ex depth0.9ex width0.6pt\kern1pt}}
\newcommand\NN{\mathbf{N}}
\newcommand\restricted{{\kern1pt\vrule height1.3ex depth0.5ex width0.6pt\kern1pt}}
\def\RR{\mathbf{R}}

\newcommand\baseleaftopol{\mathcal{B}}
\newcommand\body{\mathbf{B}}
\newcommand\CA{\mathcal{A}}
\newcommand\CC{\mathbf{C}}

\newcommand\contrf[1]{\iota(#1)}

\newcommand\Domexp{{\mathcal{D}\kern-1.5pt\mathcal{E}}}
\newcommand\eexp{{\operatorname{e}}}
\newcommand\extder{\mathrm{d}}
\newcommand\foliation{\mathcal{F}}
\def\fracp#1#2{\frac{\partial #1}{\partial #2}}
\def\ie{i.e.}
\font\klein=cmr5 at 4pt
\newcommand\kleinlessthan{\hbox{\klein>}}
\newcommand\Imexp{{\mathcal{I}\kern-1.5pt\mathcal{E}}}
\newcommand\leaftopol{{\mathcal{T}_\foliation}}
\def\Liealg#1{\mathfrak#1}
\def\mapob{\ }
\newcommand\mb{{\overline m}}
\def\mo{^{-1}}
\newcommand\nbhdbase{\mathcal{V}}
\newcommand\nilpotent{\mathcal{N}}
\newcommand\Oc{{\check O}}
\newcommand\Oh{{\widehat O}}
\newcommand\Ot{{\widetilde O}}
\newcommand\QQ{\mathbf{Q}}
\def\quote#1{``#1''}
\def\scirc{\,{\raise 0.8pt\hbox{$\scriptstyle\circ$}}\,}
\newcommand\slice{plaque}
\newcommand\smooths{\mathbf{F}}
\newcommand\stresd[1]{\emph{#1}}
\newcommand\stress[1]{{\bf#1}}
\newcommand\Ub{{\overline U}}
\newcommand\Uh{{\widehat U}}
\newcommand\Ut{{\widetilde U}}
\newcommand\varphit{{\widetilde\varphi}}
\newcommand\Vh{{\widehat V}}
\newcommand\Vt{{\widetilde V}}
\newcommand\Xb{\,{\overline X}\,}

\newcommand\Yt{{\widetilde Y}}
\newcommand\ZZ{\mathbf{Z}}

\begin{document}
\setdefaultleftmargin{2.5em}{2.1em}{1.87em}{1.7em}{1em}{1em}

\author{Gijs M. Tuynman}

\title{Integrating infinitesimal (super) actions}

\address[Gijs M. Tuynman]{Laboratoire Paul Painlev\'e, U.M.R. CNRS 8524 et UFR de Math\'ematiques,
Universit\'e de Lille I, 59655 Villeneuve d'Ascq Cedex, France}
\email{FirstName[dot]LastName[at]univ-lille1[dot]fr}

\begin{abstract}  
In this paper we generalize some results of Richard Palais to the case of Lie supergroups and Lie superalgebras. 
More precisely, let $G$ be a Lie supergroup, $\Liealg g$ its Lie superalgebra and let $\rho$ be an infinitesimal action (a representation) of $\Liealg g$ on a supermanifold $M$. 
We will show that there always exists a local (smooth left) action of $G$ on $M$ such that $\rho$ is the map that associates the fundamental vector field on $M$ to an algebra element (we will say that the action integrates $\rho$).
We also show that if $\rho$ is univalent, then there exists a unique maximal local action of $G$ on $M$ integrating $\rho$.
And finally we show that if $G$ is simply connected and all (smooth, even) vector fields $\rho(X)$ are complete then there exists a global (smooth left) action of $G$ on $M$ integrating $\rho$.
Omitting all references to the super setting will turn our proofs into variations of those of Palais.

\end{abstract}

\maketitle

\section{Introduction}

Let $D\subset G\times M$ be an action domain (precise definitions will be given later) for a local smooth left-action $\Psi:D\to M$ of a Lie group $G$ on a manifold $M$ and let $\Liealg g$ be the Lie algebra of $G$ (seen as left-invariant vector fields on $G$). Associated to $X\in \Liealg g$ we can define the fundamental vector field $X^M$ on $M$ by
$$
X^M\caprestricted_m = \frac{d}{dt}\bigrestricted_{t=0} \Psi\bigl(\exp(-tX), m\bigr)
\mapob.
$$
The map $\rho:X\mapsto X^M$ from the Lie algebra $\Liealg g$ to (smooth) vector fields on $M$ is a Lie algebra homomorphism: $[X,Y]^M = [X^M,Y^M]$ (to have this without a minus sign we introduced the minus sign in the definition of the fundamental vector field $X^M$).
Said differently, the map $\rho$ is a representation of the Lie algebra $\Liealg g$ by vector fields on $M$.
We will describe this situation by saying that the (local left-) action $\Psi$ (of $G$ on $M$) \stresd{integrates} the infinitesimal action $\rho$ (of $\Liealg g$ by (smooth) vector fields on a manifold $M$, where ``infinitesimal action'' is another word for ``Lie algebra representation by vector fields'').

We thus have shown that for any local action $\Psi$ there is an infinitesimal action $\rho$ such that $\Psi$ integrates $\rho$.
Conversely, one can ask the question whether for a given infinitesimal action $\rho$ there exists a local action that integrates it. 
It is more or less obvious that if $\Psi:D\to M$ is a local action and if $D'\subset D$ is also an action domain, then the restriction of $\Psi$ to $D'$ is again a local action with the same infinitesimal action. 
This implies that uniqueness of a local action must be subordinated to the action domain.

As in the non-super case, we will show (using the same kind of arguments as in \cite{Pa57}) that for Lie superalgebras and Lie supergroups we have existence and uniqueness of local actions integrating a given infinitesimal action $\rho$. We also have the same criterium for the existence of a (unique) maximal local action  and the same criterium for the existence of a global action. 
However, in order to understand these criteria, it is useful to have examples that show what can go wrong, examples that explain the why of some subtleties of the proofs.
To do so, we start by stating some properties that a local and\slash or global action has.

So let us suppose that $\Psi:D\subset G\times M\to M$ is a local action and choose a fixed $m_o\in M$. For this $m_o$ we define $D_{m_o}\subset G$ by
$$
D_{m_o} = \{\, g\in G \mid (g,m_o)\in D\,\}
$$ 
and we consider the graph of the map $D_{m_o}\to M$, $g\mapsto \Psi(g,m_o)$, \ie, we consider the submanifold
$$
\Lambda_{m_o} = \{\,(g,m)\in G\times M \mid (g,m_o)\in D\ \&\ m=\Psi(g,m_o)\,\}
\subset G\times M
\mapob.
$$
It is immediate that the tangent space to this submanifold is given by
$$
T_{(g,m)}\Lambda_{m_o} = \{\, X^r\caprestricted_g - X^M\caprestricted_m \mid X\in \Liealg g\,\}
\mapob,
$$
where $X^r$ denotes the right-invariant vector field on $G$ whose value at the identity $e\in G$ is $X$: $X^r\caprestricted_e = X\in T_eG\cong \Liealg g$. It is also immediate that the projection $p:G\times M\to G$, $p(g,m)=g$ defines a diffeomorphism between $\Lambda_{m_o}$ and $D_{m_o}$.
But we can say more: if we define the foliation $\foliation\subset T(G\times M)$ on $G\times M$ by
$$
\foliation\caprestricted_{(g,m)} = \{\, X^r\caprestricted_g - X^M\caprestricted_m \mid X\in \Liealg g\,\}
\mapob,
$$
then $\Lambda_{m_o}$ is an integral manifold of $\foliation$. 
As $D_{m_o}$ is connected, $\Lambda_{m_o}$ is contained in the leaf $L_{(e,m_o)}$ through $(e,m_o)$ (the unique maximal connected integral manifold of $\foliation$ containing $(e,m_o)\,$).
One can show that we have the following properties:
\begin{enumerate}
\item\label{discuLambdaconnandopen}
$\Lambda_{m_o}$ is a connected and open subset of $L_{(e,m_o)}$ containing $(e,m_o)$, 

\item\label{discuLambdabijection}
$p:\Lambda_{m_o}\to D_{m_o}$ is a bijection,

\item
$D= \cup_{m_o\in M} D_{m_o}\times \{m_o\}$, and

\item
$m=\Psi(g,m_o) \quad\Longleftrightarrow\quad (g,m)\in \Lambda_{m_o} 
\quad\Longleftrightarrow\quad m=(p\restricted_{\Lambda_{m_o}})\mo(g)$.

\end{enumerate}
We thus see that $D$ and $\Psi$ are completely determined by the sets $\Lambda_{m_o}$ and that these sets in turn are more or less determined by the foliation $\foliation$ (the main constraint being that the map $\Psi$ defined by the $\Lambda_{m_o}$ should satisfy the group law).

It is more or less obvious that, given a local action $\Psi:D\to M$, any smaller action domain $D'\subset D$ also defines, by restriction, a local action $\Psi'=\Psi\caprestricted_{D'}:D'\to M$. Hence a local action will never be unique. 
On the other hand, if the leaves of $\foliation$ project injectively to $G$ (this is the notion of being univalent), then we can choose $\Lambda_{m_o}=L_{(e,m_o)}$ and the associated local action will be maximal in the sense that any other local action has a smaller action domain.

A particular case in which the leaves of $\foliation$ project injectively to $G$ is when the action is global. 
More precisely, it is relatively easy to show (or immediate) that for a global we have the equalities $D_{m_o}=G$ (which implies that $p$ maps $\Lambda_{m_o}$, and ipso facto $L_{(e,m_o)}$, onto $G$) and $L_{(e,m_o)}=\Lambda_{m_o}$.
It follows immediately that if a leaf of $\foliation$ passing through some $(e,m_o)$ does not project bijectively onto $G$, then a global action can not exist.
Moreover, the definition of the fundamental vector field tells us that the flow $\Phi_{X^M}$ of $X^M$ is given by
$$
\Phi_{X^M}(t,m) = \Psi\bigl(\exp(-tX), m\bigr)
\mapob.
$$
This shows that, if the action is global, then necessarily the flow of the vector field $X^M$ is complete, \ie, defined for all $t\in \RR$.
It follows that if one of the vector fields $\rho(X)$ has a non-complete flow, then there can not exist a global action integrating $\rho$.

\begin{definition}[examp1]{{Example}}
Let $G=\RR$ be the additive Lie group of real numbers with coordinate $g\in \RR$, let $M=(0,1)\subset \RR$ be the open unit interval in $\RR$ with coordinate $x\in (0,1)$. The Lie algebra $\Liealg g$ of $G$ is (isomorphic to) $\RR$ and can be seen as the set of left-invariant vector fields $a\cdot \partial_g$, $a\in \RR$ on $G=\RR$. Associated to the Lie algebra element $a\cdot \partial_g$ we define the vector field $\rho(a\cdot\partial_g) = -a\cdot\lambda\cdot \partial_x$ on $M$ for some fixed constant $\lambda\in \RR$. As the Lie algebra is $1$-dimensional, this is automatically an infinitesimal action.

The foliation $\foliation$ on $G\times M$ is given by
$$
\foliation\caprestricted_{(g,x)}
=
\RR\cdot (\partial_g + \lambda\partial_x)
$$
and the leaf $L_{(g,x)}$ passing through $(g,x)$ is given by
$$
L_{(g,x)} = \{\,(g+t,x+\lambda t) \mid -x< \lambda t<1-x \,\}
\ .
$$
In this case the projection $p:G\times M\to G$ restricted to a leaf is not a bijection (not surjective), unless $\lambda=0$. Hence this infinitesimal action $\rho$ cannot be integrated to a global left-action of $G$ on $M$ unless $\lambda=0$, in which case it is obtained by the trivial action $\Psi(g,x) = x$.

As the leaves $L_{(e,x)}$ project injectively to $G$, the infinitesimal action can be integrated to a (unique) maximal local action on the domain $D\subset G\times M$ given by
\begin{align*}
D
&
=\{\, (g,x)\in G\times M\mid g\in p(L_{(e,x)}) \,\}
\\
&
=
\{\,(t,x)\in G\times M \mid -x< \lambda t<1-x\,\}
\mapob.
\end{align*}
On this domain the local action $\Psi:D\to M$ is given by
$$
\Psi(t,x) = x+\lambda t
\mapob.
$$
One could note that these formul{\ae} include the trivial case $\lambda=0$ with its global action.

\end{definition}

\begin{definition}[examp2]{Example}
If in example \recalt{examp1} we change the group $G$ to the circle group $G=\RR/\ZZ$, we can keep the rest unchanged. The formula defining $\foliation$ is unchanged, but the leaf $L_{(g,x)}$ passing through $(g,x)$ is given by 
$$
L_{(g,x)} = \{\,(g+t+\ZZ,x+\lambda t) \mid -x<\lambda t<1-x \,\}
\ .
$$
As $G\times M$ is a (finite) cylinder, it is not hard to see that the leaves are parts of helices with slope $\lambda$. 
\begin{figure}[htb]
\null\hfil
\includegraphics[width=0.3\textwidth]{\drawingpath 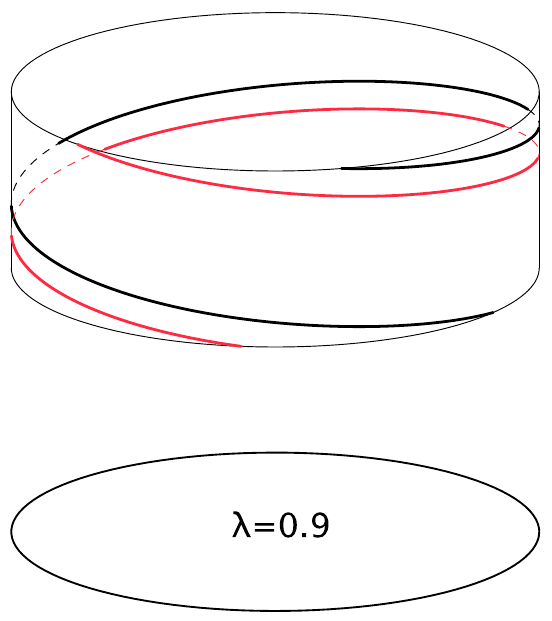}
\hfil\hfil\hfil
\includegraphics[width=0.3\textwidth]{\drawingpath 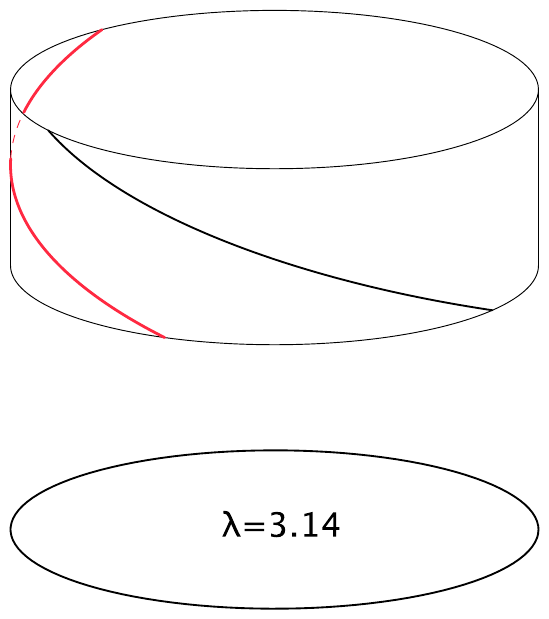}
\hfil
\null
\end{figure}

\noindent
The projection $p:G\times M\to G$ restricted to a leaf will be surjective for $\vert\lambda\vert<1$, but then it will not be injective unless $\lambda=0$. On the other hand, for $\vert\lambda\vert\ge1$ it will be injective but not surjective. Apart from the trivial case $\lambda=0$, this infinitesimal action of $\Liealg g$  thus never can be integrated to a global left-action of $G$ on $M$.

For $\vert \lambda\vert\ge1$ the leaves $L_{(e,x)}$ project injectively to $G$ and the infinitesimal action $\rho$ can be integrated to a (unique) maximal local action defined on $D\subset G\times M$ given by
\begin{align*}
D
&
=\{\, (g,x)\in G\times M\mid g\in p(L_{(e,x)}) \,\}
\\
&
=
\{\,(t+\ZZ,x)\in G\times M \mid -x<\lambda t<1-x\,\}
\end{align*}
On this domain the local action $\Psi:D\to M$ is given by
$$
\Psi(t+\ZZ,x) = x+\lambda t
\mapob,
$$
where $t\in \RR$ representing the element $t+\ZZ\in G=\RR/\ZZ$ must be chosen such that $0<x+\lambda t<1$ (for $(t+\ZZ,x)\in D$ such a choice exists and is unique).

\end{definition}

\begin{definition}[examp3]{Example}
If in example \recalt{examp2} we change the manifold $M$ to $M=\RR$, then the leaf $L_{(g,x)}$ passing through $(g,x)$ is given by 
$$
L_{(g,x)} = \{\,(g+t+\ZZ,x+\lambda t) \mid t\in \RR \,\}
\ .
$$
As now $G\times M$ is an infinite cylinder, the projection $p:G\times M$ restricted to a leaf will always be surjective onto $G$, but (apart from the case $\lambda=0\,$) it will never be injective, and thus a global action integrating $\rho$ will never exist.
On the other hand, the projection restricted to a leaf is (always) a covering map (of the circle by the real line for $\lambda\neq0$; it is the identity map for $\lambda=0$).

\end{definition}

\begin{definition}[examp4]{Example}
If in example \recalt{examp2} we change the manifold $M$ to the circle $M=\RR/\ZZ$, the formul{\ae} still do not change, but then the leaf $L_{(g,x)}$ passing through $(g,x)$ is given by 
$$
L_{(g,x)} = \{\,(g+t+\ZZ,x+\lambda t+\ZZ) \mid t\in \RR \,\}
\ .
$$
And now the topological nature of a leaf depends upon the value of $\lambda\in \RR$: for $\lambda\in \QQ$ it will be a circle and for $\lambda\in \RR\setminus \QQ$ it will be the real line. 
\begin{figure}[htb]
\null\hfil
\includegraphics[width=0.3\textwidth]{\drawingpath 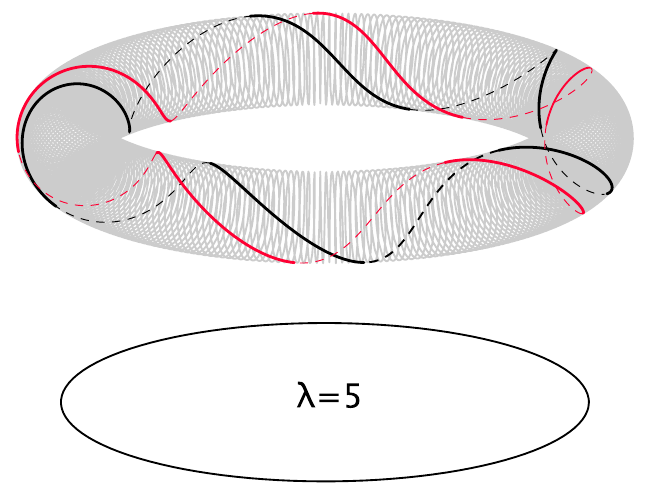}
\hfil\hfil\hfil
\includegraphics[width=0.3\textwidth]{\drawingpath 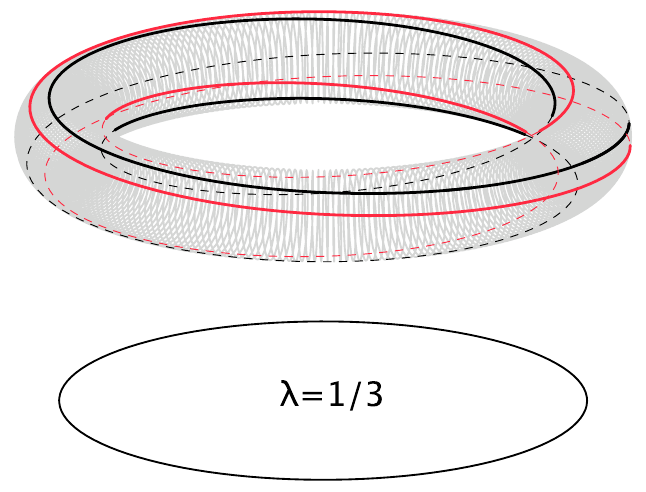}
\hfil
\null
\end{figure}

\noindent
In all cases the projection $p:G\times M$ restricted to a leaf will be surjective onto $G$. However, for $\lambda\in \RR\setminus\QQ$ it will be (equivalent to) the standard universal covering map from $\RR$ to the circle $\RR/\ZZ$. For $\lambda\in \QQ$ it will be a finite covering map from the circle to itself, which is a bijection if and only if $\lambda\in \ZZ$. It follows that for $\lambda\notin \ZZ$ there does not exist a global left-action of $G$ on $M$ integrating this infinitesimal action. On the other hand, it is easy to see that for $\lambda=n\in \ZZ$, the (global) left-action $\Psi:G\times M\to M$ defined by
$$
\Psi(g,x) = (x+ng)
$$
integrates the infinitesimal action.

\end{definition}

\begin{definition}[examp5]{Example}
If in examples \recalt{examp3} and \recalt{examp4} we change the group $G$ to $G=\RR$, then the projection $p:G\times M\to G$, when restricted to a leaf, will be a bijection for all $\lambda\in \RR$. And it is easy to show that the infinitesimal action of $\Liealg g \cong \RR$ on $M$ is integrated to the global left-action of $G=\RR$ on $M$ given by
\begin{multline*}
\Psi(g,x) = x+\lambda g
\text{ when $M=\RR$}
\quad\text{ or }
\\
\quad
\Psi(g,x) = x+\lambda g+\ZZ
\text{ when $M=\RR/\ZZ$.}
\end{multline*}

\end{definition}

Looking at the examples \recalt{examp1} to \recalt{examp5}, one could make the following observations. For an infinitesimal action $\rho$ on a fixed manifold $M$, the choice of a Lie group with Lie algebra $\Liealg g$ influences the existence of a (unique) maximal local action: if the group is simply connected there always is a maximal local action and if it is not, only (very) particular cases integrate to a global action.
Comparing the examples \recalt{examp1} and \recalt{examp2} with the examples \recalt{examp3}, \recalt{examp4} and \recalt{examp5} shows another phenomenon: in the last three cases the fundamental vector fields are complete and the projection of a leaf of $\foliation$ to $G$ is a covering map, whereas in the first two cases these projections are not coverings and neither are the fundamental vector fields complete (giving another argument agains the existence of a global action integrating $\rho$).
The idea that comes to mind is that the manifold $M$ in the first two examples is ``too small'' and that one should enlarge it in such a way that the vector fields become complete. And that is exactly what is done in the last three examples.
This line of reasoning suggests that, if the vector fields given by the infinitesimal action are not complete, then we search for a bigger manifold on which we still have an infinitesimal action, but one on which the vector fields are complete. 
And then we look at the simply connected Lie group associated to $\Liealg g$ and we obtain a global action. Restriction to the original manifold then would give a local action.
That the situation is not as simple as this suggests is shown in our last example.

\begin{definition}[examp6]{Example}
In this last example, we change the setting to higher dimensions by looking at the classical example of two commuting vector fields whose flows do not commute \cite[p353]{LM87}.\footnote{At the same time it is an explicit case of the generic example in \cite[Thm XXI, p88]{Pa57}.}
For the manifold $M$ we take $\RR^2$ minus a square:
$$
M = \{\,(x,y)\in \RR^2\mid \max(\,\vert x\vert, \vert y\vert\,)>1\,\}
\mapob.
$$
On $M$ we define the two vector fields $X$ and $Y$ by
$$
X\caprestricted_{(x,y)} = \fracp{}{x}\bigrestricted_{(x,y)}
\quad\text{and}\quad
Y\caprestricted_{(x,y)}
=
\begin{cases}
\displaystyle\fracp{}{y}\bigrestricted_{(x,y)}
&\quad x\ge -1 
\\
\noalign{\vskip2\jot}
\displaystyle\chi'\bigl(\chi\mo(y)\bigr)\cdot \fracp{}{y}\bigrestricted_{(x,y)}
&\quad x<-1 
\end{cases}
\mapob,
$$
where $\chi:\RR\to \RR$ is a smooth map satisfying $\chi'\ge1$ everywhere as well as
$$
t\le -1 \quad\Rightarrow\quad \chi(t) = t
\qquad\text{and}\qquad
t\ge 0 \quad\Rightarrow\quad \chi(t) = t+1
\mapob.
$$
In particular $Y$ is everywhere the constant vector field $\partial_y$, except in the region to the left of the (excluded) square $\max(\,\vert x\vert, \vert y\vert\,)\le 1$, where it is a multiple of $\partial_y$, a multiple which is, on average, greater than $1$. 

\begin{figure}[htb]
\null\hfil
\includegraphics[width=0.3\textwidth]{\drawingpath 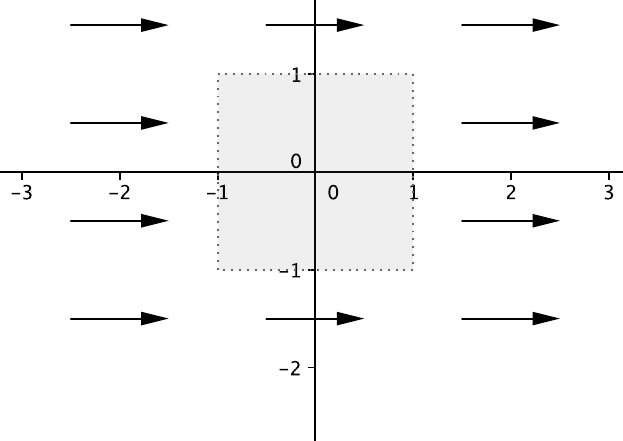}
\hfil\hfil\hfil
\includegraphics[width=0.3\textwidth]{\drawingpath 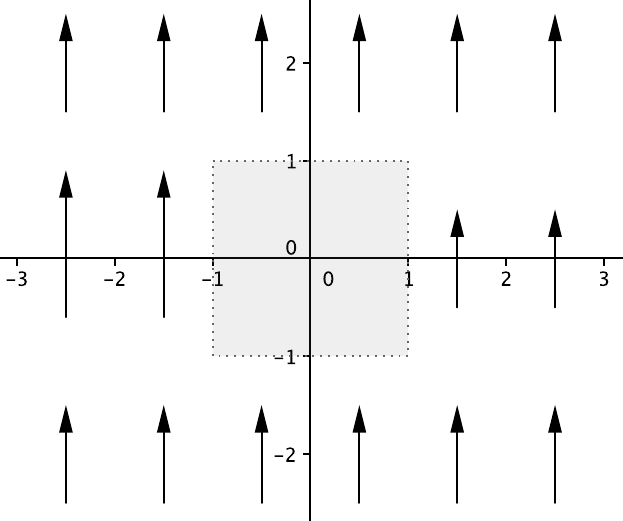}
\hfil
\null

\vskip-\baselineskip

{\null\hfil
The vector field $X$ \hfil\hfil\hfil The vector field $Y$ \hfil}
\end{figure}

It is immediate that $[X,Y]=0$, and thus these two vector fields define an infinitesimal action of the $2$-dimensional abelian Lie algebra $\Liealg g=\RR^2$ on $M$. On the other hand, it is an elementary exercise (due to the particular form of $Y$) to show that the flows $\Phi_X$ and $\Phi_Y$ (of $X$ and $Y$ respectively) are given by the formul{\ae}
\begin{multline*}
\Phi_X(s,x,y) = (x+s, y)
\quad\text{and}\quad
\\
\Phi_Y(t,x,y) = 
\begin{cases}
(x,y+t) &\quad x\ge -1
\\
\noalign{\vskip1\jot}
(x,\chi\bigl(\chi\mo(y)+t)\bigr) &\quad x<-1
\end{cases}
\mapob.
\end{multline*}
Due to the excluded square, these flows are not complete, and hence this infinitesimal action can not be integrated to a global action.
But there is worse: these flows do not commute, as we have 
$$
\Phi_X\bigl(4,\Phi_Y(4,-2,-2)\bigr) = (2,3) \neq (2,2) = \Phi_Y\bigl(4,\Phi_X(4,-2,-2)\bigr)
\mapob.
$$
The more prosaic explanation is that the flow of $Y$ runs faster (vertically) to the left of the excluded square than it does to the right, whereas the flow of $X$ always runs (horizontally) at the same speed everywhere. 
It follows that it is not possible to enlarge $M$ to a bigger manifold on which the vector fields $X$ and $Y$ become complete in such a way that the infinitesimal action integrates to a global action of the simply connected abelian Lie group $\RR^2$.

When we search for integral manifolds and leaves of the foliation $\foliation$ which is given by
$$
\foliation\caprestricted_{(s,t,x,y)} = \RR\cdot (\partial_s\caprestricted_{(s,t)} - X\caprestricted_{(x,y)}) + \RR\cdot ( \partial_t\caprestricted_{(s,t)} - Y\caprestricted_{(x,y)})
\mapob, 
$$
the manifold $G\times M$ splits into two regions $R_1$ and $R_2$ given by
\begin{align*}
R_1 
&
= \{\, (s,t,x,y)\mid x\ge -1 \text{ or } \vert y\vert \ge 1\,\}
\\
R_2
&
= \{\, (s,t,x,y)\mid x< -1 \text{ and } \vert y\vert < 1\,\}
\mapob.
\end{align*}
On $R_1$ the maximal integral manifolds are given by the equations
$$
x+s = \text{const.} \qquad\text{and}\qquad y+t=\text{const.}
$$
whereas on $R_2$ they are given by
$$
x+s = \text{const.} \qquad\text{and}\qquad \chi\mo(y)+t=\text{const.}
$$
We thus have to join these two regions to find the leaves on $G\times M$. The dependence on $s$ being easy, we concentrate on the $(t,x,y)$ dependence, which allows us to draw $2$-dimensional pictures of this $3$-dimensional situation. 
On $R_1$ the leaves (forgetting the first coordinate!) are U-formed parts of inclined planes that are put around the ``bar'' $\{t\in\RR\}\times\{\max(\,\vert x\vert, \vert y\vert\,)\le1\}$. 
On $R_2$ the leaves are curved razor blades that are more inclined vertically than the planes in $R_1$. 
It follows that a single leaf in $R_2$ connects two ``planes'' in $R_1$. 
As a consequence, a single leaf in the full space $G\times M$ is a kind of staircase winding around the bar $\{t\in\RR\}\times\{\max(\,\vert x\vert, \vert y\vert\,)\le1\}$.
\begin{figure}[htb]
\null\hfil
\includegraphics[scale=0.88]{\drawingpath 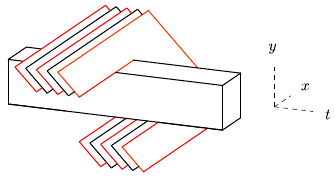}
\hfil
\includegraphics[scale=0.65]{\drawingpath 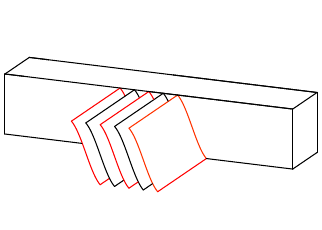}
\hfil\hfil
\includegraphics[scale=0.65]{\drawingpath 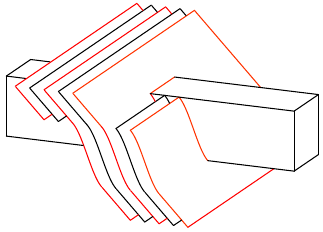}
\hfil
\null
\vskip-\baselineskip
{\null\hfil
\kern2em 5 leaves in $R_1$\hfill\hfill 
5 leaves in $R_2$ \hfill 
Just 2 leaves in $G\times M$\kern1em \hfil}
\end{figure}
It follows easily that when we restrict the projection $p:G\times M\to G$ to a leaf, it is always surjective but never injective. In particular, it will never be a covering map ($G=\RR^2$ is simply connected so any connected covering must be a bijection).

\end{definition}

We will show that these examples (in non-super geometry) are characteristic also for the super differential geometric setting, allowing us to generalize some of the results of Palais. Our main results are the following (precise definitions will follow in the main text).

\begin{proclaim}{{Theorem}}
Let $G$ be a connected Lie supergroup with Lie superalgebra $\Liealg g$ and let $\rho$ be an infinitesimal action of $\Liealg g$ on a supermanifold $M$. Then there exists an action domain $D\subset G\times M$ and a smooth local left-action $\Psi:D\to M$ integrating $\rho$.

\end{proclaim}

\begin{proclaim}{{Theorem}}
Let $G$ be a connected Lie supergroup with Lie superalgebra $\Liealg g$, let $\rho$ be an infinitesimal action of $\Liealg g$ on a supermanifold $M$ and let $\Psi_i:D\to M$, $i=1,2$ be two smooth local left-actions integrating $\rho$. Then $\Psi_1=\Psi_2$. 

\end{proclaim}

\begin{proclaim}{{Theorem}}
Let $G$ be a connected Lie supergroup with Lie superalgebra $\Liealg g$ and let $\rho$ be an infinitesimal action of $\Liealg g$ on a supermanifold $M$. If $(G,\rho)$ is univalent, then there exists a (unique) maximal smooth local left-action $\Psi:D\to M$ integrating $\rho$. 

\end{proclaim}

\begin{proclaim}{{Theorem}}
Let $G$ be a connected Lie supergroup with Lie superalgebra $\Liealg g$ and let $\rho$ be an infinitesimal action of $\Liealg g$ on a supermanifold $M$. If all vector fields $\rho(X)$ with $X\in \body \Liealg g$ are complete, then the foliation $\foliation$ is covering.

\end{proclaim}

\begin{proclaim}{{Theorem}}
Let $G$ be a simply connected Lie supergroup with Lie superalgebra $\Liealg g$ and let $\rho$ be an infinitesimal action of $\Liealg g$ on a supermanifold $M$. If all vector fields $\rho(X)$ with $X\in \body \Liealg g$ are complete, then there exists a (unique) smooth global left-action $\Phi:G\times M\to M$ integrating $\rho$.

\end{proclaim}

We will work with the geometric $H^\infty$ version of DeWitt supermanifolds, which is equivalent to the theory of graded manifolds of Leites and Kostant (see \cite{DW84}, \cite{Ko77}, \cite{Le80}, \cite{Ro07} \cite{Tu04}, \cite{Va04}) and any reader using a (slightly) different version of supermanifolds should be able to translate the results to her\slash his version of supermanifolds (transposing the proofs might be another problem altogether).
Appendix A gives some notational conventions as well as an extremely succinct overview of our version of $H^\infty$ supermanifolds.
For the moment it suffices to note that the basic graded ring is denoted $\CA$ and that $\body$ denotes the body map $\body:\CA \to \RR$ (which extends to all other super objects).

Our proofs use the same ingredients as in \cite{Pa57}, but not necessarily in the same order.
In particular we rely more heavily upon topological arguments.
There are two reasons that prohibit copying directly the proofs in \cite{Pa57}. First of all, the standard DeWitt topology in super differential geometry is not Hausdorff, so any argument based on this criterion can not be used.
We circumvent this problem by using the existence of the flow of a smooth vector field on a supermanifold (with its maximal flow domain). 
But as this property is deduced in \cite{Pa57}, we have to turn around some of the arguments.
A second reason is that for supermanifolds we have to change our viewpoint from separate objects to families.
This is because in the $H^\infty$ category of supermanifolds, fixing some variables in a smooth map does not yield (in general) a smooth map in the remaining variables.
Only when one chooses real values for some the variables (i.e., values in the body of the source) will it be guaranteed that the resulting map is smooth in the remaining variables.
This peculiarity is due to the difference between the $H^\infty$ and $G^\infty$ categories of supermanifolds. In the $H^\infty$ category the only constant functions in the graded ring of super smooth functions are reals, whereas in the $G^\infty$ category all elements of $\CA$ belong to the graded ring of smooth functions (as constant functions).
In \S\ref{sectionappendixfoliations} this difference is discussed in more detail and the reader unfamiliar with this subtlety is urged to read at least this subsection.

In the context of this paper, this peculiarity implies that not all vector fields $\rho(X)$ will be smooth (only when $X$ belongs to the body of the Lie superalgebra). Moreover, the individual maps $m\mapsto \Psi(g,m)$ for fixed $g\in G$ (i.e., the action of the fixed element $g$ on $M$) will in general not be smooth, even when $\Psi$ itself is smooth. 
A similar problem occurs with foliations: not all leaves will have the structure of an immersed submanifold. These problems are easily solved by looking systematically at families of maps, instead of individual maps (this should not come as a surprise for people working with the sheaf-theoretic\slash ringed spaces version of supermanifolds). A particularly useful example is the family of fundamental vector fields $X^M=\rho(X)$. Instead of looking at these vector fields one by one, one could look at the single vector field $Z_M$ on $\Liealg g \times M$ defined by
$$
Z_M\caprestricted_{(X,g,m)} = 0\caprestricted_X + \rho(X)\caprestricted_m
\mapob.
$$
In standard differential geometry, if $Z_M$ is a smooth vector field, then all vector fields $\rho(X)$ will be smooth. 
In super differential geometry however, even if $Z_M$ is smooth, not all individual vector fields $\rho(X)$ will be smooth. 
That will be guaranteed only when $X$ belongs to $\body\Liealg g$, the body of $\Liealg g$.
On the other hand however, when we are only interested in topological properties, restricting attention to individual members of a family does not pose any problem. 
For instance, the map $\Psi_g:M\to M$ given by $\Psi_g(m) = \Psi(g,m)$ will always be continuous, even when it is not (allowed to be called) differentiable. This explains why we will rely more heavily on topological arguments. 

We now stress that our arguments applies as well to the ordinary differential geometric context as to the super context. 
Most of the time the reader will not even be aware that there is a difference.
And when there is a difference, it suffices to suppress the words\slash symbols ``super,'' ``body,'' ``homogeneous'' and ``$\body$'' and to replace $\CA$ and $\CA_0$ by $\RR$.
This means that the proofs we give here can be seen as an alternative to\slash variation of the proofs given in \cite{Pa57}, but proofs that generalize directly to the case of supermanifolds.
In the case of ordinary differential geometry, our reasoning can be simplified slightly, using only differential geometric arguments, as some of the detours are made especially for the super differential geometric context.


We finish this introduction with the remark that if $M$ is a real analytic supermanifold (see also the end of \S\ref{subsecAmfds}) and if the vector fields $\rho(X)$ are analytic, then the (local) actions obtained in our proofs are automatically analytic. The (simple) reason is that the proofs that the actions are smooth depend only upon the fact that the flow of the vector field $Z_M$ derived from the $\rho(X)$ is smooth. And if these vector fields are analytic, then so will be the flow of $Z_M$.

\section{Generalities on foliations}

\begin{definition}{{Definition}}
Let $f:M\to N$ be a map between two topological spaces. We will say that $f$ is a \stresd{local homeomorphism} if for any $m\in M$ there exists an open neighbourhood $U\subset M$ of $m$ and an open neighbourhood $V\subset N$ of $f(m)$ such that $f:U\to V$ is a homeomorphism. 

\end{definition}

\begin{definition}{{Definition}}
Let $\foliation\subset TM$ be a foliation of rank $k$, \ie, an involutive subbundle of rank $k$. Let $(x_1, \dots, x_n)$ be a local system of coordinates on the open subset $U\subset M$. 
We will say that the local coordinates $(x_1, \dots, x_n)$ are adapted to the foliation if on $U$ the foliation is spanned by the tangent vectors $\partial_{x_1}, \dots, \partial_{x_k}$:
$$
\forall m\in U
\quad:\quad
\foliation\caprestricted_m
=
\Bigl\{\ \sum_{i=1}^k \alpha_i\cdot \partial_{x_i}\caprestricted_m \mid \alpha_i\in \CA\ \Bigr\}
\mapob.
$$
By Frobenius' theorem, around every point there exists local coordinates adapted to the foliation.
Now let $(x_1, \dots, x_n)$ be a local coordinate system on $U\subset M$ adapted to the foliation and choose $a_i\in \CA_{\varepsilon_i}$, $i=k+1, \dots, n$ (where $\varepsilon_i$ indicates the parity of the coordinate $x_i$). We then can define the \slice{} $U_{a_{\kleinlessthan k}}$ by
$$
U_{a_{\kleinlessthan k}}
=
\{\,m\in U \mid x_{k+1}(m) = a_{k+1}, \dots, x_n(m) = a_n\,\}
\mapob.
$$

\end{definition}

\begin{proclaim}{{Trivial Lemma}}
Let $(x_1, \dots, x_n)$ be local coordinates on $U\subset M$, let $a_i\in \CA_{\varepsilon_i}$ be arbitrary and let $V\subset M$ be any open subset. Then $(x_1, \dots, x_n)$ are local coordinates on $U\cap V$ and we have the equality
$$
(U\cap V)_{a_{\kleinlessthan k}} = U_{a_{\kleinlessthan k}} \cap V
\mapob.
$$

\end{proclaim}

\begin{proclaim}[slicesbasistopology]{{Lemma}}
The collection $\baseleaftopol$ of all connected components of all \slice{}s of all local coordinate systems adapted to the foliation forms the basis of a topology for $M$, a topology that is finer than the original topology of $M$.

\end{proclaim}

\langbewijs
\begin{preuve}
To prove that $\baseleaftopol$ is the basis for a topology, we have to prove the property
$$
\forall B_1,B_2\in \baseleaftopol\ \forall m\in B_1\cap B_2\ \exists B_3\in \baseleaftopol : m\in B_3\subset B_1\cap B_2
\mapob.
$$
So suppose $B_1,B_2\in \baseleaftopol$ and $m\in B_1\cap B_2$. By definition there exist local coordinate systems $(x_1, \dots, x_n)$ on an open $U\subset M$, local coordinates $(y_1, \dots, y_n)$ on an open $V\subset M$ both adapted to the foliation, and $a_i,b_i\in \CA_{\varepsilon_i}$, $i>k$, such that:
\begin{multline*}
B_1 \text{ a connected component of } U_{a_{\kleinlessthan k}}
\quad\text{and}\quad
\\
B_2 \text{ a connected component of } V_{b_{\kleinlessthan k}}
\mapob.
\end{multline*}
We claim that $B_3$ defined as the connected component of $U_{a_{\kleinlessthan k}}\cap V = (U\cap V)_{a_{\kleinlessthan k}}$ containing $m$ satisfies the requirement. It obviously belongs to $\baseleaftopol$ and we obviously have the inclusion $B_3\subset B_1$ as it is connected and included in $U_{a_{\kleinlessthan k}}$.
To prove that it is included in $V_{b_{\kleinlessthan k}}$, we note that, by definition of coordinate systems, there exist smooth maps $\psi_i$ on $U\cap V$ such that $y_i=\psi_i(x_1, \dots, x_n)$ on $U\cap V$. The fact that both coordinate systems are adapted to the foliation implies that we must have
\begin{multline*}
\Bigl(\,\fracp{\psi_j}{x_i}(m)\,\Bigr)_{i,j=1}^k
\text{ an invertible matrix}
\quad\text{and}\quad
\\
\fracp{\psi_j}{x_i}(m)
=
0
\text{ for $i\le k$ and $j>k$.}
\end{multline*}
In particular, the functions $y_j=\psi_j(x_1, \dots, x_n)$ with $j>k$ are constant on $B_3$ as it is connected and contained in $U_{a_{\kleinlessthan k}}$. But because $m\in V_{b_{\kleinlessthan k}}$, we must have $y_j(m) = b_j$ for $j>k$, and thus $B_3\subset V_{b_{\kleinlessthan k}}$. As $B_3$ is connected and contains $m$, we thus must have the inclusion $B_3\subset B_2$ as wanted.

To show that the associated topology is finer than the original topology, it suffices to make two remarks. First that we have the obvious equality
$$
U = \dcup_{a_{k+1}, \dots, a_n} U_{a_{\kleinlessthan k}}
\mapob.
$$
And second that we can take arbitrary small coordinate systems adapted to the foliation, simply by taking intersections with arbitrary open subsets (of the original topology). 
\end{preuve}
\eindelangbewijs 

\kortbewijs
\begin{preuve}
A straightforward consequence of Frobenius' theorem.
\end{preuve}
\eindekortbewijs

\begin{definition}[definitionofleaftopology]{{Definition}}
The topology on $M$ defined by the basis given in \recalt{slicesbasistopology} is called the \stresd{leaf topology} of $M$ and denoted as $\leaftopol$. In non-graded (non-super) geometry, the connected components of $M$ with respect to this topology are exactly the leaves of the foliation, \ie, the immersed maximal integral manifolds of the foliation. In super geometry some (most) of the connected components do not have the structure of an immersed submanifold, and thus are not leaves according to the official definition. That is why, in the super differential geometric context, we ``have to'' use the leaf topology: the collection of all immersed maximal integral manifolds do not fill up the manifold $M$. Via the connected components for the leaf topology we have acces to ``all'' leaves. But then we have to forego differential geometric arguments, as these connected components do not in general have the structure of an immersed supermanifold.

\end{definition}

\begin{proclaim}[localleaftopisrelativetop]{{Lemma}}
Let $(x_1, \dots, x_n)$ be a local system of coordinates on the open set $U$ adapted to the foliation. Then the topology on a \slice{} $U_{a_{\kleinlessthan k}}$ induced by the leaf topology is the same as the the one induced by the original topology on $M$.

\end{proclaim}

\langbewijs
\begin{preuve}
As the leaf topology is finer than the original topology, it follows immediately that the topology induced on $U_{a_{\kleinlessthan k}}$ by the original topology is included in the topology induced on $U_{a_{\kleinlessthan k}}$ by the leaf topology.

For the converse, choose $m\in U_{a_{\kleinlessthan k}}$ and a basic open neighbourhood $B_2$ of $m$ for the leaf topology, \ie, $B_2$ is the connected component of some \slice{} $V_{b_{\kleinlessthan k}}$. In the proof of \recalt{slicesbasistopology} we have shown that $B_3$, the connected component of $(U\cap V)_{a_{\kleinlessthan k}}$ containing $m$, is included in the intersection $U_{a_{\kleinlessthan k}} \cap B_2$. As $(x_1, \dots, x_n)$ are local coordinates on $U\cap V$ adapted to $\foliation$, it follows that there exists an open set $m\in W\subset U\cap V$ of the form $W=W_1\times W_2$ with $(x_1, \dots, x_k)$ local coordinates on $W_1$ and $(x_{k+1}, \dots, x_n)$ local coordinates on $W_2$. Moreover, we may assume that $W_1$ is connected. We thus have:
$$
W_1\times \{(a_{k+1}, \dots, a_n)\} 
=
W_{a_{\kleinlessthan k}}
\subset (U\cap V)_{a_{\kleinlessthan k}}
\mapob.
$$
As $m$ belongs to the connected set $W_1\times \{(a_{k+1}, \dots, a_n)\} 
=
W_{a_{\kleinlessthan k}}$, it follows that we have the inclusions
$$
W\cap U_{a_{\kleinlessthan k}}
=
W_{a_{\kleinlessthan k}} 
\subset 
B_3
\subset U_{a_{\kleinlessthan k}} \cap B_2
\mapob.
$$
This shows that any open neighbourhood in $U_{a_{\kleinlessthan k}}$ of $m$ for the topology induced by the leaf topology contains an open neighbourhood in $U_{a_{\kleinlessthan k}}$ of $m$ for the topology induced by the original topology of $M$. And hence the topology on $U_{a_{\kleinlessthan k}}$ induced by the original topology is finer than the one induced by the leaf topology.
\end{preuve}
\eindelangbewijs

\kortbewijs
\begin{preuve}
A straightforward consequence of Frobenius' theorem.
\end{preuve}
\eindekortbewijs

\begin{definition}{{Remark}}
While it is true that the topology on a \slice{} $U_{a_{\kleinlessthan k}}$ is the same whether induced by the original topology or by the leaf topology, it is \stress{not} true that the topology on a connected component of the leaf topology is induced by the original topology. It suffices to think of the torus $(\RR/\ZZ)^2$ foliated by lines with an irrational slope. The connected components are the leaves of this foliation, which are homeomorphic to the real line. But as these lines are dense in the torus, the topology induced on such a leaf by the topology of the torus is not the same as the topology of the real line: no finite interval on the line can be the intersection of an open set in the torus with a leaf.
 
\end{definition}

\section{Global, local and infinitesimal actions}

\begin{definition}{Notation and conventions}
Starting this section, $M$ will always be a supermanifold of total dimension $n$, $G$ will be a connected Lie supergroup of total dimension $d$, $e\in G$ denoting the identity element, and $\Liealg g$ will be the Lie superalgebra of $G$, seen as the set of left-invariant vector fields on $G$ (which thus is isomorphic to the tangent space at $e\in G$). When needed, $f_1, \dots, f_d$ will be a homogeneous basis of $\Liealg g$.

The flow of a smooth vector field will always be denoted by the greek letter $\Phi$; a subscript will indicate for which vector field it is the flow and a superscript will indicate a component with respect to some direct product structure.

\end{definition}

\begin{definition}[vectorfieldsonGtimesM]{{Definitions}}
$\bullet$
For any $X\in \Liealg g$ we denote by $X^r$ the right-invariant vector field on $G$ whose value at $e$ is $X\in \Liealg g\cong T_eG$.

$\bullet$
An \stresd{action domain} is an open subset $D\subset G\times M$ satisfying the two conditions
\begin{enumerate}
\item
$\{e\}\times M\subset D$ and

\item
for all $m\in M$ the set $D_m\subset G$ defined by
$$
g\in D_m \qquad\Longleftrightarrow\qquad (g,m)\in D
$$
is connected (and thus a connected open neighbourhood of $e\in G$).

\end{enumerate}

$\bullet$
Let $M$ be a supermanifold and $G$ a connected Lie supergroup. A \stresd{smooth local left-action of $G$ on $M$}  (or simply \stresd{a local action of $G$ on $M$}) is a smooth map $\Psi:D\to M$ defined on an action domain $D\subset G\times M$ satisfying the two conditions
\begin{enumerate}
\item
for all $m\in M$: $\Psi(e,m)=m$ and

\item
if $(g,m)$ and $(hg,m)$ belong to $D$ and if $\bigl(h,\Psi(g,m)\bigr)$ belongs to $D$, then we have the equality
$$
\Psi\bigl(h,\Psi(g,m)\bigr) = \Psi(hg,m)
\mapob.
$$

\end{enumerate}
If $D=G\times M$, the action is called \stresd{global}.

$\bullet$
Let $c_{ij}^k \in \RR$ be the structure constants of $\Liealg g$ associated to the basis $f_1, \dots, f_d\,$: $[f_i,f_j] = \sum_{k=1}^d c_{ij}^k \,f_k$. 
An even (left-)linear map $\rho$ from $\Liealg g$ to the space $\Gamma(TM)$ of sections of the tangent bundle $TM$ of a supermanifold $M$ is called a \stresd{smooth representation of $\Liealg g$ on $M$} or an \stresd{infinitesimal action of $\Liealg g$ on $M$} (in the context of this paper we prefer the latter) if the vector fields $\rho(f_i)$ are smooth (necessarily of the same parity as $f_i$ because $\rho$ is assumed to be even) and satisfy the commutation relations $[\rho(f_i),\rho(f_j)] = \sum_{k=1}^d c_{ij}^k\,\rho(f_k)$. An equivalent condition would be to require that $\rho(X)$ is smooth for all $X\in \body\Liealg g$ (the $X\in \Liealg g$ having real coordinates) and such that $[\rho(X),\rho(Y)] = \rho([X,Y])$ for all $X,Y\in \body \Liealg g$.

$\bullet$
If $\rho$ is an infinitesimal action of $\Liealg g$ on $M$, then the three even smooth vector fields $Z_R$, $Z_M$ and $Z_A$ on $\Liealg g_0 \times G$, $\Liealg g_0\times M$ and $\Liealg g_0\times G\times M$ respectively are defined by
\begin{gather*}
Z_R\caprestricted_{(X,g)}
=
0\caprestricted_X + X^r\caprestricted_g
\quad,\quad
Z_M\caprestricted_{(X,m)}
=
0\caprestricted_X - \rho(X)\caprestricted_m
\\
\qquad\text{and}\qquad
Z_A\caprestricted_{(X,g,m)}
=
0\caprestricted_X + X^r\caprestricted_g - \rho(X)\caprestricted_m
\mapob.
\end{gather*}
Their flows are denoted by $\Phi_R$, $\Phi_M$ and $\Phi_A$ respectively. 
It is immediate that $\Phi_R$ is defined on the whole of $\CA_0 \times \Liealg g_0\times G$ (\ie, $Z_R$ is complete) and that it is given in terms of the exponential map by
$$
\Phi_R(t,X,g) = \bigl(X, \exp(tX)g \bigr)
\mapob.
$$
(Nota Bene: it might be better to say that the exponential map is defined by this equation in terms of the flow of $Z_R$, but that is of less importance here.) 
On the other hand, the flow $\Phi_M$ is defined on an open subset $W_M\subset \CA_0 \times \Liealg g_0 \times M$ satisfying the condition that for each $(X,m)\in \Liealg g_0\times M$ the set
$$
I_{(X,m)} = \{\,t\in \CA_0 \mid (t,X,m)\in W_M\,\}
=
\CA_0\times \{(X,m)\}\cap W_M \subset \CA_0
$$
or equivalently
$$
I_{(X,m)} \times \{(X,m)\}
=
\bigl(\CA_0\times \{(X,m)\}\bigr)\cap W_M \subset \CA_0
$$
is connected and contains $0$ (the domain of definition of the maximal integral curve passing at $t=0$ through $(X,m)$, but in super differential geometry not all these separate curves are differentiable).

\end{definition}

\begin{definition}{{Remarks}}
$\bullet$
The data of a (left-)linear map $\rho$ from $\Liealg g$ to the space of sections of the tangent bundle $TM$ such that the vector fields $\rho(f_i)$ are smooth of the same parity as $f_i$ is completely encoded in the single even smooth vector field $Z_M$:
$$
Z_M\caprestricted_{(X,m)} 
=
0_X - \sum_{i=1}^d X^i \cdot \rho(f_i)\caprestricted_m
\mapob,
$$
where $X= \sum_{i=1}^d X^i \cdot f_i$ is the decomposition of $X\in \Liealg g$ with respect to the given basis. Since the $\rho(f_i)$ are smooth, the last equality immediately shows that $Z_M$ is indeed smooth on the product $\Liealg g_0\times M$. With a slight abuse of notation, one could say that we have $Z_M\caprestricted_{(X,m)} = -\rho(X)\caprestricted_m$, \ie, $Z_M$ at $(X,m)$ equals $-\rho(X)$ at $m$.

$\bullet$
We can extend the idea of the global vector field $Z_M$ to encode the commutation property. To do so, we consider the two global smooth vector fields $Z_1$ and $Z_2$ on $\Liealg g\times \Liealg g\times M$ defined by
$$
Z_1\caprestricted_{(X,Y,m)} = 0\caprestricted_X + 0\caprestricted_Y + \rho(X)\caprestricted_m
\quad\text{,}\quad
Z_2\caprestricted_{(X,Y,m)} = 0\caprestricted_X + 0\caprestricted_Y + \rho(Y)\caprestricted_m
\mapob.
$$
With the same abuse of notation as above, we thus can write $Z_1\caprestricted_{(X,Y,m)} = \rho(X)\caprestricted_m$ and $Z_2\caprestricted_{(X,Y,m)} = \rho(Y)\caprestricted_m$. The commutation condition can then be stated as saying that we should have the equality $[Z_1,Z_2]\caprestricted_{(X,Y,m)} = \rho([X,Y])\caprestricted_m$. The passage via the smooth vector fields $Z_1$ and $Z_2$ avoids the problem that the individual vector fields $\rho(X)$ need not be smooth, and thus that their commutator is not defined.

$\bullet$
The domain of definition of the flow of a vector field is a particular case of an action domain for the group $G=\CA_0$ and the flow itself is a particular case of a local action of $\CA_0$ on the supermanifold.

\end{definition}


\begin{proclaim}[algebrarepassociatedtolocalaction]{Lemma}
Let $\Psi:D\to M$ be a local action of a Lie supergroup $G$ on a supermanifold $M$ and let $\Liealg g$ be the Lie superalgebra of $G$. 
Then the map $\rho:\Liealg g\to \Gamma(TM)$ defined by
$$
\rho(X)\caprestricted_m = -\contrf{X^r\caprestricted_e+0\caprestricted_m}{T\Psi}
$$
is an infinitesimal action of $\Liealg g$ on $M$.

\end{proclaim}

\begin{preuve}
That $\rho$ is left-linear and even is immediate because $T\Psi$ is.
If $X$ belongs to $\body \Liealg g$, the right-invariant vector field $X^r$ on $G$ is smooth.
It follows immediately that the vector field $\Xb\caprestricted_{(g,m)} = X^r\caprestricted_g+0_m$ on $G\times M$ is smooth, and thus its restriction to $D\subset G\times M$ is smooth. As $T\Psi$ is smooth, the composite map
\begin{multline*}
m\mapsto (e,m)\mapsto \Xb\caprestricted_{(e,m)} 
\\
\mapsto -\contrf{\Xb\caprestricted_{(e,m)}}{T\Psi} = -\contrf{X^r\caprestricted_e+0\caprestricted_m}{T\Psi} = \rho(X)\caprestricted_m
\end{multline*}
is smooth.

We now claim that we have the stronger equality
\begin{equation}\label{strongerresultlocalactionrho}
\contrf{\Xb\caprestricted_{(g,m)}}{T\Psi} = -\rho(X)\caprestricted_{\Psi(g,m)}
\end{equation}
for all $(g,m)\in D$, not only for all $(e,m)$.
To prove this, we need the fact that $\Psi$ is a local action. 
If we forget for the moment that it is only local, we can state the action property as the equality $\Psi\bigl(h, \Psi(g,m)\bigr) = \Psi(hg,m)$, both sides of which can be seen as a composition of two smooth maps:
\begin{equation}\label{actionequalitytwomaps}
\begin{aligned}
(h,g,m) 
&
\mapsto \bigl(h,\Psi(g,m)\bigr)\mapsto \Psi\bigl(h, \Psi(g,m)\bigr) 
\quad\text{and}
\\
(h,g,m)
&
\mapsto(hg,m)\mapsto\Psi(hg,m)
\mapob.
\end{aligned}
\end{equation}
When we apply the tangent map of these two smooth maps to the tangent vector $X^r\caprestricted_e+0\caprestricted_g+0\caprestricted_m$, we get (for $\Psi\bigl(h, \Psi(g,m)\bigr)$)
$$
X^r\caprestricted_e+0\caprestricted_g+0\caprestricted_m
\mapsto
X^r\caprestricted_e + 0\caprestricted_{\Psi(g,m)}
\mapsto -\rho(X)\caprestricted_{\Psi(g,m)}
$$
and (for $\Psi(hg,m)$)
$$
X^r\caprestricted_e+0\caprestricted_g+0\caprestricted_m
\mapsto
X^r\caprestricted_g + 0\caprestricted_m = \Xb\caprestricted_{(g,m)}
\mapsto
\contrf{\Xb\caprestricted_{(g,m)}}{T\Psi}
\mapob,
$$
proving our claim.

However, as $\Psi$ is not (necessarily) defined on the whole of $G\times M$, we have to be careful with the domains of definition of these two maps.
We thus choose $(g_o,m_o)\in D$ arbitrarily.
As $D$ is an action domain, there exists an open neighbourhood $V_1$ of $e\in G$ and $U_1$ of $\Psi(g_o,m_o)$ such that $V_1\times U_1\subset D$.
We now recall that right-translation $R_g$, $h\mapsto hg$ is always a homeomorphism of $G$ (it is a diffeomorphism when $g$ belongs to $\body G$).
It follows that $R_g(V_1) = V_1{\cdot} g$ is an open neighbourhood of $g$.
As $\Psi$ is smooth and thus in particular continuous, there exists an open neighbourhood $W\subset D$ of $(g_o,m_o)$ such that $\Psi(W)\subset U_1$.
By taking a smaller $V_1$ if necessary, we may suppose that there exists an open neighbourhood $U_2$ of $m_o$ such that $V_1{\cdot} g_o\times U_2\subset W$.
Finally, as multiplication is smooth, there exists a neighbourhood $V_2$ of $e\in G$ such that $V_2\cdot V_2\subset V_1$, and thus in particular $V_2\subset V_1$.

With these preparations we claim that both maps in \recalf{actionequalitytwomaps} are defined (and thus smooth) on $V_2\times (V_2{\cdot} g_o)\times U_2$.
We thus take $(h,g,m)\in V_2\times (V_2{\cdot} g_o)\times U_2$ and we investigate both maps.
For the first map we have $(g,m)\in (V_2{\cdot} g_o)\times U_2\subset W$ and thus (by definition of $U_2$ and $W$)
\begin{multline*}
(g,m)\in (V_2{\cdot} g_o)\times U_2
\quad\Longrightarrow\quad
(g,m)\in  W
\quad\Longrightarrow\quad
\Psi(g,m)\in U_1
\mapob.
\end{multline*}
But then $\bigl(h,\Psi(g,m)\bigr)\in V_2\times U_1\subset V_1\times U_1\subset D$, and thus $\Psi$ applies.
For the second map we have $(h,g)\in V_2\times (V_2{\cdot} g_o)$ and thus $hg\in V_2{\cdot} V_2{\cdot} g_o\subset V_1{\cdot} g_o$. But then $(hg, m) \in (V_1{\cdot} g_o)\times U_2\subset D$, proving (again) that $\Psi$ applies. This finishes the proof of the stronger equality, as we have shown that it is valid in an open neighbourhood of an arbitrary point $(g_o,m_o)\in D$.

But this stronger equality says that the vector fields $\Xb$ on $D$ and $-\rho(X)$ on $M$ are related by $T\Psi$. 
As related smooth vector fields have related commutators, it follows immediately that for $X,Y\in \body \Liealg g$ the (smooth) vector fields $[\Xb,\,\overline Y \,]$ and $[\,-\rho(X),-\rho(Y)\,] = [\,\rho(X),\rho(Y)\,]$ are related by $T\Psi$. 
But $[\Xb,\,\overline Y \,] = -\overline{[X,Y]}$ (the commutator of right-invariant vector fields is the opposite of the commutator of the corresponding left-invariant vector fields), proving that we have the equality $\rho([X,Y]) = [\rho(X),\rho(Y)]$ for $X,Y\in \body \Liealg g$ (and thus in particular for the basis elements).
\end{preuve}

\begin{definition}{Definition}
Let $\Psi:D\to M$ be a local action of $G$ on $M$ and let $\rho:\Liealg g\to \Gamma(TM)$ be the associated infinitesimal action of $\Liealg g$ on $M$ according to \recalt{algebrarepassociatedtolocalaction}. We then say that the local action $\Psi$ \stresd{integrates} the infinitesimal action $\rho$.

\end{definition}

\begin{proclaim}[linkflowsZAZRandZM]{{Lemma}}
Using ingredients and notation as in \recalt{vectorfieldsonGtimesM}, we have the following properties.
\begin{enumerate}
\item\label{flowofZMinM}
There exists a smooth function $\Phi_M^M:W_M \to M$ such that the flow $\Phi_M$ is given by
$$
\Phi_M(t,X,m)
=
\bigl(X,\Phi_M^M(t,X,m)\bigr)
\mapob.
$$

\item\label{atzerodefinedforallt}
$\CA_0\times \{0\} \times M\subset W_M$ and for all $(t,m)\in \CA_0\times M$ we have $\Phi_M^M(t,0,m) = m$.

\item\label{flowZAintermsofflowZW}
The domain of definition $W_A\subset \CA_0\times \Liealg g_0\times G\times M$ of the flow $\Phi_A$ is given by
\begin{equation}\label{expressionforWsubA}
W_A
=
\{\,(t,X,g,m)\mid (t,X,m)\in W_M\,\}
\end{equation}
and its flow $\Phi_A$ by
\begin{equation}\label{expressionforPhisubA}
\Phi_A(t,X,g,m)
=
\bigl(X,\exp(tX)g, \Phi_M^M(t,X,m) \bigr)
\mapob.
\end{equation}
\end{enumerate}
\end{proclaim}

\langbewijs
\begin{preuve}
$\bullet$ (\ref{flowofZMinM}):
Since $\Phi_M$ is smooth, the composition with the projections onto either $\Liealg g_0$ of $M$ is smooth. As the vector field $Z_M$ is zero in the direction of $X$, the $X$-component of $\Phi_M$ must be constant, proving (\ref{flowofZMinM}).

\smallskip

$\bullet$ (\ref{atzerodefinedforallt}):
If we denote by $j: M\to \Liealg g_0 \times M$ the injection $j(m) = (0,m)$, it is immediate that $Tj$ intertwines the zero vector field on $M$ and $Z_M$ on $\Liealg g_0\times M$. 
It follows that $j$ intertwines their flows. 
But the flow of the zero vector field is defined for all time and is given by the identity: $\Phi_0(t,m) = m$. 
Hence 
$$
\Phi_M(t,0,m) = \Phi_M\bigl(t,j(m)\bigr) = j\bigl(\Phi_0(t,m)\bigr) = j(m)=(0,m)
$$ 
is also defined for all $t\in \CA_0$, \ie, for all $t$ and all $m$ we have $(t,0,m)\in W_M$ and $\Phi_M^M(t,0,m) = m$.

\smallskip

$\bullet$ (\ref{flowZAintermsofflowZW}):
If we denote by $p_{12}$ and $p_{13}$ the projections $p_{12}:\Liealg g_0\times G\times M \to \Liealg g_0\times G$, $(X,g,m)\mapsto (X,g)$ and $p_{13}:\Liealg g_0\times G\times M \to \Liealg g_0\times M$, $(X,g,m)\mapsto (X,m)$, then it is immediate that their tangent maps intertwine the vector fields:
$$
Tp_{12}(Z_A\caprestricted_{(X,g,m)}) = Z_R\caprestricted_{(X,g)}
\qquad\text{and}\qquad
Tp_{13}(Z_A\caprestricted_{(X,g,m)}) = Z_M\caprestricted_{(X,m)}
\mapob.
$$
It follows that these maps intertwine their flows:
\begin{align*}
p_{12}\bigl(\Phi_A(t,X,g,m)\bigr)
&
=
\Phi_R\bigl(t,p_{12}(X,g,m)\bigr)
\\
&
=
\Phi_R(t,X,g)
=
\bigl(X,\exp(tX)g\bigr)
\\
\noalign{\noindent and}
p_{13}\bigl(\Phi_A(t,X,g,m)\bigr)
&
=
\Phi_M\bigl(t,p_{13}(X,g,m)\bigr)
=
\Phi_M(t,X,m)
\\
&
=
\bigl(X,\Phi_M^M(t,X,m)\bigr)
\mapob.
\end{align*}
As the flow $\Phi_A$ necessarily is of the form
$$
\Phi_A(t,X,g,m) = \bigl( \Phi_{A}^{\Liealg g_0}(t,X,g,m), \Phi_{A}^{G}(t,X,g,m), \Phi_{A}^M(t,X,g,m) \bigr)
$$
for smooth functions $\Phi_{A}^{\Liealg g_0}$, $\Phi_{A}^{G}$, $\Phi_{A}^M$ with values in $\Liealg g_0$, $G$ and $M$ respectively, the expression \recalf{expressionforPhisubA} for $\Phi_A$ follows immediately. But it also shows that the domain of definition $W_A$ of $\Phi_A$ contains at least the set on the right hand side of \recalf{expressionforWsubA}. And as the flow can not be defined for values of $t\in \CA_0$ for which the projection is not defined, $W_A$ must be given by this expression.
\end{preuve}
\eindelangbewijs

\kortbewijs
\begin{preuve}
A direct consequence of the fact that the vector fields $Z_M$ and $Z_R$ are related to the vector field $Z_A$ by the appropriate (canonical) projections and the fact that the flow of $Z_R$ is given by the formula
$$
\Phi_R(t,X,g) = \bigl(X, \exp(tX)g \bigr)
\mapob.
$$
\vskip-\baselineskip
\end{preuve}
\eindekortbewijs

\begin{definition}[nonsupercollectionofflows]{{Remark}}
If we look at a submanifold $\{X\} \times M\subset \Liealg g_0\times M$, then the restriction of the vector field $Z_M$ is tangent to this submanifold and equals the vector field $-\rho(X)$. It follows that we can interpret the map $(t,m) \mapsto \Phi_M(t,X,m)$ as the flow of this vector field. In the non-graded case this means that we can group the flows of all these vector fields together to form the flow of $Z_M$, the only bonus of looking at $Z_M$ being that we automatically have a smooth dependence on $X$. On the other hand, in the super case, not all subsets $\{X\} \times M$ are genuine submanifolds and not all $\rho(X)$ are smooth vector fields. And thus in the super case, the passage via $Z_M$ is obligatory. Analogous remarks hold for the vector fields $Z_A$ and $Z_R$.

\end{definition}

\begin{proclaim}[completenessofZA]{{Lemma}}
The vector field $Z_A$ on $\Liealg g_0\times G\times M$ is complete if and only if the vector field $Z_M$ on $\Liealg g_0\times M$ is complete, which is the case if and only if all vector fields $\rho(X)$ with $X\in \body \Liealg g_0$ are complete on $M$, which is the case if and only if all vector fields $\rho(X)$ with $X\in \body \Liealg g_0$ are complete on $\body M$.

\end{proclaim}

\langbewijs
\begin{preuve}
That $Z_A$ is complete if and only if $Z_M$ is complete is a direct consequence of \recaltt{flowZAintermsofflowZW}{linkflowsZAZRandZM}. 
On the other hand, the domain of definition of the flow of a smooth vector field on a supermanifold is completely determined by the domain of definition of the flow of the body of the vector field on the body of the supermanifold \cite[V.4.10]{Tu04}. As we have
$$
\body Z_M\caprestricted_{(X,m)}
=
0\caprestricted_{\body X} - \rho(\body X)\caprestricted_{\body m}
\mapob,
$$
and as these objects live in ordinary differential geometry, we can apply the argument of \recalt{nonsupercollectionofflows} and conclude that the flow of $\body Z_M$ consists of the family of flows for the separate (smooth!) vector fields $\rho(\body X)$. It follows that $Z_M$ is complete if and only if all vector fields $\rho(X)$, $X\in \body \Liealg g_0$ are complete on $\body M$. Reversing the argument on the domain of a flow, this will be the case if and only if all (smooth) vector fields $\rho(X)$, $X\in \Liealg g_0$ are complete on $M$.
\end{preuve}
\eindelangbewijs

\kortbewijs
\begin{preuve}
That $Z_A$ is complete if and only if $Z_M$ is complete is a direct consequence of \recaltt{flowZAintermsofflowZW}{linkflowsZAZRandZM}. 
The other statement is a straightforward consequence of the argument given in \recalt{nonsupercollectionofflows} concerning grouping flows and the fact that for supermanifolds the domain of the flow of an even vector field is completely determined by the domain of its body part, which is a vector field in ordinary differential geometry.
\end{preuve}
\eindekortbewijs

\begin{proclaim}[linkbetweenPsiandPhiM]{Lemma}
Let $\Psi:D\to M$ be a local action of $G$ on $M$, let $W_D\subset \CA_0\times \Liealg g_0\times M$ be defined as 
$$
W_D=\bigl\{\,(t,X,m) \mid \bigl(\exp(tX),m\bigr)\in D\,\bigr\}
$$
and let $W_\Psi\subset W_D$ be the largest action domain (for the abelian group $\CA_0$) contained in $W_D$.
Then we have the two properties
\begin{enumerate}
\item\label{transferofactiondomain}
$\CA_0\times \{0\}\times M\subset W_\Psi\subset W_M$ and

\item\label{linkPsiandPhiMM}
for all $(t,X,m)\in W_\Psi$ we have $\Psi\bigl(\exp(tX),m\bigr) = \Phi_M^M(t,X,m)$.

\end{enumerate}

\end{proclaim}

\begin{preuve}
As $D$ is open, it follows immediately (by continuity of $\exp$ among others) that $W_D$ is open in $\CA_0\times \Liealg g_0\times M$. Moreover, as $D$ is an action domain, we have $\{0\}\times \Liealg g_0\times M\subset W_D$.
This implies immediately that there exists an action domain (with respect to the abelian group $\CA_0$) contained in $W_D$, hence the largest exists too.

$\bullet$ (\ref{transferofactiondomain}):
For any $m\in M$ we have $(0,0,m)\in W_D$ which is open. Hence there exist open neighbourhoods $I\subset \CA_0$ of $0\in \CA_0$, $V\subset \Liealg g_0$ of $0\in \Liealg g_0$ and $U\subset M$ of $m$ such that $I\times V\times U\subset W_D$.
By definition of the DeWitt topology (on $\CA_0$), there exists $\varepsilon>0$ such that we have
$$
\vert \body t\vert<\varepsilon
\qquad\Longrightarrow\qquad
t\in I
\mapob.
$$
Now let $t_o\in \CA_0$ be arbitrary. 
In order to show $(t_o,0,m)\in W_\Psi$, we define $\delta=\max(1,\vert\body t_o\vert)>0$ (we need the max with $1$ to avoid dividing by $0$ lateron) and  the set $V'\subset \Liealg g_0$ by
$$
X\in V' 
\qquad\Longleftrightarrow\qquad
\frac{2\delta X}{\varepsilon}\in V
\mapob.
$$
It follows immediately that $V'$ is an open neighbourhood of $0\in \Liealg g_0$.
Moreover, for all $t\in \CA_0$ verifying $\vert \body t\vert <2\delta$ and all $X\in V'$ we have
$$
tX = \frac{\varepsilon t}{2\delta} \cdot \frac{2\delta X}{\varepsilon} 
\qquad\text{and}\qquad
\frac{\varepsilon t}{2\delta} \in I
\quad,\quad 
\frac{2\delta X}{\varepsilon} \in V
\mapob.
$$
It follows that we have the inclusion
$$
\{\,\vert \body t\vert <2\delta\,\} \times V'\times U\subset W_D
\mapob.
$$
But then the set
$$
W_\Psi \cup \{\,\vert \body t\vert <2\delta\,\} \times V'\times U
$$
also is an action domain (for the group $\CA_0$), because $\{\,\vert \body t\vert <2\delta\,\}$ is connected and $V'$ and $U$ open. Hence by maximality of $W_\Psi$ we must have the inclusion
$$
\{\,\vert \body t\vert <2\delta\,\} \times V'\times U\subset W_\Psi
\mapob.
$$
As $\vert \body t_o\vert \le \delta<2\delta$, we thus have shown that $(t_o,0,m)\in W_\Psi$ as wanted.

$\bullet$ (\ref{linkPsiandPhiMM}):
We define the (auxiliary) map $\Phi:W_\Psi\to \Liealg g_0\times M$ by
$$
\Phi(t,X,m) = \bigl(X,\Psi(\exp(tX),m)\bigr)
$$
and we note that we have the initial condition $\Phi(0,X,m) = (X,m)$ (because $\Psi(e,m)=m$).
We now claim that this $\Phi$ also satisfies the differential equation
$$
\contrf{\partial_t\caprestricted_{(t,X,m)}}{T\Phi} = Z_M\caprestricted_{\Phi(t,X,m)}
\mapob.
$$
We start with the observation that the $M$-component of $\Phi$ is the composition of several maps:
$$
(t,X,m)\mapsto (tX,m)\mapsto \bigl(\exp(tX),m\bigr) \mapsto \Psi\bigl(\exp(tX),m\bigr)
\mapob.
$$
We thus follow what happens to the tangent vector $\partial_t\caprestricted_{(t,X,m)}$:
\begin{multline*}
\qquad
\partial_t\caprestricted_{(t,X,m)}
\mapsto
X\caprestricted_{(tX,m)}
\mapsto
X^r\caprestricted_{(\exp(tX),m)}
\equiv
X^r\caprestricted_{\exp(tX)} + 0\caprestricted_{m}
\\
\noalign{\vskip1\jot}
\text{and}\quad
X^r\caprestricted_{\exp(tX)} + 0\caprestricted_{m}
\mapsto
-\rho(X)\caprestricted_{\Psi(\exp(tX),m)}
\mapob,\qquad
\end{multline*}
where the last line follows from \recalf{strongerresultlocalactionrho}, a result obtained in the proof of \recalt{algebrarepassociatedtolocalaction}.
Adding the zero tangent vector in the direction of $X$ immediately proves our claim.
Uniqueness of the flow of a vector field then tells us that we must have $W_\Psi\subset W_M$ and that $\Phi$ and $\Phi_M$ coincide on $W_\Psi$.
\end{preuve}

\section{The foliation and uniqueness}

\begin{definition}[foliationonGtimesM]{Definition}
Let $\rho$ be an infinitesimal action of $\Liealg g$ on $M$.
Then $\foliation$ is the distribution (subbundle) $\foliation\subset T(G\times M)$ spanned by the smooth vector fields $f_i^r-\rho(f_i)$. More precisely:
\begin{align*}
\foliation\caprestricted_{(g,m)}
&
=
\Bigl\{\ 
\sum_{i=1}^d X_i \cdot\bigl(\, f_i^r\caprestricted_g - \rho(f_i)\caprestricted_m \, \bigr)
\mid
X_i\in \CA
\ \Bigr\}
\\
&
=
\bigl\{\, X^r\caprestricted_g - 
\rho(X)\caprestricted_m
\mid
X\in \Liealg g
\, \bigr\}
\mapob.
\end{align*}
Denoting by $p:G\times M\to G$ the canonical projection $p(g,m)=g$, it is immediate that the tangent map $T_{(g,m)}p : T_{(g,m)}(G\times M) \to T_gG$ is an isomorphism from $\foliation\caprestricted_{(g,m)}$ to $T_gG$.

\end{definition}

\begin{proclaim}{{Lemma}}
The distribution $\foliation\subset T(G\times M)$ is involutive, \ie, $\foliation$ is a foliation on $G\times M$.

\end{proclaim}

\langbewijs
\begin{preuve}
For the right-invariant vector fields $f_i^r$ on $G$ we have the equalities $[f_i^r, f_j^r] = -\sum_{k=1}^d c_{ij}^k \, f_k^r$. It then follows immediately from the fact that $\rho$ is an infinitesimal action and the use of the minus sign in the definition of $\foliation$ that it is involutive.
\end{preuve}
\eindelangbewijs

\begin{definition}{{Terminology}}
We will talk about \quote{leaf-open} subsets and \quote{leaf-con\-ti\-nuous} maps whenever we equip $G\times M$ with the leaf topology $\leaftopol$  \recalt{definitionofleaftopology} associated to the foliation $\foliation$ defined in \recalt{foliationonGtimesM} and when this space appears alone or in a direct product.

\end{definition}

\begin{proclaim}[locadaptedcoordinates]{Lemma}
Let $(g,m)\in G\times M$ be arbitrary and let $E$ be a graded vector space of the same dimension as $M$. Then there exists an open neighbourhood $\Ut\subset G\times M$ of $(g,m)$, an open neighbourhood $V\subset G$ of $g$, an open subset $O\subset E_0$ and a diffeomorphism $\varphi : \Ut\to V\times O$ with the following properties.
\begin{enumerate}
\item
$V$ is a connected coordinate neighbourhood with (local) coordinates $(y_1, \dots, y_d)$.

\item
$(y_1, \dots, y_d, x_{d+1}, \dots, x_{d+n})$ is a (local) coordinate system on $\Ut$ (via $\varphi$) adapted to the foliation $\foliation$, where $x_{d+1}, \dots, x_{d+n}$ are (global) coordinates on $E_0$.

\item
On $\Ut$ we have the equality $p=p_1\scirc \varphi$, where $p_1:V\times O\to V$ is the canonical projection.
$$
\begin{matrix}
& \Ut 
\\
\rlap{\kern1em\llap{\raise1.5ex\hbox{$\scriptstyle\varphi$}\kern0pt}$\swarrow$} && 
\llap{$\searrow$\rlap{\raise1.5ex\hbox{$\scriptstyle p$}}\kern0.5em}
\\
\noalign{\vskip2\jot}
V\times O\quad & \underset{p_1}{{-}\!{-}\!{\to}} &\quad V
\end{matrix}
$$

\end{enumerate}

\end{proclaim}

\begin{preuve}
A sloppy but rather short proof would be the following argument.
Let $U$ be an open neighbourhood of $(g,m)$ with local coordinates $x_1, \dots$, $x_{d+n}$ adapted to the foliation $\foliation$ and let $y_1, \dots, y_d$ be local coordinates in an open neighbourhood of $g\in G$. 
As the foliation is spanned by the tangent vectors $\partial_{x_i}$, $i\le d$ and as the tangent map of the projection $p:G\times M \to G$ is an isomorphism from $\foliation_{(g,m)}$ to $T_gG$, it follows that the $d\times d$ matrix
$$
\Bigl(\fracp{y_i}{x_j}(g,m)\Bigr)_{i,j=1}^d
$$
is invertible. 
By the inverse function theorem it follows that there exists an open neighbourhood $\Ut\subset U$ of $(g,m)$ such that $y_1, \dots, y_d$, $x_{d+1}, \dots, x_{d+n}$ forms a local system of coordinates (on $\Ut\subset G\times M$) adapted to the foliation. 
Without loss of generality (by taking smaller $\Ut$ and $V$ if necessary) we may assume that $\Ut$ is of the form $V\times O$ with $(x_{d+1}, \dots, x_{d+n})$ local coordinates on $O$, where $O$ is an open (coordinate) subset of a the even part of a graded vector space of the same graded dimension as $M$.

The sloppiness of this argument lies in the fact that we are mixing subsets of the manifold with the coordinate charts (which are images by maps from subsets of the manifold to sets in a (graded) vector space).
If we want to be a bit more precise, we start by introducing a second graded vector space $F$, this time of the same graded dimension as $G$.
Then any coordinate system on $G\times M$ is a bijective map $\varphit:\Ut\to \Ot\subset (F\times E)_0$ from an open set $\Ut\subset G\times M$ to an open set $\Ot$ in the even part of the graded vector space $F\times E$.
Writing the image with respect to a homogeneous basis $e_1, \dots, e_d, e_{d+1}, \dots, e_{d+n}$ of $F\times E$ (with $e_1, \dots, e_d$ a basis of $F$ and $e_{d+1}, \dots, e_{d+n}$ one of $E$) as
$$
\varphit(g,m) = \sum_{i=1}^{d+n} \varphit_i(g,m)\cdot e_i
\equiv
\sum_{i=1}^{d+n} x_i(g,m)\cdot e_i
$$
then provides us with the coordinate functions $x_i$ on the coordinate chart $\Ut$.
Similarly, a system of local coordinates on $G$ is a bijective map $\chi:V\to \Vt\subset F_0$ from an open subset $V\subset G$ to an open subset $\Vt\subset$ in the even part of $F_0$. 
And then the local coordinates $y_i$ are defined by
$$
\chi(g) = \sum_{i=1}^d \chi_i(g)\cdot e_i \equiv \sum_{i=1}^d y_i(g)\cdot e_i
\mapob.
$$
Replacing $\Ut$ by $\Ut\cap p\mo(V)$ and $\Ot$ by $\varphit\bigl(p\mo(V)\bigr)\cap \Ot=\varphit\bigl( \Ut\cap p\mo(V)\bigr)$, we may assume that $p(\Ut)\subset V$ and thus that the map
$$
\chi\scirc p\scirc \varphit\mo : \Ot\to \Vt
$$
is well defined.

The fact that the coordinates $x_1, \dots, x_{d+n}$ are adapted to the foliation and that the foliation projects bijectively onto the tangent space of $G$ then translates as the fact that the $d\times d$ matrix
$$
\Bigl(\fracp{(\chi_i\scirc p\scirc \varphi\mo)}{x_j}\bigl(\varphi(g,m)\bigr)\Bigr)_{i,j=1}^d
$$
is invertible.
Introducing the map $\rho:\Ot\to \Vt\times E_0$ by
\begin{multline*}
\quad
\rho(x_1, \dots, x_{d+n})
=
\bigl(\,(\chi_1\scirc p\scirc \varphit\mo)(x_1, \dots, x_{d+n}), \dots, 
\\
(\chi_d\scirc p\scirc \varphit\mo)(x_1, \dots, x_{d+n}), 
x_{d+1}, \dots, x_{d+n}\,\bigr)
\quad
\end{multline*}
it follows that the Jacobian matrix of $\rho$ at $\varphit(g,m)$ is invertible. 
By the inverse function theorem we deduce the existence of open subsets $\Ot'\subset \Ot$ and $\Ot''\subset \Vt\times E_0$ such that 
$$
\rho:\Ot'\to \Ot''
$$
is a diffeomorphism.
Taking a smaller $\Ot''$, we may assume that this open subset is of the form $\Ot''=\Vt''\times O$ for some open subset $\Vt''\subset \Vt$ and an open subset $O\subset E_0$.
\begin{equation*}
\begin{matrix}
\Ut & \overset{\varphit}{{-}\!{-}\!{-}\!{\to}} & \Ot & & \Ot & \supset & \Ot'
\\
\llap{$\scriptstyle p$}\downarrow &&\downarrow\rlap{$\scriptstyle \chi\scirc p \scirc \varphit\mo$} &\qquad\qquad & \llap{$\scriptstyle \rho$}\downarrow && \downarrow \rlap{$\scriptstyle \rho$}
\\
V & \overset{\chi}{{-}\!{-}\!{-}\!{\to}} & \Vt && \Vt\times E_0 & \supset & \Ot''=\Vt''\times O
\end{matrix}
\mapob.
\end{equation*}

With these preparations we replace $\Ut$ by $\varphit\mo(\Ot')$ and $V$ by $\chi\mo(\Vt'')$ and we define (with the new $\Ut$ and $V$) $\varphi : \Ut\to V\times O$ by
$$
\varphi = (\chi\mo\times id)\scirc \rho \scirc \varphit
\mapob.
$$
It then follows immediately that $\varphi$ is a local coordinate system and that we have
\begin{align*}
p_1\scirc \varphi
&
=
p_1\scirc (\chi\mo\times id)\scirc \rho \scirc \varphit
\\
&
=
\chi\mo \scirc (p_1\scirc \rho)\scirc \varphit
=
\chi\mo\scirc (\chi\scirc p\scirc \varphit\mo) \scirc \varphit
=
p
\mapob.
\end{align*}
Moreover, as $\varphit$ is adapted to the foliation and as the maps $\chi\mo\times id$ and $\rho$ send the subspace $\mathrm{span}(\partial_{x_1}, \dots, \partial_{x_d}) \subset T_{(x_1, \dots, x_{d+n})}(F_0\times E_0)$ to itself, the map $\varphi$ is again a coordinate system adapted to the foliation.
\end{preuve}

\begin{proclaim}[projectionlocalleafhomeo]{{Proposition}}
The map $p:G\times M\to G$ is a local leaf-homeo\-mor\-phism, and thus in particular leaf-continuous.

\end{proclaim}

\langbewijs
\begin{preuve}
What we have to show is that for all $(g,m)\in G\times M$ there exists a leaf-open neighbourhood $W$ of $(g,m)$ and an open neighbourhood $V$ of $g$ such that $p:W\to V$ is a homeomorphism (when $W$ is equipped with the topology induced by the leaf topology). 
So choose $(g,m)\in G\times M$ and let $\Ut$, $V$, $O$ and $\varphi$ be as in \recalt{locadaptedcoordinates}.
It follows immediateley that the \slice{} $W\equiv \Ut_{a_{\kleinlessthan d}}$ (which is a leaf-open neighbourhood of $(g,m)$ for a suitable choice of $a$) is given by
$$
\Ut_{a_{\kleinlessthan d}} = \varphi\mo\bigl( V \times \{(a_{d+1}, \dots, a_{d+n})\} \bigr)
\mapob.
$$
As the topology on $V\times O$ is the product topology, the map
$$
p_1:V\times \{(a_{d+1}, \dots, a_{d+n})\}\to V
$$ 
is a homeomorphism. But $\varphi$ is a diffeomorphism, hence 
$$
p:\Ut_{a_{\kleinlessthan d}} \equiv \varphi\mo\bigl( V \times \{(a_{d+1}, \dots, a_{d+n})\} \bigr) \to V
$$ 
is a homeomorphism when $\Ut_{a_{\kleinlessthan d}}$ is equipped with the topology induced by the topology of $G\times M$. As this is the same topology as the one induced by the leaf-topology \recalt{localleaftopisrelativetop}, it follows that $p:\Ut_{a_{\kleinlessthan d}} \to V$ is a leaf-homeomorphism.
\end{preuve}
\eindelangbewijs

\kortbewijs
\begin{preuve}
By Frobenius' theorem and the inverse function theorem we may assume that (locally) we have a coordinate system $(x_1, \dots, x_{d+n})$ on $G\times M$ adapted to the foliation and such that $(x_1, \dots, x_d)$ are local coordinates on $G$. The result then follows immediately from \recalt{localleaftopisrelativetop}.
\end{preuve}
\eindekortbewijs

\begin{proclaim}[righttranslpreservesleaves]{{Lemma}}
The map $R:(G\times M) \times G\to G\times M$ defined as $R\bigl((g,m),h\bigr) = (gh,m)$ is leaf-continuous. In particular for a fixed $h\in G$ the map $R_h:G\times M\to G\times M$ defined as $R_h(g,m) = (gh,m)$ is a leaf-homeomorphism. 

\end{proclaim}

\langbewijs
\begin{preuve}
Choose $\bigl((g,m),h\bigr)\in (G\times M)\times G$ and a basic leaf-open neighbourhood $V_{b_{\kleinlessthan d}}^c$ of $(gh,m)$, which means that $V\subset G\times M$ is an open neighbourhood of $(gh,m)$ with local coordinates $y_1, \dots, y_{d+n}$ adapted to the foliation (and as before, the superscript ${}^c$ indicates that we take the appropriate connected component). 
As $R$ is smooth, it is in particular continuous (for the original topology), so there exist open neighbourhoods $U\subset G\times M$ and $W\subset G$ of $(g,m)$ and $h$ respectively such that $R(U\times W)\subset V$. 
By taking a smaller $U$ if necessary, we may assume without loss of generality that there exists local coordinates $x_1, \dots, x_{d+n}$ on $U$ adapted to the foliation. 
By taking a smaller $W$ if necessary, we may also assume that $W$ is connected.

Now consider the vector field $X_P$ on $(G\times M)\times G$ defined as
$$
X_P\caprestricted_{((g,m),h)} = f_i^r\caprestricted_g - \rho(f_i)\caprestricted_m + 0\caprestricted_h
\mapob,
$$
where $f_i$ is one of the basis elements of the Lie superalgebra $\Liealg g$.
Then, as $f_i^r$ is right-invariant, it follows immediately that the tangent map of $R$ produces the image
$$
\contrf{X_P\caprestricted_{((g,m),h)}}{TR}
=
f_i^r\caprestricted_{gh} - \rho(f_i)\caprestricted_m
\mapob.
$$
With a slight abuse of notation, this means that $TR$ maps the foliation $\foliation$ (extended to the product $(G\times M)\times G$) to itself. If we now denote by $z_i$ local coordinates on $W\subset G$, then we can write (for any $1\le i\le n+d$):
$$
y_i = R_i(x, z)
\mapob.
$$
But on $U\times W$ the foliation is spanned by $\partial_{x_i}$ for $i\le d$ and on $V$ it is spanned by $\partial_{y_i}$ for $i\le d$. The fact that $R$ maps the foliation to itself thus implies that we must have
$$
\fracp{R_i}{x_j}(x,z) = 0 
\qquad\text{for $i>d$ and $j\le d$.}
$$
And thus $y_i$, $i>d$ is constant on the connected set $U_{a_{\kleinlessthan d}}^c \times W$, which implies that we have the inclusion
$$
R(U_{a_{\kleinlessthan d}}^c \times W) \subset V_{b_{\kleinlessthan d}}^c
$$
as wanted.
\end{preuve}
\eindelangbewijs

\kortbewijs
\begin{preuve}
This is a straightforward consequence of Frobenius' theorem and the fact that the generating vector fields $X^r-\rho(X)$ of the foliation are invariant under right-translation: the vector fields $X^r$ are right-invariant by definition, and the vector fields $\rho(X)$ are not affected.
\end{preuve}
\eindekortbewijs

\begin{proclaim}[flowisleafcontinuous]{{Proposition}}
The flow $\Phi_A$ is leaf-continuous.

\end{proclaim}

\langbewijs
\begin{preuve}
What we have to show is the following. Suppose we are given a point $(t_o,X_o,g_o,m_o)\in W_A\subset  \CA_0\times \Liealg g_0 \times G\times M$, an open set $W_1\subset \Liealg g_0$ and a basic open set $V_{b_{\kleinlessthan d}}^c\subset G\times M$ for the leaf topology (where the superscript ${}^c$ here and in the sequel indicates that we take the appropriate connected component) such that $\Phi_A(t_o,X_o,g_o,m_o)\in W_1\times V_{b_{\kleinlessthan d}}^c$. Then we have to find open neighbourhoods $I_0\subset \CA_0$ of $t_o$, $W_0\subset \Liealg g_0$ of $X_o$ and $(U_o)_{a_{\kleinlessthan d}}^c\subset G\times M$ (a basic open set for the leaf topology) of $(g_o,m_o)$ such that we have the inclusion
$$
\Phi_A(I_0 \times W_0\times (U_0)_{a_{\kleinlessthan d}}^c) \subset W_1 \times V_{b_{\kleinlessthan d}}^c
\mapob.
$$
To prove this, we proceed in two steps. In the first step we show that it is true for all points in $W_A$ with $t$ \quote{sufficiently small.} And in the second step we show that it is true for all points in $W_A$.

For the first step, we start with a local coordinate system $x_1, \dots, x_{d+n}$ on an open set $U\subset G\times M$ adapted to the foliation. 
It follows immediately that the vector field $Z_A$ has the local expression on $\Liealg g_0 \times U$
$$
Z_A\caprestricted_{(X_i, x_j)} = \sum_{i,j=1}^d X_i\cdot F_{ij}(x)\cdot \fracp{}{x_j}\bigrestricted_{(X_i,x_j)}
\mapob,
$$
where $X_i$ denote the coordinates of $X\in \Liealg g_0$ with respect to the fixed basis $f_1, \dots, f_n$ of $\Liealg g$: $X=\sum_{i=1}^d X_i \,f_i$, and where $F_{ij}$ are smooth functions on $G\times M$ (for fixed $i$ they represent the coefficients of the smooth vector field $f_i^r - \rho(f_i)$ on $G\times M$ with respect to the coordinate system $x$).
Local existence and uniqueness of the solutions of ordinary differential equations then tells us that for any $(X_o,g_o,m_o)\in \Liealg g_0\times U$ there exists a connected open neighbourhood $I_o\subset \CA_0$ of $0$, a connected open neighbourhood $W_o\subset \Liealg g_0$ of $X_o$ and an open neighbourhood $U_o\subset U$ of $(g_o,m_o)$ such that there exists a (unique) local flow 
$$
\Phi_A:I_o\times W_o\times U_o\to \Liealg g_0\times U
$$
for the restriction of the vector field $Z_A$ to $\Liealg g_0\times U$. As the coordinates $X_i$ and $x_j$ can be used as well in the source as in the target space of this local flow $\Phi_A$, and because the components of $Z_A$ in the direction of $\partial_{X_i}$ and $\partial_{x_j}$ for $j>d$ are zero, we can write this local flow sloppily as (using that $I_o$ is connected)
$$
\Phi_A(t,X_i,x_j) = \bigl(X_i, x_j(t)\bigr)
\qquad\text{with}\qquad
x_j(t) = x_j
\text{ for $j>d$.}
$$
Now let $(t,X,g,m)\in I_o\times W_o\times U_o$ be arbitrary and let $W_1\times V_{b_{\kleinlessthan d}}^c$ be a basic open neighbourhood of $\Phi_A(t,X,g,m)$ for the leaf topology. As $\Phi_A(t,X,g,m)$ also belongs to $\Liealg g_0\times U$, there exists $a_{d+1}, \dots, a_{d+n}$ such that $\Phi_A(t,X,g,m)\in W_1 \times U_{a_{\kleinlessthan d}}$. 
According to the proof of \recalt{slicesbasistopology}, we thus have 
$$
\Phi_A(t,X,g,m) \in W_1\times (V\cap U)_{a_{\kleinlessthan d}}^c \subset W_1\times V_{b_{\kleinlessthan d}}^c
\mapob.
$$

As $\Phi_A$ is continuous for the usual topologies, there exist connected open neighbourhoods $I'_o\subset I_o$ of $t$, $W'_o\subset W_o$ of $X$ and $U'_o\subset U_o$ of $m$ such that we have 
$$
\Phi_A(I'_o\times W'_o\times U'_o) \subset W_1 \times (U\cap V)
\mapob.
$$ 
As $x_j(t) = x_j$ for $j>d$ we have (using connectedness!)
$$
\Phi_A(I'_o\times W'_o\times (U'_o)_{a_{\kleinlessthan d}}^c) \subset W_1 \times (U\cap V)_{a_{\kleinlessthan d}}^c
 \subset W_1\times V_{b_{\kleinlessthan d}}^c
\mapob,
$$
which shows that $\Phi_A$ is continuous for the leaf topology at $(t,X,g,m)$. 
We thus have shown that for all points $(X_o,g_o,m_o)\in \Liealg g_0\times G\times M$ there exists an open neighbourhood $I_o\times W_o\times U_o$ of $(0,X_o,g_o,m_o)\in W_A$ such that $\Phi_A$ is leaf-continuous for all points in this neighbourhood.

For the second step we fix $(X,g,m)$ and define $C\subset \CA_0$ as 
\begin{multline*}
C
=
\{\,t\in \CA_0 \mid (t,X,g,m)\in W_A 
\\
\text{ and } \Phi_A \text{ leaf-continuous at } (t,X,g,m)\,\}
\mapob.
\end{multline*}
We then note that according to the first step, $C$ contains on open interval containing $0$. So let $t\in \overline C$ such that $(t,X,g,m)\in W_A$. 
We now apply the first step to the point $\Phi_A(t,X,g,m)$ and conclude that there exist open subsets $I_o\subset \CA_0$ containing $0$, $W_o\subset \Liealg g_0$ and $U_o\subset G\times M$ such that $\Phi_A(t,X,g,m) \in W_o\times U_o$ and such that $\Phi_A$ is leaf-continuous at all points in $I_o\times W_o\times U_o\subset W_A$.

As $\Phi_A$ is continuous, $\Phi_A\mo(W_o\times U_o)\subset W_A$ is open, and thus there exists an open interval $I_2\subset \CA_0$ containing $t$ such that we have the inclusion $\Phi_A(I_2 \times \{(X,g,m)\})\subset W_o\times U_o$ and (of course) $I_2 \times \{(X,g,m)\} \subset W_A$. 
As $I_o$ is open and contains $0$, we may assume (by taking a smaller $I_2$ if needed) that we also have the inclusion 
$$
I_2 - I_2 \equiv \{\, s-s' \mid s,s'\in I_2\,\} \subset I_o
\mapob.
$$
Now take $t''\in I_2$ arbitrary and note that, because $t$ is in the closure of $C$, there exists $t'\in C\cap I_2$. And thus in particular $\Phi_A$ is leaf-continuous at $(t',X,g,m)$. 
On the other hand, we have $\Phi_A(t',X,g,m)\in W_o\times U_o$ and $t''-t'\in I_2-I_2\subset I_o$ and thus $\bigl( t''-t',\Phi_A(t',X,g,m)\bigr)\in I_o \times W_o\times U_o \subset W_A$. We thus can apply the group law to obtain the equality
$$
\Phi_A(t'',X,g,m) = \Phi_A\bigl(t''-t', \Phi_A(t',X,g,m) \bigr)
\mapob.
$$
But by the same token we have that $\Phi_A$ is leaf-continuous at $(t',X,g,m)$ and at $\bigl( t''-t',\Phi_A(t',X,g,m)\bigr)$. 
And thus by composition of maps, $\Phi_A$ is leaf-continuous at $(t'',X,g,m)$. 
It follows that $C$ is open and closed in ${\bigl(\CA_0 \times\{(X,g,m)\}\bigr)} \cap W_A$, which is connected. And thus we have the equality $C = \bigl(\CA_0 \times\{(X,g,m)\}\bigr) \cap W_A$. The final conclusion is that $\Phi_A$ is leaf-continuous at all points of $W_A$.
\end{preuve}
\eindelangbewijs

\kortbewijs
\begin{preuve}
The proof proceeds in two steps. In the first step one proves that $\Phi_A$ is leaf-continuous for all points $(t,X,g,m)$ with $t$ sufficiently small, \ie, in an open neighbourhood of $\{0\}\times \Liealg g_0 \times G\times M$ and in the second step that it is leaf-continuous on the whole domain of definition of $\Phi_A$.

For the first step we choose local coordinates $x_1, \dots, x_{d+n}$ on $G\times M$ adapted to the foliation, to which we add coordinates $X_i$ on $\Liealg g_0$. 
Local existence and uniqueness of differential equations then tell us that the flow $\Phi_A$ exists in a neighbourhood of $t=0$, and the form of the vector field $Z_A$ implies that the only coordinates that depend upon $t$ are $x_1, \dots, x_d$. It follows in a straightforward way that $\Phi_A$ is leaf-continuous for all points in an open neighbourhood of $(0,X_o,g_o,m_o)$.

For the second step we argue with integral curves, for which we fix $(X,g,m)$. By the first step, the set $C\subset \CA_0$ of values for $t$ such that $\Phi_A$ is leaf-continuous at $(t,X,g,m)$ contains an open neighbourhood of $t=0$. 
Now let $t\in \overline C$ be such that $(t,X,g,m)\in W_A$. Again by the first step, $\Phi_A$ is leaf-continuous for all points in an open neighbourhood, say $\tilde V$ of $\bigl(0,\Phi_A(t,X,g,m)\bigr)$. 
But $\Phi_A$ is continuous, hence there exists an open neighbourhood, say $\tilde U$ of $(t,X,g,m)$ such that its image by $\Phi_A$ is contained in $\tilde V$. But $C\cap \tilde U$ contains points arbitrary close to $t$. Writing
$$
\Phi_A(t'',X,g,m) = \Phi_A\bigl(t''-t', \Phi_A(t',X,g,m) \bigr)
\mapob,
$$
it follows in a straightforward manner that $\Phi_A$ is leaf-continuous for all points in an open neighbourhood of $t\subset \CA_0$. Hence $C$ is both open and closed, it is non-empty and connected (by defintion of the domain of a flow), hence it is the \quote{full} integral curve through $(X,g,m)$, proving that $\Phi_A$ is leaf-continuous at all points of $W_A$.
\end{preuve}
\eindekortbewijs

\begin{definition}{Convention}
We now fix once and for all an open neighbourhood $\Domexp\subset \Liealg g_0$ of $0\in \Liealg g_0$ such that $\exp:\Domexp\to G$ is a diffeomorphism onto its image $\Imexp=\exp(\Domexp)$, which thus is an open neighbourhood of $e\in G$.

\end{definition}

\begin{definition}[defoflocalleafcomponents]{Definition\slash Construction}
For any $(g,m)\in G\times M$ and any open neighbourhood $V\subset \Imexp$ of $e\in G$ such that
$$
\{1\}\times \exp\mo(V)\times \{(g,m)\}\subset W_A
$$ 
(which is equivalent to the requirement $\{1\}\times \exp\mo(V)\times \{m\}\subset W_M$),
we define the map $\psi_{g,m,V}:V\cdot g\to G\times M$ by
\begin{align*}
\psi_{g,m,V}(h) 
&
= 
(p_{23}\scirc \Phi_A)\bigl(1, \exp\mo(hg\mo), g,m\bigr)
\\
&
\overset{\text{\rlap{\hss\recalt{expressionforPhisubA}\vrule depth3pt width0pt}}}{=}
\ \Bigl(\ \exp\bigl(1\cdot \exp\mo(hg\mo)\bigr)\cdot g\ \ ,\ \  \Phi_M^M(1, \exp\mo(hg\mo),m)\ \Bigr)
\\
&
=
\bigl(\ h\ ,\ \Phi_M^M(1, \exp\mo(hg\mo),m)\ \bigr)
\mapob,
\end{align*}
where $p_{23}:\Liealg g_0\times G\times M\to G\times M$ denotes the canonical projection $p_{23}(X,g,m) = (g,m)$.
Such maps do indeed exist: we know that $(1,0,g,m)\in W_A$ \recaltt{atzerodefinedforallt}{linkflowsZAZRandZM}. As $W_A$ is open, there exists in particular an open neighbourhood $U$ of $0\in \Liealg g_0$ such that $\{1\}\times U\times \{(g,m)\}$ is contained in $W_A$. By taking a smaller $U$ if necessary, we may assume that $U\subset \Domexp$ and then it suffices to take $V=\exp(U)$.

\end{definition}

\begin{proclaim}[localleafcomponents]{Lemma}
Let $g$, $m$, $V$ and $\psi_{g,m,V}$ be as in \recalt{defoflocalleafcomponents}. Then $U_{g,m,V} = \psi_{g,m,V}(V{\cdot} g)$ is leaf-open in $G\times M$ and $\psi_{g,m,V}:V\cdot g\to U_{g,m,V}$ is a leaf-homeomorphism with $p$ (or more precisely $p\restricted_{U_{g,m,V}}$) as its inverse.

\end{proclaim}

\begin{preuve}
The map $h\mapsto hg\mo$ is a homeomorphism from $V\cdot g$ to $V$ (right translation is always a homeomorphism; it is a diffeomorphism when $g$ belongs to $\body G$) and $\exp\mo$ is a diffeomorphism and hence a homeomorphism from $V\subset \Imexp$ onto its image in $\Domexp$.
As we have the product topology on the triple product $\CA_0\times \Liealg g_0\times (G\times M)$ (with the leaf topology on the third factor), the canonical injection $X\mapsto (1,X,g,m)$ is leaf-continuous.
For the same reason the projection $p_{23}:\Liealg g_0\times (G\times M) \to G\times M$ is leaf continuous. And finally the map $\Phi_A$ is leaf continuous \recalt{flowisleafcontinuous}. And thus $\psi_{g,m,V}$ is leaf-continuous as composition of leaf-continuous maps.

As we obviously have the equality
$$
p\scirc \psi_{g,m,V} =id\caprestricted_{V\cdot g}
\mapob,
$$
and as $\psi_{g,m,V}$ and $p$ are leaf-continuous \recalt{projectionlocalleafhomeo}, it follows that $p:U_{g,m,V}\to V\cdot g$ is a leaf-homeomor\-phism with $\psi_{g,m,V}$ as its inverse. 

Now choose $x\in U_{g,m,V}$ and define $h=p(x)$. Since $p$ is locally a leaf-homeo\-mor\-phism \recalt{projectionlocalleafhomeo}, there exists a leaf-open neighbourhood $U'$ of $x$ and an open neighbourhood $V'\subset V\cdot g$ of $h$ such that $p:U'\to V'$ is a leaf-homeomorphism. 
Since $\psi_{g,m,V}$ is leaf-continuous, $V''=V'\cap \psi_{g,m,V}\mo(U')\subset V'$ is an open neighbourhood of $h$. Now $\psi_{g,m,V}(V'')\subset U'\cap U_{g,m,V}$ and on $U_{g,m,V}$ the map $\psi_{g,m,V}$ is the inverse of $p$. 
But $p$ is a leaf-homeomorphism on $U'$ and thus $\psi_{g,m,V}(V'') = (p\restricted_{U'})\mo(V'')$ is leaf-open, which shows that $x$ has a leaf-open neighbourhood contained in $\psi_{g,m,V}(V\cdot g)=U_{g,m,V}$. This shows that $U_{g,m,V}$ is leaf-open.
\end{preuve}

\begin{proclaim}[opentoconnectedisconnectedcomponent]{Lemma}
Let $U\subset G\times M$ be leaf open such that $p:U\to \Ub=p(U)$ is a bijection and $\Ub$ connected. Then $\Ub$ is open in $G$ and $U$ is a leaf-connected component of $p\mo(\Ub)$.

\end{proclaim}

\begin{preuve}
That $\Ub$ is open is an immediate consequence of the fact that $p$ is a local leaf-homeomorphism. And thus $p:U\to \Ub$ is a homeomorphism, implying that $U$ is connected.
It follows that $U$ is an open connected subset of $p\mo(\Ub)$, so in order to prove that it is a leaf-connected component of $p\mo(\Ub)$, it suffices to show that $U$ is leaf-closed in $p\mo(\Ub)$.
To prove that, it suffices to find for any $(g,m)\in p\mo(\Ub)\setminus U$ a leaf-open neighbourhood $U'$ such that $U'\cap U=\emptyset$.

So let $(g,m)\in p\mo(\Ub)\setminus U$ be arbitrary.
As $(g,m)\in p\mo(\Ub)$, it follows that $g\in \Ub$ and that there exists $m'\in M$ such that $(g,m')=(p\restricted_U)\mo(g)\in U$.
Now let $V_1$ and $V_2$ be two open neighbourhoods of $e\in G$ such that the maps $\psi_{g,m,V_1}$ and $\psi_{g,m',V_2}$ are defined.
It follows that $U_{g,m,V_1}$ and $U_{g,m',V_2}$ are leaf-open neighbourhoods of $(g,m)$ and $(g,m')$ respectively \recalt{localleafcomponents}.
By taking a smaller $V_2$ if necessary, we may assume that $U_{g,m',V_2}$ is contained in the open set $U$.
If we then define $V=V_1\cap V_2$, it follows that $U_{g,m',V}$ is a leaf-open neighbourhood of $(g,m')$ contained in $U$ and $U_{g,m,V}$ a leaf-open neighbourhood of $(g,m)$.

We claim that $U_{g,m,V}$ is disjoint from $U$. So suppose we have $(h,m'')\in U\cap U_{g,m,V}$.
But $p(U_{g,m,V})=V\cdot g = p(U_{g,m',V})\subset p(U)= \Ub$ and $p$ is a bijection from $U$ to $\Ub$, and thus from $U_{g,m',V}$ to $V\cdot g$. 
Hence $(h,m'')\in U_{g,m',V}$ and we must have the equality
$$
\psi_{g,m,V}(h) = (h,m'') = \psi_{g,m',V}(h)
\mapob,
$$
\ie, we have the equality
$$
(p_{23}\scirc \Phi_A)\bigl(1, \exp\mo(hg\mo), g,m\bigr)
=
(p_{23}\scirc \Phi_A)\bigl(1, \exp\mo(hg\mo), g,m'\bigr)
\mapob.
$$
It follows immediately from the explicit form of the flow $\Phi_A$ \recalf{expressionforPhisubA} that we must have the equality
$$
\Phi_A\bigl(1, \exp\mo(hg\mo), g,m\bigr)
=
\Phi_A\bigl(1, \exp\mo(hg\mo), g,m'\bigr)
\mapob.
$$
But then we can use the property of a flow, applying $\Phi_A(-1, \ \cdot\ )$, to conclude that we have the equality
$$
\bigl(1, \exp\mo(hg\mo), g,m\bigr)
=
\bigl(1, \exp\mo(hg\mo), g,m'\bigr)
\mapob,
$$
and thus in particular $m=m'$. 
This contradicts the initial choice $(g,m)\notin U$ and thus $U'=U_{g,m,V}$ is the sought-for open neighbourhood of $(g,m)$ disjoint from $U$.
\end{preuve}

\begin{proclaim}[localactiongivesleafs]{Lemma}
Let $\Psi:D\to M$ be a local action of $G$ on $M$. When we define $D_m\subset G$ and $\psi_m:D_m\to G\times M$ for a fixed $m\in M$ by
$$
D_m = \{\,g\in G\mid (g,m)\in D\,\}
\qquad\text{and}\qquad
\psi_m(g) = \bigl( g,\Psi(g,m)\bigr)
\mapob,
$$
then $\psi_m(D_m)$ is leaf-open and $p:\psi_m(D_m)\to D_m$ is bijective.

\end{proclaim}

\begin{preuve}
It is immediate from the definition of $\psi_m$ that $p:\psi_m(D_m)\to D_m$ is bijective, as $\psi_m$ is injective and $p\scirc\psi_m = id\caprestricted_{D_m}$.
We now choose $g\in D_m$ and, without worrying for the moment about the exact domain of applicability, we make the computation
\begin{align*}
\psi_m(h)
&
=
\bigl(h, \Psi(h,m)\bigr)
\\
&
=
\bigl(\,h\,,\, \Psi\bigl(hg\mo, \Psi(g,m)\bigr)\,\bigr)
\\
&
\overset{\text{\rlap{\hss\recaltt{linkPsiandPhiMM}{linkbetweenPsiandPhiM}\vrule depth3pt width0pt}}}{=}
\quad
\bigl(\,h\,,\, \Phi_M^M(1,\exp\mo(hg\mo),\Psi(g,m)\bigr)\,\bigr)
=
\psi_{g,\Psi(g,m),V}(h)
\mapob.
\end{align*}
Now the last item is defined on the open neighbourhood $V{\cdot} g$ of $g$ (with $V$ satisfying conditions) and the first one on $D_m$. As $D_m$ is open, we can take a smaller $V$ if necessary such that both $\psi_m$ and $\psi_{g,\Psi(g,m),V}$ apply.
Remains to show that the intermediate steps are justified.

As the first and last equalities are the definitions, we start with the second equality, which uses the group property of a local action. 
For this to be true, we must have $(g,m), (h,m), \bigl(hg\mo,\Psi(g,m)\bigr)\in D$. 
As we have $g,h\in V{\cdot} g\subset D_m$, we indeed have $(g,m), (h,m)\in D$. 
Now $\bigl(e,\Psi(g,m)\bigr)\in D$, so there exists an open neighbourhood $V'$ of $e$ such that $V'\times \{\Psi(g,m)\}\subset D$. 
By taking a smaller $V$ if necessary, we may assume $V\subset V'$.
And as $h\in V{\cdot} g$, we have $hg\mo\in V\subset V'$ and thus $\bigl(hg\mo,\Psi(g,m)\bigr)\in D$ as wanted.
This justifies the second equality for a sufficiently small neighbourhood $V$ of $e\in G$.

For the third equality, we want to apply \recaltt{linkPsiandPhiMM}{linkbetweenPsiandPhiM} to 
$(1,\exp\mo(hg\mo),m)$. 
We thus have to show that this belongs to $W_\Psi$.
But $W_\Psi$ is open and contains $(1,0,m)$ \recaltt{transferofactiondomain}{linkbetweenPsiandPhiM}.
Hence there exists an open neighbourhood $U'$ of $0\in \Liealg g_0$ such that $\{1\}\times U'\times \{m\}\subset W_\Psi$.
By taking (once again) a smaller $V$ if necessary, we may assume that $\exp\mo(V)\subset U'$.
And thus we may apply \recaltt{linkPsiandPhiMM}{linkbetweenPsiandPhiM} to conclude that our computation is fully justified on this smallest $V{\cdot} g$.

But then by \recalt{localleafcomponents} $U_{g,\Psi(g,m),V}=\psi_m(V{\cdot} g)$ is a leaf-open neighbourhood of $\psi_m(g)$ contained in $\psi_m(D_m)$, proving that $\psi_m(D_m)$ is leaf-open.
\end{preuve}

\begin{proclaim}[uniquenessoflocalactions]{Proposition}
If $\Psi_i:D\to M$, $i=1,2$ are two local actions (defined on a same action domain) integrating the same infinitesimal action $\rho$, then $\Psi_1=\Psi_2$.

\end{proclaim}

\begin{preuve}
Fix $m\in M$ and define $D_{m}\subset G$ and $U_{i,m}\subset G\times M$ by
$$
D_{m} = \{\,g\in G\mid (g,m)\in D\,\}
\quad\text{and}\quad
U_{i,m} = \{\,\bigl(g,\Psi_i(g,m)\bigr) \mid g\in D_{m}\,\}
\mapob.
$$
According to \recalt{localactiongivesleafs} $U_{i,m}$ is leaf-open and $p$ projects it bijectively to $D_{m}$ which is connected (because $D$ is an action domain).
Hence by \recalt{opentoconnectedisconnectedcomponent} $U_{i,m}$ is a leaf-connected component of $p\mo(D_m)$. 
But $\Psi_i(e,m) = m$ and thus $(e,m)\in U_{1,m}\cap U_{2,m}$ and thus these two leaf-connected components of $p\mo(D_m)$ must be the same. 
As $p$ maps $U_{i,m}$ bijectively to $D_m$, it follows immediately that we have $\Psi_1(g,m)=\Psi_2(g,m)$ for all $g\in D_m$.
\end{preuve}

\section{Extending smoothness}

\begin{proclaim}[actionwillbesmooth]{{Lemma}}
Let $G$ be a connected Lie supergroup, $M$ a supermanifold and let $\Psi:D\to M$ be a (set theoretic) local left-action of $G$ on $M$.
If for all $m\in M$ there exist open $V_m\subset G$ and $U_m\subset M$ satisfying the two conditions
\begin{enumerate}
\item
$(e,m)\in V_m\times U_m\subset D$ and

\item
the restriction of $\Psi$ to $V_m\times U_m$ is smooth,

\end{enumerate}
then $\Psi$ is smooth on $D$.

\end{proclaim}

\langbewijs
\begin{preuve}
The idea to prove that $\Psi$ is smooth in a neighbourhood of an arbitrary point $(g,m)\in D$ is to write $g$ as the product of (a finite number of) ``small'' elements as
$$
g= h_N\cdot h_{N-1} \cdots h_2\cdot h_1
$$
and to use the group property to write $\Psi(hg,m)$ as
$$
\Psi(hg,m)
=
\Psi\bigl(h,\Psi\bigl(h_N,\dots,\Psi(h_1,m)\dots\bigr)\bigr)
\mapob.
$$
If we then fix the $h_i$ and let $h$ run through a neighbourhood of $e\in G$, the product $hg$ runs through a neighbourhood of $g$. If the successive stages all belong to the part where $\Psi$ is smooth, it will also be smooth when we fix the first variable. And thus in a neighbourhood of $(g,m)$ we will have written $\Psi$ as the composition of the smooth maps $m'\mapsto \Psi(h_i,m')$ and the (final) smooth map $(h,m')\mapsto\Psi(h,m')$.
The biggest obstacle to this idea is that there is no obvious way to decompose $g$ as such a product, because the size of the neighbourhood of $e\in G$ where the map $m'\mapsto\Psi(h,m')$ is smooth depends upon the point $m'$.
The major part of the proof thus will consist in constructing such a decomposition.

\medskip

But before we start with the construction of the decomposition, we first make life somewhat easier with respect to smoothness problems.
If $V$ is an open neighbourhood of $e\in G$, then $V{\cdot} g$ is an open neighbourhood of $g\in G$, just because right translation is a homeomorphism (and for $g\in \body G$ it also is a diffeomorphism).
But by definition of the DeWitt topology, any open set $O$ (in any supermanifold) satisfies the equality $O = \body\mo\bigl(\body O\bigr)$, which immediately implies that we have the equality
$$
V{\cdot} g = V{\cdot} \body g
\mapob.
$$
It follows that it suffices to prove that $\Psi$ is smooth on open neighbourhoods of the form $(V{\cdot} g)\times U$ with $g\in \body G$, $V$ an open neighbourhood of $e\in G$ and $U$ an open neighbourhood of an arbitrary $m\in M$ (it would suffice to use arbitrary $m\in \body M$, but that will not simplify our argument a jot).

We thus choose $g\in \body G$ and $m\in M$ arbitrary such that $(g,m)\in D$.
By definition of an action domain, the set $D_m$ (defined as in \recalt{localactiongivesleafs}) is open and connected. As $G$ is locally path connected, it follows that $D_m$ and $\body D_m$ are path connected.
Hence there exists a continuous map $\gamma:[0,1]\to\body D_m\subset D_m$ such that $\gamma(0)=e$ and $\gamma(1)=g$.
Associated to this curve we define, for each $t\in [0,1]$, the points $m_t$ by
$$
m_t= \Psi\bigl(\gamma(t), m\bigr)
\mapob.
$$
We now claim that there exists a finite sequence $0= t_0\le t_1\le \cdots \le t_N=1$ and open sets $V_i\subset G$ and $U_i\subset M$ satisfying the three conditions
\begin{enumerate}
\item
$(e,m_{t_i})\in V_i\times U_i\subset D$,

\item
the restriction of $\Psi$ to $V_i\times U_i$ is smooth and

\item
$\gamma(t_{i+1})\in V_i\cdot \gamma(t_i)$.

\end{enumerate}
The first two conditions are ``trivial'' because they are given by hypothesis. The problem lies in the third condition.

To prove our claim, we consider the set $A$ of all endpoints of such sequences:
\begin{multline*}
A=\{\,\tau\in [0,1] \mid \exists N\ \exists 0= t_0\le \cdots\le t_N=\tau 
\\
\text{ satisfying the $3$ conditions }\,\}
\end{multline*} 
and we claim that we have $1\in A$.
To prove this claim we define $T=\sup A$ and we assume (proof by contradiction) that either $T<1$ or $T=1\notin A$. 
We then invoke the hypothesis that there exist $V_{m_T}\ni e$ and $U_{m_T}\ni m_T$ such that $\Psi$ is smooth on $V_{m_T}\times U_{m_T}$.
As $\Psi(e,m_T)=m_T$, and as $\Psi$ is smooth on $V_{m_T}\times U_{m_T}$, there exist open neighbourhoods $e\in\Vh_{m_T}\subset V_{m_T}$ and $m_T\in \Uh_{m_T}\subset U_{m_T}$ such that 
$$
\Psi\bigl(\, \Vh_{m_T}\times \Uh_{m_T}\,\bigr) \subset U_{m_T}
\mapob.
$$
Using continuity of the multiplication and the inverse, there exists an open neighbourhood $e\in \Vt_{m_T}\subset \Vh_{m_T}$ such that $\Vt_{m_T} \cdot \Vt_{m_T}\mo\subset \Vh_{m_T}$. 
Now $\gamma$ is continuous, so there exists $\varepsilon>0$ such that 
\begin{equation}\label{aroundTineneighbourhood}
\gamma\bigl( \,(T-\varepsilon, T+\varepsilon)\cap[0,1]\,\bigr)\subset \Vt_{m_T}\cdot \gamma(T)
\mapob.
\end{equation}
By definition of $T$, there exists $\tau\in A$, $T-\varepsilon<\tau\le T$, \ie, there exist $0= t_0\le t_1\le \cdots \le t_N=\tau$ and open sets $V_i$ and $U_i$ satisfying the three conditions above. 
We then define the sets $V_{N+1}=V_{N+2}=\Vh_{m_T}$ and $U_{N+1}=U_{N+2}=U_{m_T}$ as well as the points $t_{N+1}=t_N$ and $t_{N+2}$ as
$$
t_{N+2} = 
\begin{cases}
T+\delta &\quad\text{if $T<1$ with $0<\delta<\min(1-T,\varepsilon)$}
\\
T\equiv 1 & \quad\text{if $T=1\notin A$.}
\end{cases}
$$
As we have $\Vh_{m_T}\times U_{m_T}\subset V_{m_T}\times U_{m_T}$, $\Psi$ is smooth on $V_{N+1}\times U_{N+1}=V_{N+2}\times U_{N+2}$.
Moreover, we have
$$
\gamma(t_{N+1}) = \gamma(t_N)\in V_N\cdot \gamma(t_N)
\mapob,
$$
but also (using \recalf{aroundTineneighbourhood})
$$
\gamma(\tau)\equiv\gamma(t_N),\gamma(t_{N+2})\in \Vt_{m_T}\cdot \gamma(T)
\mapob,
$$
which implies that we have
$$
\gamma(t_{N+2})\in \Vt_{m_T}\cdot \Vt_{m_T}\mo \cdot \gamma(t_N)\subset \Vh_{m_T}\cdot \gamma(t_{N})
=
V_{N+1}\cdot \gamma(t_{N+1})
\mapob.
$$
Finally we note that we have (again using \recalf{aroundTineneighbourhood})
$$
\gamma(t_N), \gamma(t_{N+2})\in \Vt_{m_T}\cdot \gamma(T)
\quad\Leftrightarrow\quad
\gamma(t_N)\gamma(T)\mo, \gamma(t_{N+2})\gamma(T)\mo\in \Vt_{m_T}
$$
and thus
\begin{align*}
m_{t_{N+1}} 
&
=
m_{t_N}
= 
\Psi\bigl(\gamma(t_N),m\bigr)
=
\Psi\bigl(\gamma(t_N)\cdot \gamma(T)\mo, \Psi(\gamma(T),m)\bigr)
\\
&
=
\Psi\bigl(\gamma(t_N)\cdot \gamma(T)\mo, m_T\bigr)
\in U_{m_T}
=
U_{N+1}
\end{align*}
and
$$
m_{t_{N+2}} 
=
\Psi\bigl(\gamma(t_{N+2})\gamma(T)\mo, m_T\bigr)
\in U_{m_T}
=
U_{N+2}
\mapob.
$$
One should note that these computation are allowed because all couples concerned belong to $D$. 
We thus have shown that the sequence extended to $t_{N+2}$ also satisfies the three conditions, and thus $t_{N+2}\in A$. In the case $T<1$ this contradicts the definition of $T$ (because $T<t_{N+2}$) and in the case $T=1$ this contradicts $1\notin A$ (because in that case $t_{N+2}=1$).
This finishes the proof of our claim.

\medskip

Once we know the existence of a sequence $0= t_0\le t_1\le \cdots\le t_N=1$ satisfying the three conditions, we define the elements $h_i\in \body G$, $1\le i\le N$ by (and remember: $\gamma$ takes its values in $\body G$ and $g=\gamma(1)$)
$$
h_i=\gamma(t_{i})\gamma(t_{i-1})\mo
\mapob,
$$
where one should note that, because $t_0=0$ and $\gamma(0) = e$, we have $h_1=\gamma(t_1)$.
It follows that we have the equality
\begin{equation}\label{gasproductofhs}
g=
h_N\cdot h_{N-1} \cdots h_1 
\mapob.
\end{equation}
Moreover, the condition $\gamma(t_{i+1})\in V_i\cdot \gamma(t_i)$ implies that we have
$$
h_i=\gamma(t_i)\gamma(t_{i-1})\mo\in V_{i-1}
\mapob.
$$
As $\Psi$ is smooth on $V_{i-1}\times U_{i-1}$ and as $h_i\in V_{i-1}$ has real coordinates (it belongs to $\body G$), it follows that the map $\Psi^{(i)}:U_{i-1}\to M$ defined by
$$
\Psi^{(i)}(m') = \Psi(h_i,m')
$$
is smooth on $U_{i-1}$.
Moreover, using the group property of $\Psi$, we have
\begin{equation}\label{recursiononmti}
m_{t_i} = \Psi\bigl(\gamma(t_i),m\bigr)
=
\Psi\bigl(\, h_i, \Psi\bigl(\gamma(t_{i-1}),m\bigr)\,\bigr)
=
\Psi^{(i)}(m_{t_{i-1}})
\mapob.
\end{equation}
Using \recalf{gasproductofhs}, we thus have
\begin{align*}
\Psi(g,m) 
&
\overset{t_N=1}{=} 
m_{t_N}
=
(\Psi^{(N)}\scirc \Psi^{(N-1)}\scirc \cdots\scirc \Psi^{(1)})(m_{t_0})
\\
&
\overset{t_0=0}{=}
(\Psi^{(N)}\scirc \Psi^{(N-1)}\scirc \cdots\scirc \Psi^{(1)})(m)
\mapob.
\end{align*}
But this expression is valid for the given fixed couple $(g,m)\in D$. In order to extend it to a neighbourhood, we first write it in the form
$$
\Psi(e\cdot g,m) 
=
\Psi\bigl(e,(\Psi^{(N)}\scirc \Psi^{(N-1)}\scirc \cdots\scirc \Psi^{(1)})(m)\bigr)
\mapob,
$$
which suggests the formula
$$
\Psi( g',m') 
=
\Psi\bigl(R_{g\mo}(g'),(\Psi^{(N)}\scirc \Psi^{(N-1)}\scirc \cdots\scirc \Psi^{(1)})(m')\bigr)
\mapob,
$$
where $R_{g\mo}h=hg\mo$ denotes right-translation over $g\mo$, which is a diffeomorphism because $g$ has real coordinates by assumption.
It now only remains to show that this expression is valid and smooth on a neighbourhood of $(g,m)$.

To find a suitable neighbourhood we first note that $\Psi$ is smooth on $V_N\times U_N$, where $V_N$ is a neighbourhood of $e\in G$ and $U_N$ an open neighbourhood of $m_{t_N}=\Psi\bigl(\gamma(t_N),m\bigr) = \Psi(g,m)$.
As the maps $\Psi^{(i)}$ and $R_{g\mo}$ are smooth (on suitable neighbourhoods), our formula will give a smooth map (by composition of smooth maps), provided we have 
$$
R_{g\mo}(g')\in V_N
\qquad\text{and}\qquad
(\Psi^{(N)}\scirc \Psi^{(N-1)}\scirc \cdots\scirc \Psi^{(1)})(m')\in U_N
\mapob.
$$
The first condition translates as $g'\in V_N\cdot g$, which is indeed an open neighbourhood of $g\in G$.
To find a neighbourhood of $m$ such that the second condition is satisfied, we are going to shrink by backward recursion the neighbourhood $U_0$.

We start with the observation that $\Psi^{(N)}:U_{N-1}\to M$ is smooth with (using \recalf{recursiononmti}) $\Psi^{(N)}(m_{t_{N-1}}) = m_{t_N}\in U_N$. Hence there exists an open neighbourhood $\Uh_{N-1}\subset U_{N-1}$ such that 
$$
m_{t_{N-1}}\in \Uh_{N-1}
\qquad\text{and}\qquad
\Psi^{(N)}(\Uh_{N-1})\subset U_N
\mapob.
$$
The next step is to note that $\Psi^{(N-1)}$ is smooth on $U_{N-2}$ and that we have (again using \recalf{recursiononmti}) $\Psi^{(N-1)}(m_{t_{N-2}}) = m_{t_{N-1}} \in \Uh_{N-1}$. Hence there exists an open neighbourhood $\Uh_{N-2}\subset U_{N-2}$ such that 
$$
m_{t_{N-2}}\in \Uh_{N-2}
\qquad\text{and}\qquad
\Psi^{(N-1)}(\Uh_{N-2})\subset \Uh_{N-1}
$$
and thus
$$
m_{t_{N-2}}\in \Uh_{N-2}
\qquad\text{and}\qquad
(\Psi^{(N)}\scirc \Psi^{(N-1)})(\Uh_{N-2})\subset U_N
\mapob.
$$
Continuing this way with decreasing $i$, we find neighbourhoods $\Uh_i \subset U_i$ of $m_{t_i}$ such that we have (for $1\le i<N$)
$$
m_{t_{i-1}}\in \Uh_{i-1}
\quad\text{and}\quad
(\Psi^{(i)}\scirc \Psi^{(i+1)}\scirc\cdots \scirc \Psi^{(N)})(\Uh_{i-1})\subset U_N
\mapob.
$$
In particular there exists an open neighbourhood $\Uh_0$ of $m=m_0=m_{t_0}$ such that 
$$
(\Psi^{(1)}\scirc \Psi^{(2)}\scirc\cdots \scirc \Psi^{(N)})(\Uh_{0})\subset U_N
\mapob.
$$
The final conclusion is that $\Psi$ is smooth on the neighbourhood $(V_N{\cdot} g) \times \Uh_0\subset D$ of $(g,m)$ and hence $\Psi$ is smooth on $D$ as desired.
\end{preuve}
\eindelangbewijs

\kortbewijs
\begin{preuve}
This is a direct consequence of the fact that $\Phi$ is a left-action and the fact that an open neighbourhood of the identity generates the connected component of $G$.
\end{preuve}
\eindekortbewijs

\section{Existence}

In the introduction we argued that a local action is completely determined by subsets $\Lambda_m$ of the leafs $L_{(e,m)}$ of $\foliation$ subject to some conditions. 
It thus seems natural to choose maximal subsets $\Lambda_m$ within the given constraints.
Unfortunately, this is not as easy as one might think.
To illustrate (some of) the problems, we look at example \recalt{examp2} with $\lambda=\tfrac23$.
Looking at the constraint that $p:\Lambda_m\to G$ must be injective, we might be tempted to define these sets by
$$
\Lambda_x = \{\, (t+\ZZ, x+\tfrac23\, t) \mid 
t\in 
D_x 
\,\}
\mapob,
$$
where the open interval $D_x\subset \RR$ of length $1$ is given by
$$
D_x = 
\begin{cases}
\bigl(-\tfrac32\,x , 1-\tfrac32\,x\bigr) & 0<x<\tfrac13
\\
\bigl(-\tfrac12, \tfrac12\bigr) & \tfrac13\le x\le \tfrac23
\\
\bigl(\,\tfrac32\,(1-x)-1,\tfrac32\,(1-x)\bigr) & \tfrac23<x<1
\mapob.
\end{cases} 
$$
This immediately gives us for the associated action domain the set
$$
D = \dcup_{x\in (0,1)} D_x\times \{x\} =  \{\, (t+\ZZ,x) \mid t\in D_x   \,\}
\mapob.
$$

\begin{figure}[htb]
\null\hfil
\includegraphics[width=0.8\textwidth]{\drawingpath 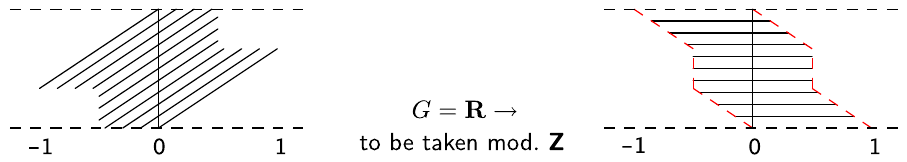}
\hfil
\null

\vskip-\baselineskip

{\null\hfil\hfil\hfil\hfil\hfil
The sets $\Lambda_x$ \kern4cm The sets $D_x$ and the domain $D$ \hfil}
\end{figure}

\noindent
And indeed $p:\Lambda_x\to D_x\subset G$ is injective (and it is not injective on any bigger set), but there does not exist a local action defined on this $D$ integrating the given infinitesimal action. 
To see why, suppose that $\Psi:D\to M$ is such a local action.
We then consider the point $x=\tfrac14$ and the ``times'' $s=t=\tfrac13+\ZZ$, for which we have $s+t=\tfrac23+\ZZ=-\tfrac13+\ZZ$. 
We then note that the points $(s,x)$ and $(s+t,x)$ belong to $D$ with $D_{1/4} = \bigl( -\tfrac38, \tfrac58\bigr)$. 
As $\bigl(t,\Psi(t,\tfrac14)\bigr)$ should lie on $\Lambda_{1/4}$, it follows immediately that we have
$$
\Psi(t,x) = \tfrac14 + \tfrac23\cdot \tfrac13 = \tfrac{17}{36}
\qquad\text{and}\qquad
\Psi(s+t,x) = \tfrac14 - \tfrac23\cdot \tfrac13 = \tfrac1{36}
\mapob.
$$
We now note that the point $\bigl(s,\tfrac{17}{36}\bigr) = \bigl(s,\Psi(t,x)\,\bigr)$ belongs to $D$ and that we have, as $\bigl(s,\Psi(s,\tfrac{17}{36})\bigr)$ should lie on $\Lambda_{17/36}$, 
$$
\Psi\bigl(s,\tfrac{17}{36}\bigr) = \tfrac{17}{36} + \tfrac23\cdot \tfrac13 = \tfrac{25}{36}
\mapob.
$$
But if $\Psi$ is a local action, it should satisfy the group law and in particular we should have
\begin{align*}
\tfrac1{36}
&
=
\Psi(s+t, x) 
=
\Psi\bigl( s, \Psi(t,x) \bigr)
=
\Psi\bigl( s,\tfrac{17}{36}  \bigr)
=
\tfrac{25}{36}
\mapob.
\end{align*}
This contradiction shows that this $D$ cannot be an action domain for a local action integrating the given infinitesimal action.

Close inspection shows that the various $\Lambda_x$ are not disjoint as one would expect.
It thus is tempting to think that reducing the $\Lambda_x$ in such a way that they become disjoint should suffice to establish the group law.
This idea is wrong on two counts.
The first reason is that it does not work as expected: if we reduce the intervals $D_x$ as follows:
$$
D_x = \bigl(\, \max\bigl(-\tfrac12,\tfrac32\,(x-1)\,\bigr)\ ,\ \min\bigl(\tfrac12, \tfrac32\,x\bigr) \,\bigr)
\mapob,
$$
then the sets $\Lambda_x$ no longer intersect,

\begin{figure}[htb]
\null\hfil
\includegraphics[width=0.8\textwidth]{\drawingpath 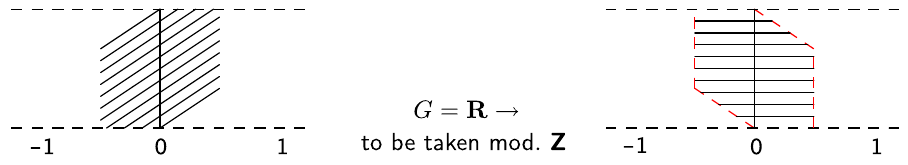}
\hfil
\null

\vskip-\baselineskip

{\null\hfil\hfil\hfil\hfil\hfil
The sets $\Lambda_x$ \kern4cm The sets $D_x$ and the domain $D$ \hfil}
\end{figure}

\noindent
but the points $(t,x)$, $(s+t,x)$ and $\bigl(s,\Psi(t,x)\bigr)$ also belong to this smaller potential action domain, so we still get a contradiction.

The second reason that the no-intersection idea is wrong is because the $\Lambda_x$ can intersect! 
To see this, consider the (mostly smaller) intervals $D_x$ defined by
$$
D_x = 
\begin{cases}
\bigl(-\tfrac32\,x , \tfrac16\bigr) & 0<x\le\tfrac19
\\
\bigl(-\tfrac16 , \tfrac7{12}\bigr) & \tfrac19< x<\tfrac29
\\
\bigl(-\tfrac16 , \tfrac16\bigr) & \tfrac29\le x\le\tfrac79
\end{cases} 
\qquad
D_x = 
\begin{cases}
\bigl(-\tfrac7{12} , \tfrac16\bigr) & \tfrac79< x<\tfrac89
\\
\bigl(-\tfrac16,\tfrac32\,(1-x)\bigr) & \tfrac89\le x<1
\mapob.
\end{cases} 
$$

\begin{figure}[htb]
\null\hfil
\includegraphics[width=0.8\textwidth]{\drawingpath 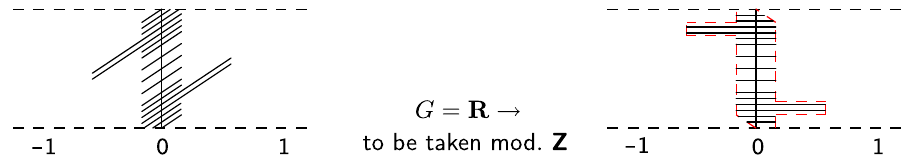}
\hfil
\null

\vskip-\baselineskip

{\null\hfil\hfil\hfil\hfil\hfil
The sets $\Lambda_x$ \kern4cm The sets $D_x$ and the domain $D$ \hfil}
\end{figure}

\noindent
With this choice, the sets $\Lambda_{1/6}$ and $\Lambda_{5/6}$ intersect (both contain the point $\bigl(\tfrac12+\ZZ, \tfrac12 \bigr)$), but one can show that the associated action domain $D$ indeed is an action domain for a local action $\Psi:D\to M$ integrating the given infinitesimal action.

The upshot of these examples is the observation that the choice of the sets $\Lambda_x$ that should give us the action domain is rather delicate, especially with respect to the group law.
The purpose of the following definition is to provide sufficiently small sets $\Lambda_x$ so that we can prove that the group law is satisfied.

\begin{definition}[defofneighbourhoodbase]{Definition\slash Notation}
Choose any norm on $\body\Liealg g_0$ with respect to which we define the open balls $B_r(x)\subset \body\Liealg g_0$ of radius $r$ and center $x$.
Using these we define the collection $\nbhdbase$ by
$$
\nbhdbase
=
\{\,\exp\bigl(\,\body\mo(B_{1/n}(0))\,\bigr)\mid n\in\NN^* \text{ and }\body\mo\bigl(B_{1/n}(0)\bigr)\subset \Domexp\,\}
\mapob.
$$
By definition of the DeWitt topology the open balls $\body\mo\bigl(B_{1/n}(0)\bigr)$ form a neighbourhood basis at $0\in \Liealg g_0$. Hence, because $\exp$ is a diffeomorphism from $\Domexp$ to $\Imexp$, $\nbhdbase$ is a neighbourhood basis at $e\in G$ for the topology on $G$.
This basis has the following properties (and any other neighbourhood basis at $e\in G$ with the same properties will do for our argument lateron):
\begin{enumerate}
\item
each element of $\nbhdbase$ is invariant under the map $g\mapsto g\mo$,

\item
each element of $\nbhdbase$ is connected, 

\item
on each element of $\nbhdbase$ the inverse exponential map (\ie, the logarithm) is a diffeomorphism onto its image,

\item
$\nbhdbase$ is totally ordered by inclusion,

\item
any subset $\nbhdbase'\subset \nbhdbase$ admits a maximal element (with respect to inclusion).

\end{enumerate}
Using this neighbourhood basis, we define\slash choose for each $m\in M$ the sets $\Vt_m''$,  $\Vt_m'$, $\Vt_m$, $U_m$, $V_m''$, $V_m'$ and $V_m$ as follows (in this order): 
\begin{itemize}
\item
$\Vt_m''$ is the maximal element in $\nbhdbase$ such that $\{1\}\times \exp\mo(\Vt_m'')\times \{m\}\subset W_M$,

\item
$\Vt_m'$ is the maximal element in $\nbhdbase$ such that $\Vt_m'\cdot \Vt_m' \subset \Vt_m''$,

\item
$\Vt_m$ is the maximal element in $\nbhdbase$ such that $ \Vt_m\cdot \Vt_m \subset \Vt_m'$,

\item
$U_m$ is an open neighbourhood of $m$ for which there exists an element $V_m''\in \nbhdbase$ such that $\{1\}\times \exp\mo(V_m'')\times U_m \subset W_M$,

\item
$V_m'$ is an element in $\nbhdbase$ such that $V_m' \cdot V_m' \subset V_m''$,

\item
$V_m$ is an element of $\nbhdbase$ such that $V_m\cdot V_m\subset V_m'$.

\end{itemize}
Such sets exists because $W_M$ is open and contains $\{1\}\times \{0\}\times M$ and because multiplication in $G$ is continuous.

\end{definition}

\begin{proclaim}[nestingDmsinVms]{Lemma}
If we have $(g,m)\in V_{m'}\times U_{m'}$, then we have $V_{m'}\subset \Vt_{m}$.

\end{proclaim}

\begin{preuve}
As we have the inclusions $V_{m'}\cdot V_{m'}\subset V_{m'}'$ and $V_{m'}'\cdot V_{m'}'\subset V_{m'}''$, we have $g\in V_{m'}''$. 
But we have 
$$
\{1\}\times \exp\mo(V_{m'}'')\times \{m\}\subset 
\{1\}\times \exp\mo(V_{m'}'')\times U_{m'}\subset W_M
\mapob,
$$
hence by maximality of $\Vt_m''$ we must have $V_{m'}''\subset \Vt_m''$.
But then we have $V_{m'}\cdot V_{m'}\subset \Vt_{m}''$, hence by maximality of $\Vt_m'$ we have $V_{m'}'\subset \Vt_m'$. Repeating this argument we obtain the desired conclusion $V_{m'}\subset \Vt_m$.
\end{preuve}

\begin{proclaim}[existenceoflocalactions]{Theorem}
The set $D=\smash{\dcup_{m\in M}} V_m\times U_m$ is an action domain and the map $\Psi:D\to M$ defined by
$$
\Psi(g,m) 
=
\Phi_M^M\bigl(1,\exp\mo(g), m\bigr)
$$
is a local action integrating $\rho$.

\end{proclaim}

\begin{preuve}
As $V_m$ belongs to $\nbhdbase$, it is a connected open neighbourhood of $e$. It follows immediately that $D$ is an action domain.
For $m\in M$ we have the inclusions
$$
V_{m} \subset V_{m}' \subset V_{m}''
$$
and thus the inclusions
$$
\{1\}\times \exp\mo(V_{m})\times U_{m}\subset \{1\}\times \exp\mo(V_{m}'')\times U_{m} \subset W_M
\mapob,
$$
which shows that $\Psi$ is well defined on $V_m\times U_m$. As $m\in M$ is arbitrary, it follows that $\Psi$ is well defined on $D$.
As $\Phi_M^M$ is smooth on $W_M$, it follows immediately that $\Psi$ is smooth (being the smooth map $\Phi_M^M$ with its first entry restricted to the real value $1$ and its second entry composed with the smooth map $\exp\mo$).

We now have to show that $\Psi$ verifies the conditions of a local action. 
We start with the observation that we have
$$
\Psi(e,m) = \Phi_M^M(1,\exp\mo(e), m) = \Phi_M^M(1,0,m) \overset{\recaltt{atzerodefinedforallt}{linkflowsZAZRandZM}}{=} m
\mapob,
$$
as required.
The last item to check is the group property. So suppose we have $(g,m)$,  $(hg,m)$ and $\bigl(h,\Psi(g,m)\bigr)\in D$. We then have to show the equality 
$$
\Psi(hg,m) = \Psi\bigl(h,\Psi(g,m)\bigr)
\mapob.
$$
By definition of $D$ there exist $m_1,m_2,m_3\in M$ such that we have
\begin{multline*}
(g,m)\in V_{m_1}\times U_{m_1}
\quad,\quad
(hg,m)\in V_{m_2}\times U_{m_2}
\\
\quad\text{and}\quad
\bigl(h,\Psi(g,m)\bigr)\in V_{m_3}\times U_{m_3}
\mapob.
\end{multline*}
According to \recalt{nestingDmsinVms} we thus have the inclusions
$$
V_{m_1}, V_{m_2}\subset \Vt_m
\qquad\text{and}\qquad
V_{m_3}\subset \Vt_{\Psi(g,m)}
\mapob.
$$
It follows that we have (using that the elements of $\nbhdbase$ are invariant under inverse)
$$
hg,g\mo\in \Vt_m
\qquad\text{and thus}\qquad
h=hgg\mo\in \Vt_m\cdot \Vt_m \subset \Vt_m'
\mapob.
$$
Now let $\Vt_o$ be the smallest of $\Vt_m'$ and $\Vt_{\Psi(g,m)}$ (which exists because $\nbhdbase$ is totally ordered by inclusion).
As we have $\Vt_o\subset \Vt_m'$ and $g\in \Vt_m\subset \Vt_m'$, we have the inclusion $\Vt_o\cdot g\subset \Vt_m''$.
It follows that the map \recalt{defoflocalleafcomponents} 
$$
F_1=\psi_{e,m,\Vt_m''}
$$ 
is defined on $\Vt_o\cdot g$.
On the other hand, we have the inclusion $\Vt_o\subset \Vt_{\Psi(g,m)}\subset \Vt_{\Psi(g,m)}'$ and thus the map 
$$
F_2=\psi_{g,\Psi(g,m),\Vt_{\Psi(g,m)}'}
$$ 
is also defined on $\Vt_o\cdot g\subset \Vt_{\Psi(g,m)}'\cdot g$.
According to \recalt{localleafcomponents} and \recalt{opentoconnectedisconnectedcomponent} we thus have homeomorphisms $F_i:\Vt_o\cdot g\to U_i = F_i(\Vt_o)$, where the $U_i$ are leaf-connected components of $p\mo(\Vt_o\cdot g)$.
But we have 
\begin{align*}
F_1(g) 
&
= \bigl(g,\Phi_M^M(1,\exp\mo(g),m)\bigr) = \bigl(g,\Psi(g,m)\bigr)
\\
F_2(g) 
&
= \bigl(g,\Phi_M^M(1,\exp\mo(gg\mo), \Psi(g,m))\bigr) = \bigl(g,\Psi(g,m)\bigr)
\mapob.
\end{align*}
The two leaf-connected components $U_1$ and $U_2$ thus have a point in common, hence are the same. 
It follows that we have the equality $F_1=F_2$. 
As we have $h\in \Vt_m' \cap \Vt_{\Psi(g,m)}$, we have $h\in \Vt_o$ and thus $F_1$ and $F_2$ are in particular defined at the point $hg\in \Vt_o\cdot g$, which means that we have in particular the equality
\begin{align*}
\bigl(hg,\Psi\bigl(h,\Psi(g,m)\bigr)\bigr) = F_2(hg)
=
F_1(hg) = \bigl(hg, \Psi(hg,m)\bigr)
\mapob.
\QEDici
\end{align*}

\end{preuve}

\section{Maximal and global actions}

\begin{definition}{Definitions}
Let $G$ be a Lie supergroup, $\Liealg g$ its Lie superalgebra and let $\rho$ be an infinitesimal action of $\Liealg g$ on a supermanifold $M$. 

$\bullet$
Following \cite{Pa57}, we will say that $(G,\rho)$ is \stresd{univalent} when the restriction of $p:G\times M\to G$ to any leaf of $\foliation$ (a leaf-connected component of $G\times M$) is injective.

$\bullet$
We will say that a local action $\Psi:D\to M$ integrating $\rho$ is \stresd{maximal} if for any local action $\Psi':D'\to M$ integrating $\rho$ we have $D'\subset D$.

\end{definition}

\begin{proclaim}[uniquenessofmaximallocalaction]{Lemma}
Let $\rho$ be an infinitesimal action of $\Liealg g$ on $M$.
If there exists a maximal local action integrating $\rho$, then it is unique.

\end{proclaim}

\begin{preuve}
Let $\Psi_i:D_i\to M$, $i=1,2$ be two maximal local actions integrating $\rho$. Then by definition of maximality we must have $D_1=D_2$. And then by \recalt{uniquenessoflocalactions} we must have $\Psi_1=\Psi_2$.
\end{preuve}

\begin{proclaim}[projectionofleavesisactiondomain]{Lemma}
Let $\rho$ be an infinitesimal action of $\Liealg g$ on $M$.
Then the set $D\subset G\times M$ defined as
$$
D=\dcup_{m\in M} p(L_{(e,m)})\times \{m\}
$$
is an action domain, where $L_{(e,m)}$ is the leaf-connected component of $G\times M$ containing $(e,m)$.

\end{proclaim}

\begin{preuve}
According to the definition \recalt{vectorfieldsonGtimesM}, we have to prove three affirmations: that $D$ is open in $G\times M$, that it contains $\{e\}\times M$ and that all sets $D_m$ are connected.
As $D_m$ is defined by $D_m\times \{m\} = G\times \{m\}\cap D$, it follows immediately that we have $D_m = p(L_{(e,m)})$. Now $p$ is leaf-continuous and $L_{(e,m)}$ is connected, so $D_m$ is connected as wanted.
At the same time we have $e=p(e,m)\in p(L_{(e,m)})=D_m$ and thus $(e,m)\in D$ for all $m\in M$, proving that we have the inclusion $\{e\}\times M\subset D$. Remains the hard part.

To prove that $D$ is open, we have to show that for $(g_o,m_o)\in D$ we can find open neighbourhoods $V_o\subset G$ of $g_o$ and $U_o\subset M$ of $m_o$ such that $V_o\times U_o\subset D$. 
Now $(g_o,m_o)\in D$ means that there exists $\mb_o\in M$ such that $(g_o,\mb_o)\in L_{(e,m_o)}$. 
We thus have to show that for all $g\in V_o$ and all $m\in U_o$ we can find $\mb\in M$ such that $(g,\mb)$ lies on the leaf $L_{(e,m)}$ through $(e,m)$.
Intuitively this means that we have to show that neighbouring leaves are not suddenly (much) smaller. 
\begin{figure}[htb]
\null\hfil
\includegraphics[width=0.9\textwidth]{\drawingpath 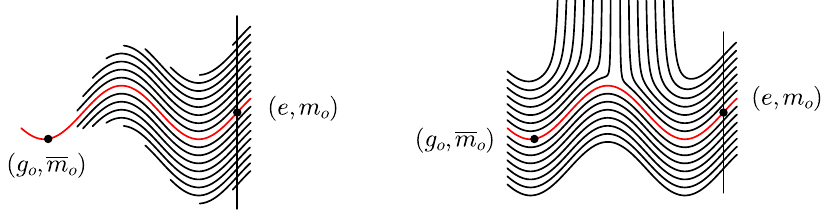}
\hfil
\null

\vskip-\baselineskip

{\null\hfil
Two situations that have to be excluded \hfil}
\end{figure}%
It is rather obvious that the pictures shown are not regular foliations, but there might be other, more vicious situations. And the property we have to prove is not a local one: we have to link leaves in a neighbourhood of $(g_o,\mb_o)$ to leaves in a neighbourhood of $(e,m_o)$.

The idea will be that we consider a path from $(e,m_o)$ to $(g_o,\mb_o)$ and that we (try to) follow neighbouring leaves along this path.
As each leaf-connected component is locally path-connected, it is globally path-connected. Hence there exists a continuous (with respect to the leaf-topology) map $\gamma:[0,1]\to L_{(e,m_o)}\subset G\times M$ such that $\gamma(0) = (e,m_o)$ and $\gamma(1) = (g_o,\mb_o)$.
According to \recalt{locadaptedcoordinates} there exists for each point $\gamma(t)\in G\times M$ an open neighbourhood $\Ut_t$ of $\gamma(t)$ with some special properties. 
As $[0,1]$ is compact and as $\gamma\mo(\Ut_t)$, $t\in [0,1]$ is an open cover of $[0,1]$, there exists $N\in\NN^*$ such that for all $1\le k\le N$ there exists $t_k\in [0,1]$ such that $\gamma\bigl(\,[(k-1)/N, k/N]\,\bigr)\subset \Ut_{t_k}$.

Changing notation slightly and resuming (some of) the results of \recalt{locadaptedcoordinates}, we thus have $N$ open neighbourhoods $\Ut^{(k)}$ with local coordinates $y_i^{(k)}, x_{j}^{(k)}$ adapted to $\foliation$, open neighbourhoods $V^{(k)}\subset G$ (which we may assume to be connected!) with local coordinates $y_i^{(k)}$, open sets $O^{(k)}\subset E_0$ with coordinates $x_j^{(k)}$ and diffeomorphisms $\varphi^{(k)}:\Ut^{(k)}\to V^{(k)}\times O^{(k)}$ with the properties 
\begin{enumerate}
\item
$\gamma\bigl(\, [\,(k-1)/N,k/N\,] \,\bigr) \subset \Ut^{(k)}$ and thus in particular
$\gamma\bigl((k-1)/N\bigr)$, $\gamma(k/N) \in \Ut^{(k)}$,

\item
$p=p_1\scirc \varphi^{(k)}$, 

\item
for all $ (g,m),(g',m')\in \Ut^{(k)}$ we have the implication (in general it is not an equivalence!)
\begin{multline*}
\qquad
(p_2\scirc \varphi^{(k)})(g,m) = (p_2\scirc \varphi^{(k)})(g',m')
\qquad\Longrightarrow\qquad
\\
(g,m) \text{ and } (g',m') \text{ lie on the same leaf,}
\qquad
\end{multline*}

\end{enumerate}
where $p_1:V^{(k)}\times O^{(k)}\to V^{(k)}$ and $p_2:V^{(k)}\times O^{(k)}\to O^{(k)}$ denote the canonical projections on the first and second factor.
We thus have in particular $\gamma(k/N)\in \Ut^{(k)}\cap \Ut^{(k+1)}$. 
Associated to these data we define the points $g^{(k)}\in V^{(k)}$ and $b^{(k)}\in O^{(k)}$ by the equalities
\begin{gather*}
g^{(k)} = p\bigl(\gamma(k/N)\bigr)
\qquad\text{and}\qquad
\varphi^{(k)}\bigl(\gamma(k/N)\bigr) = (g^{(k)}, b^{(k)})
\mapob.
\end{gather*}
As the leaf-connected components of $\Ut^{(k)}$ are the sets $V^{(k)}\times \{x\}$ (or rather their inverse image under $\varphi^{(k)}$) and as $\gamma$ is leaf-continuous, it follows that $\gamma\bigl(\,[(k-1)/N,k/N]\,\bigr)$ is contained in a single one of these \slice{}s.
This implies in particular that we have the equality
\begin{equation}\label{samebforsuccessivevarphi}
\varphi^{(k)}\bigl(\, \gamma\bigl((k-1)/N\bigr)\,\bigr)=
(g^{(k-1)}, b^{(k)})
\mapob.
\end{equation}

As we have $\gamma(k/N)\in \Ut^{(k)}\cap \Ut^{(k+1)}$ we have the ``change of coordinates'' map
$$
\varphi_{k+1,k} \equiv \varphi^{(k+1)}\scirc \bigl(\varphi^{(k)}\bigr)\mo 
:
\varphi^{(k)}(\Ut^{(k)}\cap \Ut^{(k+1)}) \to \varphi^{(k+1)}(\Ut^{(k)}\cap \Ut^{(k+1)})
\mapob.
$$
Because of the equality $p=p_1\scirc \varphi^{(k)}$, it follows that we have
$$
\varphi_{k+1,k}(g,\mu) = (g,\mu')
\mapob,
$$
with $\mu'$ a priori depending on $g$ and $\mu$.

Using the local coordinates $y_i^{(k)},x_j^{(k)}$ and $y_i^{(k+1)},x_j^{(k+1)}$ respectively, we can write the map $\varphi_{k+1,k}$ in terms of these coordinates as
\begin{align*}
\bigl(\,y_i^{(k+1)}\,,\,x_j^{(k+1)}\,\bigr) 
&
=
\Bigl(\,\,y_i^{(k+1)}(y_{i'}^{(k)}, x_{j'}^{(k)})\,\,,\,\,x_j^{(k+1)}(y_{i'}^{(k)}, x_{j'}^{(k)}) \,\,\Bigr) 
\\
&
= \varphi_{k+1,k}(y_{i'}^{(k)}, x_{j'}^{(k)})
\mapob.
\end{align*}
As both sets of local coordinates are adapted to $\foliation$, \ie, the $\partial_{y_{i'}^{(k)}}$ and $\partial_{y_i^{(k+1)}}$ both span $\foliation$, we must have 
\begin{equation}\label{Jacsubmatrixnull}
\fracp{x_j^{(k+1)}}{y_{i'}^{(k)}} = 0
\mapob,
\end{equation}
\ie, the functions $x_j^{(k+1)}$ are locally independent of $y_{i'}^{(k)}$ and in particular in a neighbourhood of $(g^{(k)}, b^{(k)})$. 
It follows that there exists an open neighbourhood $\Vh^{(k)}$ of $g^{(k)}$, an open neighbourhood $\Oh^{(k)}$ of $b^{(k)}$ and a (``reduced'') smooth map $\varphi_{k+1,k}^{r}:\Oh^{(k)}\to O^{(k+1)}$ such that the restriction of $\varphi_{k+1,k}$ to $\Vh^{(k)}\times \Oh^{(k)}\subset \varphi^{(k)}\bigl(\Ut^{(k)}\cap \Ut^{(k+1)}\bigr)$ is given by
\begin{equation}\label{locallinkleafcharts}
\varphi_{k+1,k}\restricted_{\Vh^{(k)}\times \Oh^{(k)}}:
(g,\mu)\mapsto \bigl(g,\varphi_{k+1,k}^{r}(\mu)\bigr)
\mapob.
\end{equation}
\begin{figure}[htb]
\null\hfil
\includegraphics[width=0.9\textwidth]{\drawingpath 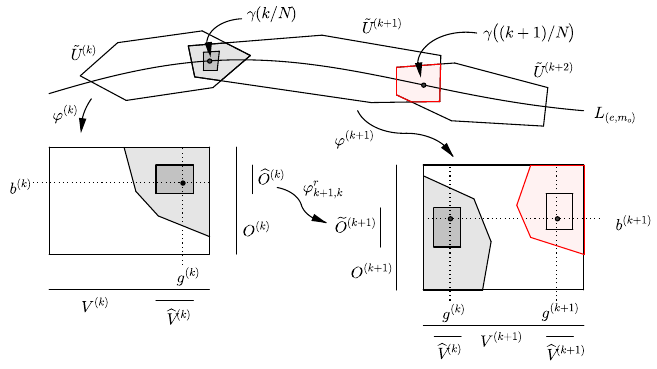}
\hfil
\null
\end{figure}%
But $\varphi_{k+1,k}$ is a diffeomorphism, hence its restriction to $\Vh^{(k)}\times \Oh^{(k)}$ is a diffeomorphism onto its image, which obviously is $\Vh^{(k)}\times \varphi_{k+1,k}^{r}(\Oh^{(k)})$. It follows that $\varphi_{k+1,k}^{r}$ is a diffeomorphism from $\Oh^{(k)}$ onto its open image $\Ot^{(k+1)}=\varphi_{k+1,k}^{r}(\Oh^{(k)}) \subset O^{(k+1)}$.
We finally invoke \recalf{samebforsuccessivevarphi} to deduce that we have 
$$
\varphi_{k+1,k}(g^{(k)},b^{(k)}) = (g^{(k)}, b^{(k+1)})
$$
and thus $\Ot^{(k+1)}$ is an open neighbourhood of $b^{(k+1)}$.
By taking successive preimages, we define the open set $\Oc^{(1)}\subset \Oh^{(1)}$ by
$$
\Oc^{(1)}
=
\bigl[\,(\varphi_{2,1}^{r})\mo \scirc (\varphi_{3,2}^{r})\mo \scirc \cdots \scirc (\varphi_{N,N-1}^{r})\mo \,\bigr]\bigl(\,\Ot^{(N)}\,\bigr)
\mapob,
$$
and then inductively $\Oc^{(k+1)}=\varphi_{k+1,k}^{r}(\Oc^{(k)}) \subset \Oh^{(k+1)}\cap \Ot^{(k+1)}$.
It follows in particular that $\varphi_{k+1,k}^{r}$ is a diffeomorphism from $\Oc^{(k)}$ onto $\Oc^{(k+1)}$.

With these preparations we are finally able to define our neighbourhoods $U_o\subset M$ of $m_o$ and $V_o\subset G$ of $g_o$ as required in the beginning of this proof.
For $V_o$ we simply define it as $\Vh^{(N)}$, whereas we define $U_o$ implicitly by
$$
\{e\}\times U_o = \bigl(\,\{e\}\times M\,\bigr)\ \cap\  \bigl(\varphi^{(1)}\bigr)\mo\bigl(\,\Vh^{(1)}\times \Oc^{(1)}\,  \bigr)
\mapob.
$$
To prove that these neighbourhoods will do, we simply check the requirement. We thus take $(g,m)\in V_o\times U_o$ and we define the sequence of points $\mu^{(k)}\in \Oc^{(k)}$ as follows.
We start with the observation that, by definition of $U_o$, we have $\varphi^{(1)}(e,m)\in \Vh^{(1)}\times \Oc^{(1)}$, and thus we can define $\mu^{(1)}\in \Oc^{(1)}$ by the equation
$$
\varphi^{(1)}(e,m) = (e,\mu^{(1)})
\mapob.
$$
Next we define inductively 
$$
\mu^{(k+1)} = \varphi_{k+1,k}^{r}(\mu^{(k)})
\mapob.
$$
Associated to these points we define the points $m^{(k)}\in M$ by
\begin{equation}\label{recursivedefofmuk}
\varphi^{(k)}(g^{(k)},m^{(k)}) = (g^{(k)}, \mu^{(k)})
\end{equation}
and, last but not least, the point $\mb$ by
$$
\varphi^{(N)}(g,\mb) = (g,\mu^{(N)})
\mapob.
$$
Combining \recalf{locallinkleafcharts} and \recalf{recursivedefofmuk} we immediately have the equality
$$
\varphi^{(k+1)}(g^{(k)},m^{(k)}) = (g^{(k)}, \mu^{(k+1)})
\mapob.
$$
We thus have the scheme of points
{\def\vertskip{3\jot}
$$
\kern-2.5em
\begin{matrix}
\scriptstyle(e,m)
\\
\noalign{\vskip3\jot}
\downarrow\rlap{\raise2pt\hbox{$\scriptstyle\varphi^{(1)}$}}
\\
\noalign{\vskip3\jot}
\scriptstyle(e,\mu^{(1)})
\\
\noalign{\vskip\vertskip}
\end{matrix}
\kern1.5em
\begin{matrix}
&\rlap{\hss$\scriptstyle(g^{(1)},m^{(1)})$} 
\\
\noalign{\vskip3\jot}
&\smash{\llap{\raise6pt\hbox{$\scriptstyle\varphi^{(1)}$}\kern-2pt$\swarrow$\kern1.5em}} \kern0.1em \smash{\rlap{\kern1.5em$\searrow$\kern1pt\raise6pt\hbox{$\scriptstyle\varphi^{(2)}$}}}
\\
\noalign{\vskip3\jot}
\scriptstyle(g^{(1)}, \mu^{(1)}) & \smash{\underset{\scriptstyle id\times\varphi_{2,1}^{r}}{{-}\!{-}\!{\to}}} & \scriptstyle(g^{(1)},\mu^{(2)})
\\
\noalign{\vskip\vertskip}
\end{matrix}
\kern0.5em \cdots\kern0.5em
\begin{matrix}
&\rlap{\hss$\scriptstyle(g^{(N-1)},m^{(N-1)})$} 
\\
\noalign{\vskip3\jot}
&\smash{\llap{\raise6pt\hbox{$\scriptstyle\varphi^{(N-1)}$}\kern-2pt$\swarrow$\kern1.5em}} \kern0.1em \smash{\rlap{\kern1.5em$\searrow$\kern1pt\raise6pt\hbox{$\scriptstyle\varphi^{(N)}$}}}
\\
\noalign{\vskip3\jot}
\scriptstyle(g^{(N-1)}, \mu^{(N-1)}) & \smash{\underset{\rlap{\hss$\scriptstyle id\times\varphi_{N,N-1}^{r}$}}{{-}\!{-}\!{\to}}}  & \scriptstyle(g^{(N-1)},\mu^{(N)})
\\
\noalign{\vskip\vertskip}
\end{matrix}
\kern1.5em
\begin{matrix}
\scriptstyle(g,\mb)
\\
\noalign{\vskip3\jot}
\llap{\raise2pt\hbox{$\scriptstyle \varphi^{(N)}$}}\downarrow
\\
\noalign{\vskip3\jot}
\scriptstyle(g,\mu^{(N)})
\\
\noalign{\vskip\vertskip}
\end{matrix}
$$
}
When we now apply property (iii) of the local charts $\Ut$ as cited above, it follows immediately that the sequence of points $(e,m)$, $(g^{(k)}, m^{(k)})$ and $(g,\mb)$ lie on the same leaf. In other words, for $(g,m)\in V_o\times U_o$ we have found a point $\mb$ such that $(g,\mb)$ lies on the leaf passing through $(e,m)$ as required.
This proves the inclusion $V_o\times U_o\subset D$ and thus $D$ is open.
\end{preuve}

\begin{proclaim}[univalentthusmaximallocalint]{Theorem}
Let $(G,\rho)$ be univalent. Then there exists a unique maximal local action $\Psi:D\to M$ integrating $\rho$, where $D$ is defined according to \recalt{projectionofleavesisactiondomain} by
$$
D=\dcup_{m\in M} p(L_{(e,m)})\times \{m\}
\mapob.
$$

\end{proclaim}

\begin{preuve}
According to \recalt{projectionofleavesisactiondomain}, $D$ is indeed an action domain.
We thus have to prove existence, uniqueness and maximality of $\Psi$.
To prove maximality, assume that $\Psi':D'\to M$ is any local action integrating $\rho$.
As before, we introduce the sets $D_m, D'_m\subset G$ as
\begin{align*}
D_m
&
=\{\,g\in G\mid (g,m)\in D\,\} = p(L_{(e,m)})
\\ 
D'_m 
&
= \{\, g\in G\mid (g,m)\in D'\,\}
\mapob.
\end{align*}
Writing $\psi'_m(g)=\bigl(g,\Psi'(g,m)\bigr)$, it follows from \recalt{localactiongivesleafs} that $\Lambda=\psi'_m(D'_m)$ is leaf-open and $p$ projects it bijectively to $D'_m$, which is connected by definition of an action domain.
Hence by \recalt{opentoconnectedisconnectedcomponent} $\Lambda$ is a leaf-connected component of $p\mo(D'_m)$ and thus in particular leaf-connected. As it contains $(e,m)=\psi'(e)$, it follows that we have the inclusion $\Lambda\subset L_{(e,m)}$ and thus $D'_m=p(\Lambda)\subset p(L_{(e,m)})=D_m$.
This proves the inclusion $D'\subset D$.
It then follows from \recalt{uniquenessoflocalactions} that $\Psi$ and $\Psi'$ coincide on $D'\subset D$, simply because the restriction of $\Psi$ to $D'$ is again a local action integrating $\rho$. Hence $\Psi$ is maximal. Uniqueness of this maximal local action follows immediately (or from \recalt{uniquenessofmaximallocalaction}).

Remains existence.
We start with a set-theoretic definition of the action $\Psi:D\to M$. We thus choose $(g,m)\in D$ and we want to define $\Psi(g,m)\in M$.
By definition of $D$ we have $g\in p(L_{(e,m)})$ and because $(G,\rho)$ is univalent, the projection $p:L_{(e,m)}\to D_m$ is a bijection.
Hence there exists a unique $m'\in M$ such that $(g,m')\in L_{(e,m)}$. We then define $\Psi(g,m) = m'$:
$$
\Psi(g,m)=m'
\quad\Leftrightarrow\quad
(g,m') \in L_{(e,m)}
\quad\Leftrightarrow\quad
(g,m') = (p\restricted_{L_{(e,m)}})\mo(g)
\mapob.
$$
To show that this is a left-action, choose $g,h\in G$ and $m\in M$ such that $(g,m)$, $(hg,m)$, $\bigl(h,\Psi(g,m)\bigr)\in D$. Denoting $m'=\Psi(g,m)$ as above, the definition of $\Psi$ tells us that there exists $m''$ such that $(h,m'')\in L_{(e,m')}$:
$$
\Psi(g,m) = m' \text{ , } \Psi(h,m') = m''
\ \Leftrightarrow\ 
(g,m') \in L_{(e,m)} \text{ , } (h,m'')\in L_{(e,m')}
\mapob.
$$
According to \recalt{righttranslpreservesleaves}  $R_g(L_{(e,m')})$ is a leaf-connected component. Moreover, it contains $(g,m')=R_g(e,m')$. Hence $R_g(L_{(e,m')}) = L_{(e,m)}$. But then we have
$$
R_g(h,m'')\in R_g(L_{(e,m')}) = L_{(e,m)}
$$ 
and thus
\begin{multline*}
(hg,m'') \in L_{(e,m)}
\qquad\Longleftrightarrow\qquad
\\
\Psi(hg,m) = m'' = \Psi(h,m') = \Psi\bigl(h,\Psi(g,m)\bigr)
\mapob.
\end{multline*}
As we obviously have $\Psi(e,m) = m$, we have shown that the map $\Psi$ is a (set-theoretic) left-action of $G$ on $M$.

In order to show that $\Psi$ is smooth, we want to apply \recalt{actionwillbesmooth}. 
To do so, we choose for each $m\in M$ the open sets $U_m\subset M$ and $V_m''$ defined in \recalt{defofneighbourhoodbase} and we note that, according to \recalt{localleafcomponents}, for any $\mb\in U_m$ the map
$$
\psi_{e,\mb,V_m''}:V_m''\to G\times M
$$
is a leaf-homeomorphism onto its image $U_{e,\mb,V_m''}$ (with $p$ as its inverse), that this image is leaf-connected and contains $(e,\mb) = \psi_{e,\mb,V_m''}(e)$. Hence $U_{e,\mb,V_m''}$ is contained in the leaf $L_{(e,\mb)}$. 
But by assumption of univalentness, $p$ is a injective when restricted to $L_{(e,\mb)}$. And thus we have the equality
$$
\forall g\in V_m'' : (p\restricted_{U_{e,\mb,V_m''}})\mo(g) = (p\restricted_{L_{(e,\mb)}})\mo(g)
\mapob.
$$
But this implies that for $g\in V_m''$ we have
\begin{align*}
\bigl(g,\Psi(g,\mb)\Bigr) 
&
= (p\restricted_{L_{(e,\mb)}})\mo(g)
=
(p\restricted_{U_{e,\mb,V_m''}})\mo(g)
\\
&
=
\psi_{e,\mb,V_m''}(g)
=
\bigl(g,\Phi_M^M(1,\exp\mo(g), \mb)\bigr)
\mapob.
\end{align*}
It follows that $\Psi$ restricted to $V_m''\times U_m$ is the composition of the smooth map $\exp\mo$ with $\Phi_M^M(1,\cdot)$ (fixing the time parameter to the real value $1$ does not affect smoothness).
Hence $\Psi$ is smooth on $V_m''\times U_m$. This shows that the condition of \recalt{actionwillbesmooth} is satisfied, and thus $\Psi$ is a smooth local left-action of $G$ on $M$.
\end{preuve}

\begin{proclaim}[completehencecovering]{{Proposition}}
If the flow of all smooth vector fields $\rho(X)$, $X\in \body \Liealg g_0$ is complete, then the map $p:G\times M\to G$ is a covering map when we equip $G\times M$ with the leaf topology.

\end{proclaim}

\begin{preuve}
Let $g\in G$ be fixed, then $p\mo(\{g\}) = \{g\}\times M$. We will show that there exist an open neighbourhood $V$ of $g$ and leaf-open sets $U_m\subset G\times M$ for each $m\in M$ such that $p\mo(V) = \cup_{m\in M} U_m$, such that $m\neq m'\Rightarrow U_m\cap U_{m'} = \emptyset$ and such that $p:U_m\to V$ is a leaf-homeomorphism. 
As these properties are exactly the conditions for $p$ to be a covering map (for the leaf topology), we then will have shown our result.

As all vector fields $\rho(X)$, $X\in \body \Liealg g$ are complete, the vector field $Z_A$ is complete \recalt{completenessofZA}, which means in particular that we have, for all $m\in M$, the inclusion $\{1\}\times \Domexp\times \{(g,m)\}\subset W_A = \CA_0\times \Liealg g_0\times (G\times M)$ (where $\Domexp=\exp\mo(\Imexp)$ is the fixed neighbourhood of $0\in \Liealg g_0$ on which the exponential map is a diffeomorphism).
It follows that all maps $\psi_{g,m,\Imexp}:\Imexp\cdot g\to G\times M$ \recalt{defoflocalleafcomponents} are defined.
We now claim that the open neighbourhood $V=\Imexp\cdot g$ of $g$ together with the sets $U_m=\psi_{g,m,\Imexp}(V) \equiv U_{g,m,\Imexp}$ \recalt{localleafcomponents} satisfy our requirements.

According to \recalt{localleafcomponents} $U_{g,m,\Imexp}$ is leaf-open and $p:U_{g,m,\Imexp}\to V$ is a leaf-homeo\-mor\-phism.
It thus remains to show that we have $p\mo(V) = \cup_{m\in M} U_{g,m,\Imexp}$ and the implication $m\neq m'\Rightarrow U_{g,m,\Imexp} \cap U_{g,m',\Imexp}=\emptyset$. Both these properties follow from the group property of a flow. For the first, choose $x\in p\mo(V)$ and define $h=p(x)\in V$. 
This means that $hg\mo\in \Imexp$ and thus $X=\exp\mo(hg\mo)$ is well defined.
With this $X$ we compute, using \recalf{expressionforPhisubA} and the fact that the flow is complete,
$$
p_2\bigl(\Phi_A(-1,X,x)\bigr)
=
\exp(-1\cdot X)h
=
(hg\mo)\mo h
=
g
\mapob.
$$
It follows that there exists $m\in M$ such that
$$
\Phi_A(-1,X,x) = (X,g,m)
\mapob.
$$
But then we can use the group property to compute:
\begin{align*}
(X,x)
&
=
\Phi_A\bigl(1,\Phi_A(-1, X,x)\bigr)
=
\Phi_A(1,X,g,m)
\\
&
=
\Phi_A\bigl(1,\exp\mo(hg\mo), g,m\bigr)
\mapob,
\end{align*}
from which it follows immediately (by applying $p_{23}$ to both sides) that we have $x=\psi_{g,m,\Imexp}(h) \in U_{g,m,\Imexp}$. Hence we have the equality $p\mo(V) = \cup_{m\in M} U_{g,m,\Imexp}$. 

For the last property, suppose $x\in \psi_{g,m,\Imexp}(V) \cap \psi_{g,m',\Imexp}(V)$ and define $h=p(x)$. By definition of $\psi_{g,m,\Imexp}$ and $\psi_{g,m',\Imexp}$, the equality $\psi_{g,m,\Imexp}(h) = x= \psi_{g,m',\Imexp}(h)$ implies that we (also) have the equality
$$
\Phi_A\bigl(1,\exp\mo(hg\mo),g,m\bigr) = \Phi_A\bigl(1,\exp\mo(hg\mo),g,m'\bigr)
\mapob.
$$
Applying $\Phi_A(-1, \cdot)$ to both sides (and using the group property of a flow) then tells us that we have the equality
$$
\bigl(\exp\mo(hg\mo),g,m\bigr) = \bigl(\exp\mo(hg\mo),g,m'\bigr)
\mapob,
$$
and thus in particular $m=m'$.
\end{preuve}

\begin{proclaim}[integratinginfaction]{{Theorem}}
Let $G$ be a connected and simply connected Lie supergroup with Lie superalgebra $\Liealg g$, let $M$ be a smooth supermanifold and let $\rho$ be an infinitesimal action of $\Liealg g$ on $M$. If all vector fields $\rho(X)$, $X\in \body\Liealg g_0$ are complete, then there exists a unique smooth (global) left-action of $G$ on $M$ integrating $\rho$.

\end{proclaim}

\begin{preuve}
According to \recalt{completehencecovering} the map $p:G\times M\to G$ is a covering map when $G\times M$ is equipped with the leaf topology.
Hence the restriction of $p$ to any leaf-connected component of $G\times M$ also is a covering map. 
But $G$ is simply connected, so any connected covering must be a bijection.
This implies that $(G,\rho)$ is univalent, but at the same time that we have $p(L_{(e,m)}) = G$.
According to \recalt{univalentthusmaximallocalint} we thus have a unique maximal local action integrating $\rho$ defined on $D=G\times M$, \ie, a global action.
\end{preuve}

\section{Examples}

The simplest examples of Lie superalgebras are the ones of dimension $1\vert 0$ and $0\vert 1$. The first one has a single even basis vector and the second a single odd one.
An infinitesimal action of a Lie superalgebra of dimension $1\vert 0$ thus is completely determined by a single {even} smooth vector field $X$ on a supermanifold $M$.
A local action of the corresponding Lie supergroup $\CA_0$ (with addition as operation) is nothing more nor less than a (local) flow of the generating vector field $X$.
And in this context, a global action is synonym to a complete vector field.

The situations is barely more complicated for an infinitesimal action of a Lie superalgebra of dimension $0\vert1$ with the single odd basis vector $f_1$. It is completely determined by a single odd smooth vector field $X=\rho(f_1)$ on a supermanifold $M$, but with the additional condition that the autocommutator $[X,X]$ is zero. 
This condition (which is automatically satisfied for even vector fields) is needed because any Lie superalgebra of dimension $0\vert n$ is abelian, and thus we must have $0=\rho(0)=\rho([f_1,f_1]) = [\rho(f_1),\rho(f_1)] = [X,X]$.
As $\body X$ is the zero vector field, it has (trivially) a complete flow, and thus the infinitesimal action integrates to a global action of the corresponding Lie supergroup $\CA_1$ (again with addition as operation). And indeed, any smooth odd vector field $X$ satisfying $[X,X]=0$ admits a global flow with an odd \quote{time} parameter $\tau\in \CA_1$ \cite[V.4.15-18]{Tu04}.

In \cite{MoSV93} a more general theory of integrating vector fields is developped, in which any (not necessarily homogeneous) vector field $X=X_0+X_1$ (with $X_0$ even and $X_1$ odd) can be integrated. However, no group-like properties exist for this more general case, except for $3$ exceptional cases: 
\begin{enumerate}
\item
$[X_0,X_1]=[X_1,X_1]=0$, 

\item
$[X_0,X_1]=aX_1$ and $[X_1,X_1]=0$ with $a\in \RR\setminus\{0\}$,

\item\label{specialcaseiiiof11superalgebra}
$[X_0,X_1]=0$ and $[X_1,X_1]=aX_0$ with $a\in \RR\setminus\{0\}$.

\end{enumerate}
But these three possibilities correspond exactly to the three types of Lie superalgebras of dimension $1\vert1$. 
It thus seems natural to consider these three cases, not as a single non-homogeneous vector field on $M$, but as an infinitesimal action of a $1\vert1$-dimensional Lie superalgebra $\Liealg g$: if $f_1$ is the even basis vector and $f_2$ the odd basis vector, then we write $\rho(f_1) = X_0$ and $\rho(f_2)=X_1$. The commutation relations between the vector fields then garantee that $\rho$ is an infinitesimal action of $\Liealg g$ (with the corresponding structure of a Lie superalgebra of course).
In order to illustrate the general theory of this paper in a not completely trivial example, we consider the case (\ref{specialcaseiiiof11superalgebra}).

\begin{definition}{The Lie supergroup}
We consider the Lie supergroup $G=\CA_0\times \CA_1$ with the even coordinate $a\in \CA_0$, the odd coordinate $\alpha\in \CA_1$, and  with multiplication defined by
$$
(a,\alpha)\cdot (b,\beta) = (a+b+2\alpha\beta, \alpha+\beta)
\mapob.
$$
The tangent vectors $\partial_a, \partial_\alpha$ at $e=(0,0)$ extend to left-invariant vector fields $f_1,f_2$ on $G$ as
$$
f_1\caprestricted_{(a,\alpha)}=\fracp{}{a}\bigrestricted_{(a,\alpha)}
\qquad\text{and}\qquad
f_2\caprestricted_{(a,\alpha)}=\Bigl(\fracp{}{\alpha} -2\alpha\fracp{}{a}\Bigr)\bigrestricted_{(a,\alpha)}
\mapob.
$$
It follows immediately that we have the commutation relations $[f_1,f_1]=0=[f_1,f_2]$ and $[f_2,f_2]=-4f_1$.

\end{definition}

\begin{definition}{The infinitesimal action}
Next we consider the supermanifold $M=(\CA_0)^2\times (\CA_1)^2$ of dimension $2\vert 2$ with global coordinates $(x,y,\xi,\eta)$ with $x,y$ even and $\xi,\eta$ odd.
On $M$ we define the even vector field $\rho(f_1)$ and the odd vector field $\rho(f_2)$ by
\begin{align*}
\rho(f_1)
&=
y\,\fracp{}{x} - x\,\fracp{}{y} + \eta\,\fracp{}{\xi}  - \xi\, \fracp{}{\eta}
\\
\rho(f_2)
&=
(\eta-\xi)\,\fracp{}{x} - (\eta+\xi)\,\fracp{}{y} + (y-x)\,\fracp{}{\xi} - (x+y)\,\fracp{}{\eta}
\mapob.
\end{align*}
It is elementary to check that these vector fields have the following commutation relations:
$$
[\rho(f_1),\rho(f_1)] = 0
\quad,\quad
[\rho(f_1),\rho(f_2)]=0
\quad\text{and}\quad
[\rho(f_2),\rho(f_2)] = -4\rho(f_1)
\mapob,
$$
and thus $\rho$ can be interpreted as an infinitesimal action of $\Liealg g$ on $M$.

\end{definition}

When we look at $\body \rho(f_1)$ and $\body\rho(f_2)$, we find $\body\rho(f_2)=0$ and 
$$
\body \rho(f_1)
=
x\,\fracp{}{y} - y\,\fracp{}{x}
\mapob.
$$
The first is trivially complete, and the second has the complete flow given by
$$
\Phi(t,x,y) = \bigl( \,x\cos t - y \sin t\ ,\ x\sin t + y\cos t\,\bigr)
\mapob.
$$
According to \recalt{integratinginfaction}, this infinitesimal action thus integrates to a (unique) global action of $G$ on $M$.
To compute this action, we can either use \recalt{existenceoflocalactions}, which requires the integration of a vector field on a $6$-dimensional space, or we can use (the proof of) \recalt{univalentthusmaximallocalint}, which requires the determination of the leaves of a $2$-dimensional foliation in a $6$-dimensional space. We will provide both computations, so the reader can judge which one (s)he prefers.

\begin{definition}{The action via the flow of $Z_M$}
When we want to use \recalt{existenceoflocalactions}, we not only have to compute the flow of the vector field $Z_M$, but also the exponential map $\exp:\Liealg g_0\to G$, for which we have to integrate the vector field $Z_R$.
As the right-invariant vector fields associated to the tangent vectors $\partial_a$ and $\partial_\alpha$ are given by
$$
f_1^r\caprestricted_{(a,\alpha)}=\fracp{}{a}\bigrestricted_{(a,\alpha)}
\qquad\text{and}\qquad
f_2^r\caprestricted_{(a,\alpha)}=\Bigl(\fracp{}{\alpha} +2\alpha\fracp{}{a}\Bigr)\bigrestricted_{(a,\alpha)}
\mapob,
$$
it follows that the vector field $Z_R$ on $\Liealg g_0\times G\cong (\CA_0\times \CA_1) \times (\CA_0\times \CA_1)$ with coordinates $(z,\zeta, a,\alpha)$ is given by
$$
Z_R\caprestricted_{(z,\zeta,a,\alpha)}
=
z\cdot f_1^r + \zeta\cdot f_2^r
=
(z+2\zeta\alpha) \fracp{}{a} + \zeta\fracp{}{\alpha}
\mapob.
$$
As we know that the flow $\Phi_R$ does not move the algebra elements, we can write it in terms of these coordinates as
$$
\Phi_R(t,z,\zeta,a,\alpha) = \bigl(z,\zeta, \Phi_R^a(t,z,\zeta,a,\alpha), \Phi_R^\alpha(t,z,\zeta,a,\alpha)\bigr)
\mapob.
$$
In terms of these coefficients, the flow equation $\contrf{\partial_t\caprestricted_{(t,z,\zeta,a,\alpha)}}{\Phi_R} = Z_R\caprestricted_{\Phi_R(t,z,\zeta,a,\alpha)}$ becomes the two equations
$$
\fracp{\Phi_R^a}{t}(t,z,\zeta,a,\alpha)
=
z + 2\zeta\,\Phi_R^\alpha(t,z,\zeta,a,\alpha)
\quad\text{and}\quad
\fracp{\Phi_R^\alpha}{t}(t,z,\zeta,a,\alpha)
=
\zeta
\mapob.
$$
Developping the coefficients $\Phi_R^a$ and $\Phi_R^\alpha$ with respect to the odd coordinates as
\begin{align*}
\Phi_R^a(t,z,\zeta,a,\alpha)
&
=
\Phi_{R,0}^a(t,z,a) + \zeta\,\alpha\,\Phi_{R,2}^{a}(t,z,a)
\\
\Phi_R^\alpha(t,z,\zeta,a,\alpha)
&
=
\zeta\,\Phi_{R,\zeta}^{\alpha}(t,z,a) + \alpha\,\Phi_{R,\alpha}^\alpha(t,z,a)
\mapob,
\end{align*} 
these two equations translate as the system of four ordinary differential equations (in real coordinates)
\begin{align*}
\fracp{\Phi_{R,0}^a}{t}(t,z,a)
&
=
z 
&
\fracp{\Phi_{R,2}^a}{t}(t,z,a)
&
=
2\,\Phi_{R,\alpha}^\alpha(t,z,a)
\\
\fracp{\Phi_{R,\alpha}^\alpha}{t}(t,z,a)
&
=
0
&
\fracp{\Phi_{R,\zeta}^\alpha}{t}(t,z,a)
&
=
1
\end{align*}
Now the initial conditions for this flow are 
$$
\Phi_R^a(0,z,\zeta,a,\alpha)= a
\qquad\text{and}\qquad
\Phi_R^\alpha(0,z,\zeta,a,\alpha)=\alpha
\mapob,
$$
which translates as the four initial conditions
\begin{align*}
\Phi_{R,0}^a(0,z,a) 
&
=
a
&
\Phi_{R,\alpha}^\alpha(0,z,a) 
&
= 1
\\
\Phi_{R,\zeta}^\alpha(0,z,a)
&
=0 
&
\Phi_{R,2}^a(0,z,a) 
&
=0
\mapob.
\end{align*}
Taking these initial conditions into account, we obtain the solutions
\begin{align*}
\Phi_{R,0}^a(t,z,a) 
&
= a+tz
&
\Phi_{R,\alpha}^\alpha(t,z,a) 
&
= 
1
\\
\Phi_{R,\zeta}^\alpha(t,z,a) 
&
= t
&
\Phi_{R,2}^a(t,z,a) 
&
= 2t
\mapob,
\end{align*}
which gives us the flow as
$$
\Phi_R(t,z,\zeta,a,\alpha)
=
\bigl(z,\zeta, a+t(z+2\zeta\alpha), \alpha+t\zeta \bigr)
\mapob,
$$
and finally the exponential map as
$$
\bigl( (z,\zeta), \exp(z,\zeta)\bigr)
=
\Phi_R(1,(z,\zeta),(0,0)\bigr)
=
\bigl( (z,\zeta), (z,\zeta)\bigr)
\mapob,
$$
\ie,
\begin{equation}\label{exponentialinexplicitexamp}
\exp(z\,f_1+\zeta\,f_2)
\cong
\exp(z,\zeta) = (z,\zeta)
\mapob.
\end{equation}

The next step is to integrate the vector field $Z_M$ on $\Liealg g_0\times M$ given by
\begin{align*}
Z_M
&
=
-z\,\rho(f_1) - \zeta\,\rho(f_2)
\\
&
=
\bigl(\zeta(\xi-\eta)-zy\bigr)\fracp{}{x}  + \bigl(zx + \zeta(\eta+\xi)\bigr) \fracp{}{y}  
\\
&
\kern7em
+\bigl(\zeta(x-y) - z\eta\bigr) \fracp{}{\xi} + \bigl(\zeta(x+y) + z\xi\bigr)\fracp{}{\eta}
\mapob.
\end{align*}
As for the flow of $Z_R$, we know that $\Phi_M$ does not move the algebra elements, so we can write it in terms of the coordinates as
$$
\Phi_M(t,z,\zeta,x,y,\xi,\eta)
=
(\, z,\zeta, \Phi_M^x, \Phi_M^y, \Phi_M^\xi, \Phi_M^\eta \,)
\mapob,
$$
with $\Phi_M^x, \Phi_M^y, \Phi_M^\xi, \Phi_M^\eta$ four functions depending upon all the coordinates $t$, $z$, $\zeta$, $x$, $y$, $\xi$, $\eta$.
In terms of these functions, the flow equation $\contrf{\partial_t}{T\Phi_M} = Z_M$ becomes the system of four equations (not writing the dependence on the coordinates $t$, $z$, $\zeta$, $x$, $y$, $\xi$, $\eta$)
\begin{align*}
\fracp{\Phi_M^x}{t}
&
=
\zeta\,(\, \Phi_M^\xi - \Phi_M^\eta\,) - z\,\Phi_M^y
&
\fracp{\Phi_M^y}{t}
&
=
z\,\Phi_M^x + \zeta\,(\,\Phi_M^\eta+\Phi_M^\xi\,)
\\
\fracp{\Phi_M^\xi}{t}
&
=
\zeta\,(\, \Phi_M^x - \Phi_M^y\,) - z\,\Phi_M^\eta
&
\fracp{\Phi_M^\eta}{t}
&
=
\zeta\,(\, \Phi_M^x+\Phi_M^y\,) + z\,\Phi_M^\xi
\mapob.
\end{align*}
We then develop these four functions with respect to their odd coordinates as
\begin{align*}
\Phi_M^x
&
=
\Phi_{M,0}^x + \xi\eta\,\Phi_{M,\zeta}^x + \eta\zeta\,\Phi_{M,\xi}^x + \zeta\xi\,\Phi_{M,\eta}^x
\\
\Phi_M^y
&
=
\Phi_{M,0}^y + \xi\eta\,\Phi_{M,\zeta}^y + \eta\zeta\,\Phi_{M,\xi}^y + \zeta\xi\,\Phi_{M,\eta}^y
\\
\Phi_M^\xi
&
=
\xi\,\Phi_{M,\xi}^\xi + \eta\,\Phi_{M,\eta}^\xi + \zeta\,\Phi_{M,\zeta}^\xi + \xi\eta\zeta\,\Phi_{M,3}^\xi
\\
\Phi_M^\eta
&
=
\xi\,\Phi_{M,\xi}^\eta + \eta\,\Phi_{M,\eta}^\eta + \zeta\,\Phi_{M,\zeta}^\eta + \xi\eta\zeta\,\Phi_{M,3}^\eta
\mapob,
\end{align*}
where the sixteen functions $\Phi_{M,-}^{-}$ are functions of the four even coordinates $t,z,x,y$ only.
In terms of these functions, the system of four equations becomes the system of sixteen ordinary differential equations
\begin{align*}
\fracp{\Phi_{M,0}^x}{t}
&
=
- z\,\Phi_{M,0}^y
&
\fracp{\Phi_{M,\xi}^x}{t}
&
=
\Phi_{M,\eta}^\eta - \Phi_{M,\eta}^\xi - z\,\Phi_{M,\xi}^y
\\
\fracp{\Phi_{M,\zeta}^x}{t}
&
=
- z\,\Phi_{M,\zeta}^y
&
\fracp{\Phi_{M,\eta}^x}{t}
&
=
\Phi_{M,\xi}^\xi - \Phi_{M,\xi}^\eta - z\,\Phi_{M,\eta}^y
\end{align*}
\begin{align*} 
\fracp{\Phi_{M,0}^y}{t}
&
=
z\,\Phi_{M,0}^x
&
\fracp{\Phi_{M,\xi}^y}{t}
&
=
- \Phi_{M,\eta}^\eta - \Phi_{M,\eta}^\xi + z\,\Phi_{M,\xi}^x
\\
\fracp{\Phi_{M,\zeta}^y}{t}
&
=
z\,\Phi_{M,\zeta}^x
&
\fracp{\Phi_{M,\eta}^y}{t}
&
=
\Phi_{M,\xi}^\eta + \Phi_{M,\xi}^\xi + z\,\Phi_{M,\eta}^x
\end{align*}
\begin{align*} 
\fracp{\Phi_{M,\xi}^\xi}{t}
&
=
-z\,\Phi_{M,\xi}^\eta
&
\fracp{\Phi_{M,\zeta}^\xi}{t}
&
=
\Phi_{M,0}^x - \Phi_{M,0}^y - z\,\Phi_{M,\zeta}^\eta
\\
\fracp{\Phi_{M,\eta}^\xi}{t}
&
=
-z\,\Phi_{M,\eta}^\eta
&
\fracp{\Phi_{M,3}^\xi}{t}
&
=
\Phi_{M,\zeta}^x - \Phi_{M,\zeta}^y - z\,\Phi_{M,3}^\eta
\end{align*}
\begin{align*} 
\fracp{\Phi_{M,\xi}^\eta}{t}
&
=
z\,\Phi_{M,\xi}^\xi
&
\fracp{\Phi_{M,\zeta}^\eta}{t}
&
=
\Phi_{M,0}^x + \Phi_{M,0}^y + z\,\Phi_{M,\zeta}^\xi
\\
\fracp{\Phi_{M,\eta}^\eta}{t}
&
=
z\,\Phi_{M,\eta}^\xi
&
\fracp{\Phi_{M,3}^\eta}{t}
&
=
\Phi_{M,\zeta}^x + \Phi_{M,\zeta}^y + z\,\Phi_{M,3}^\xi
\end{align*}
The initial conditions at $t=0$ are $\Phi_M^x = x$, $\Phi_M^y = y$, $\Phi_M^\xi = \xi$ and $\Phi_M^\eta=\eta$. In terms of the decomposition this gives us the initial conditions (all functions of $(t,z,x,y)$ taken at $t=0$)
\begin{align*}
\Phi_{M,0}^x &= x
&
\Phi_{M,\zeta}^x &= 0
&
\Phi_{M,\xi}^x &= 0
&
\Phi_{M,\eta}^x &= 0
\\
\Phi_{M,0}^y &= y
&
\Phi_{M,\zeta}^y &= 0
&
\Phi_{M,\xi}^y &= 0
&
\Phi_{M,\eta}^y &= 0
\\
\Phi_{M,\xi}^\xi &= 1
&
\Phi_{M,\eta}^\xi &= 0
&
\Phi_{M,\zeta}^\xi &= 0
&
\Phi_{M,3}^\xi &= 0
\\
\Phi_{M,\xi}^\eta &= 0
&
\Phi_{M,\eta}^\eta &= 1
&
\Phi_{M,\zeta}^\eta &= 0
&
\Phi_{M,3}^\eta &= 0
\mapob.
\end{align*}
Taking these initial conditions into account, the differential equations can be solved successively to yield (more or less in this order)
\begin{align*}
\Phi_{M,0}^x &= x\,\cos tz - y\,\sin tz
&
\Phi_{M,0}^y &= x\,\sin tz + y \,\cos tz
\\
\Phi_{M,\zeta}^x &= 0
&
\Phi_{M,\zeta}^y &= 0
\\
\Phi_{M,\xi}^\xi &= \cos tz
&
\Phi_{M,\xi}^\eta &= \sin tz
\\
\Phi_{M,\eta}^\xi &= -\sin tz
&
\Phi_{M,\eta}^\eta &= \cos tz
\\
\Phi_{M,\xi}^x &= t\,(\cos zt + \sin zt)
&
\Phi_{M,\xi}^y &= -t\,(\cos zt -\sin zt)
\\
\Phi_{M,\eta}^x &= -t\,(\sin tz - \cos tz)
&
\Phi_{M,\eta}^y &= t\,(\sin tz + \cos tz)
\\
\Phi_{M,\zeta}^\xi &= t\,(x-y)\,\cos tz - t\,(x+y)\,\sin tz
&
\Phi_{M,3}^\xi &= 0
\\
\Phi_{M,\zeta}^\eta &= t\,(x-y)\,\sin tz + t\,(x+y)\,\cos tz
&
\Phi_{M,3}^\eta &= 0
\mapob.
\end{align*}
Putting these functions back together, we find
\begin{align*}
\Phi_M^x
&=
x\,\cos tz - y\,\sin tz
+\eta\zeta t\,(\cos zt + \sin zt)
-\zeta\xi t\,(\sin tz - \cos tz)
\\
\Phi_M^y
&=
x\,\sin tz + y \,\cos tz
- \eta\zeta t\,(\cos zt -\sin zt)
+\zeta\xi t\,(\sin tz + \cos tz)
\\
\Phi_M^\xi
&=
\xi \cos tz - \eta \sin tz + t\zeta \bigl((x-y)\,\cos tz - (x+y)\,\sin tz\bigr)
\\
\Phi_M^\eta
&=
\xi \sin tz + \eta \cos tz + t\zeta\bigl((x-y)\,\sin tz + (x+y)\,\cos tz\bigr)
\end{align*}
When we write this in the suggestive form
$$
\begin{pmatrix}
\Phi_M^x \\ \Phi_M^y
\end{pmatrix}
=
\begin{pmatrix}
\cos tz & -\sin tz \\ \sin tz & \cos tz
\end{pmatrix}
\cdot
\left[
\begin{pmatrix} x \\ y \end{pmatrix}
+
t\zeta
\begin{pmatrix} 1&-1\\ 1&1 \end{pmatrix}
\cdot 
\begin{pmatrix}  \xi \\ \eta \end{pmatrix}
\right]
$$
$$
\begin{pmatrix}
\Phi_M^\xi \\ \Phi_M^\eta
\end{pmatrix}
=
\begin{pmatrix}
\cos tz & -\sin tz \\ \sin tz & \cos tz
\end{pmatrix}
\cdot
\left[
\begin{pmatrix} \xi \\ \eta \end{pmatrix}
+
t\zeta
\begin{pmatrix} 1&-1 \\ 1&1 \end{pmatrix}
\cdot
\begin{pmatrix} x \\ y \end{pmatrix}
\right]
\mapob,
$$
it becomes obvious that combining the variables in complex ones as $x+iy$ and $\xi + i\eta$ gives us the much simpler form
$$
\Phi_M
:
\begin{pmatrix} x+iy \\ \xi+i\eta \end{pmatrix}
\mapsto
\eexp^{itz}\cdot \begin{pmatrix} 1 & (1+i)t\zeta \\ (1+i)t\zeta & 1 \end{pmatrix}
\cdot
\begin{pmatrix}x+iy \\ \xi+i\eta \end{pmatrix}
\mapob.
$$
And then it is more or less obvious that we can combine even and odd variables together to write this as the single multiplication
$$
\Phi_M
:
(x+\xi)+i(y+\eta) \mapsto \eexp^{itz}\bigl(1+(1+i)t\zeta\bigr)\bigl( (x+\xi)+i(y+\eta) \bigr)
\mapob.
$$
We then can use the formula $\Psi(g,m) = \Phi_M^M\bigl(1,\exp\mo(g),m)$ together with the formula \recalt{exponentialinexplicitexamp} for the exponential map to find the explicit expression for the global action $\Psi$ as\footnote{Actually, this example was constructed starting with this formula, which is based upon \cite[VI.4.14]{Tu04}}
$$
\Psi\bigl((a,\alpha), (x+\xi)+i(y+\eta) \bigr)
=
\eexp^{ita}\bigl(1+(1+i)\alpha\bigr)\bigl( (x+\xi)+i(y+\eta) \bigr)
\mapob.
$$

\end{definition}

\begin{definition}{The action via the foliation $\foliation$}
The foliation $\foliation$ on $G\times M$ is generated by the two vector fields $f_1^r - \rho(f_1)$ and $f_2^r - \rho(f_2)$, \ie, by
\begin{gather*}
X=\fracp{}{a} - y\,\fracp{}{x} + x\,\fracp{}{y} - \eta\,\fracp{}{\xi}  + \xi\, \fracp{}{\eta}
\\
\noalign{\noindent and}
Y=\fracp{}{\alpha} +2\alpha\fracp{}{a}
+
(\xi-\eta)\,\fracp{}{x} + (\eta+\xi)\,\fracp{}{y} + (x-y)\,\fracp{}{\xi} + (x+y)\,\fracp{}{\eta}
\end{gather*}
Now in general it is (very) hard to find the explicit coordinates adapted to the foliation, but in this case the (local!) procedure suggested by the abstract proof (change the generating vector fields by an invertible matrix function to obtain commuting vector fields, integrate these vector fields and use their flow to create a new coordinate system adapted to the foliation) works quite well and even gives global results.

The first step thus is to abelianize the vector fields $X$ and $Y$, which here can be done globally by posing $\Yt=Y - 2\alpha\, X$:
\begin{multline*}
\quad
\Yt
=
\fracp{}{\alpha} 
+
(\xi-\eta+2y\alpha)\,\fracp{}{x} 
+ (\eta+\xi-2x\alpha)\,\fracp{}{y} 
\\
+ (x-y+2\alpha\eta)\,\fracp{}{\xi} + (x+y-2\alpha\xi)\,\fracp{}{\eta}
\mapob.
\quad
\end{multline*}
For the vector fields $X$ and $\Yt$ (which still generate $\foliation$ everywhere) we have $[X,X] = [X,\Yt] = [\Yt,\Yt] = 0$. Both vector fields thus can be integrated (using the same techniques as before) and their global flows are given by
\begin{align*}
\Phi_X(t,a,\alpha,x,\xi,y,\eta)
&
=
(a+t, \alpha, x\,\cos t - y\,\sin t, x\,\sin t + y \,\cos t , 
\\
&
\kern3em
\xi\,\cos t - \eta \,\sin t, \xi \,\sin t + \eta\,\cos t)
\\
\Phi_\Yt(\tau,a,\alpha,x,\xi,y,\eta)
&
=
\bigl(a, \alpha+\tau, x+\tau(\xi - \eta +2y\alpha), \xi + \tau(x-y+2\alpha\eta) , 
\\
&
\kern3em
y+\tau(\eta+\xi-2x\alpha), \eta+\tau(x+y - 2\alpha\xi)\bigr)
\end{align*}
The final step is to use their flow to introduce the change of coordinates 
$$
(a,\alpha,x,y,\xi,\eta) = \psi(a',\alpha',x',y',\xi',\eta')
$$ 
with $\psi$ defined by
\begin{align*}
\psi(a',\alpha', 
&
x', y',\xi',\eta')
=
\Phi_X\bigl(a',\Phi_\Yt(\alpha',0,0,x',y',\xi',\eta')\bigr)
\\
\noalign{\vskip1\jot}
&
=
\Phi_X\bigl(a',0, \alpha', x'+\alpha'(\xi' - \eta'), y'+\alpha'(\eta'+\xi') , 
\\
&
\kern10em
\xi' + \alpha'(x'-y'), \eta'+\alpha'(x'+y')\bigr)
\\
\noalign{\vskip1\jot}
&
=
(a', \alpha', (x'+\alpha'(\xi' - \eta'))\,\cos a' - (y'+\alpha'(\eta'+\xi'))\,\sin a', 
\\
&
\kern5em
(x'+\alpha'(\xi' - \eta'))\,\sin a' + (y'+\alpha'(\eta'+\xi')) \,\cos a' , 
\\
&
\kern5em
(\xi' + \alpha'(x'-y'))\,\cos a' - (\eta'+\alpha'(x'+y')) \,\sin a', 
\\
&
\kern6em
(\xi' + \alpha'(x'-y')) \,\sin a' + (\eta'+\alpha'(x'+y'))\,\cos a')
\end{align*}
When we remember that it could be useful to regroup all coordinates in a single complex one, we note that we can simplify this expression as
$$
\psi\bigl((a',\alpha'), x'+iy'+\xi'+i\eta'\bigr)
=
\bigl((a',\alpha'), \eexp^{ia'}\,\bigl(1+(1+i)\alpha'\bigr)(x'+iy'+\xi'+i\eta')\bigr)
\mapob.
$$
This allows for an easy inversion to find the new coordinates as
\begin{align*}
\bigl((a',\alpha'),x'+iy'+\xi'&+i\eta'\bigr)
=
\psi\mo\bigl((a,\alpha),x+iy+\xi+i\eta\bigr) 
\\
&
= 
\bigl((a,\alpha),\eexp^{-ia}\bigl(1-(1+i)\alpha\bigr)(x+iy+\xi+i\eta)\bigr) 
\mapob.
\end{align*}
It is elementary to show that in these new coordinates, the vector fields $X$ and $\Yt$ are given as $\partial_{a'}$ and $\partial_{\alpha'}$.
The primed coordinates are thus adapted to the foliation and, moreover, they are global coordinates with the first two being global coordinates on $G$, providing an illlustration of \recalt{locadaptedcoordinates}.

According to (the proof of) \recalt{univalentthusmaximallocalint}, the (global) action is defined by
$$
\Psi(g,m) = \mb
\qquad\Longleftrightarrow\qquad
(g,\mb) \in L_{(e,m)}
\mapob.
$$
But lying on the same leaf means having the same primed coordinates $(x',y',\xi',\eta')$, as these are adapted to the foliation.
As we have
$$
\psi\mo(e,m) = \psi\mo\bigl((1,0), x+iy+\xi+i\eta)
=
\bigl((1,0), x+iy+\xi+i\eta)
\mapob,
$$
it follows that we must have, for $g=(a,\alpha)$,
$$
\psi\mo(g,\mb)=\bigl((a,\alpha), x+iy+\xi+i\eta)
$$
and thus
$$
(g,\mb) =
\psi\bigl((a,\alpha), x+iy+\xi+i\eta\bigr) 
=
\bigl((a,\alpha), \eexp^{ia}\bigl(1 + (1+i)\alpha\bigr)(x+iy+\xi+i\eta) \bigr)
\mapob,
$$
giving the action (again) as
$$
\Psi\bigl((a,\alpha), x+iy+\xi+i\eta\bigr)
=
\eexp^{ia}\bigl(1 + (1+i)\alpha\bigr)(x+iy+\xi+i\eta)
\mapob.
$$

\end{definition}

\section{Appendix: an overview of \texorpdfstring{$\CA$}{A}-manifold theory}

This appendix is a very short overview of $\CA$-manifold theory, intended for readers with some familiarity with supermanifold theory, but not with this approach. It will not always be complete, and sometimes it might be slightly besides the truth, but it intends to give the gist of the theory, not the most precise formulation (which can be found in \cite{Tu04}).

\subsection{The basic graded commutative ring}
\vrule width0pt
\smallskip

The starting point of $\CA$-manifold theory is a graded-commutative ring $\CA$. We will fix it as being the exterior algebra of an infinite dimensional (real) vector space $V$: 
$$
\CA=\bigwedge V = \bigoplus_{k=0}^\infty \bigwedge{}{\kern-3pt\raise1.5ex\hbox{$\scriptstyle k$}\kern2pt} V =  \Bigl(\, \bigoplus_{k=0}^\infty \bigwedge {\kern-3pt\raise1.5ex\hbox{$\scriptstyle 2k$}\kern2pt} V \,\Bigr) \oplus \Bigl(\, \bigoplus_{k=0}^\infty \bigwedge{\kern-3pt\raise1.5ex\hbox{$\scriptstyle 2k+1$}\kern2pt}V \,\Bigr)
=
\CA_0\oplus \CA_1
\mapob.
$$
We denote by $\body:\CA\to \CA$ (sic) the canonical projection onto the direct summand $\RR\equiv\bigwedge^0V\subset \bigoplus_{k=0}^\infty\bigwedge^kV$; we will call $\body$ the \stresd{body map} and the image $\body(a)\in \RR$ the \stresd{body of $a\in \CA$}.
The set of nilpotent elements $\nilpotent=\bigoplus_{k=1}^\infty \bigwedge^kV$ (an ideal) is a supplement (in $\CA$) to $\RR=\body(\CA)=\bigwedge^0V$. 
We thus can identify the quotient $\CA/\nilpotent$ (which is a field) with $\body(\CA)\subset \CA$ and the canonical projection $\CA\mapsto \CA/\nilpotent$ with the body map $\body$.
Another feature of (this) $\CA$ (due to the fact that $V$ is infinite dimensional) is that for any $n\in \NN^*$ there exist elements $\xi_1, \dots, \xi_n\in \CA_1=\bigoplus_{k=0}^\infty \bigwedge^{2k+1}V$ such that the product $\xi_1\cdots\xi_n \equiv \xi_1\wedge\cdots\wedge \xi_n$ is non-zero.

\begin{definition}{Remark}
Actually, the particular form of $\CA$ is not that important, as long as it has the above mentioned features: a supplement to the nilpotent elements which is isomorphic to $\RR$ and for any $n\in \NN^*$ elements $\xi_1,\dots, \xi_n\in \CA_1$ whose product is non-zero. And even the last condition is slightly stronger than strictly needed: if an $\CA$-manifold has odd dimension $q$, then there should exist $q+1$ odd elements whose product is non-zero (and thus $V$ should have at least dimension $q+1$). However, as one does not wish to be restricted in the choice of the odd dimension, it is preferable to have this condition for any $n\in \NN^*$, not only for $q+1$. Hence the choice of an infinite dimensional vector space $V$.
\end{definition}

In the sequel we will never have any use for the particular form of $\CA$, so the symbol $V$ will no longer be reserved for the vector space whose exterior algebra is $\CA$.

\subsection{\texorpdfstring{$\CA$}{A}-vector spaces}
\vrule width0pt
\smallskip

As $\CA$ is not commutative, there is a difference between left- and right-modules over $\CA$. We will first concentrate on graded bi-modules, meaning an abelian group $M$ which is at the same time a left- and a right-module over $\CA$, and which splits as a direct sum $M=M_0\oplus M_1$ satisfying the conditions
$$
a\in \CA_\alpha \text{ and } m\in M_\beta
\quad\Rightarrow\quad
a\cdot m =(-1)^{\alpha\beta}m\cdot a \in M_{\alpha+\beta}
\mapob,
$$
where the gradings $\alpha,\beta$ should be seen as belonging to $\ZZ/2\ZZ$ and thus $1+1=0$.

Again because $\CA$ is not commutative, there is a difference between left-linear homomorphisms and right-linear ones, even for graded bi-modules. 
A left-linear morphism that is even (meaning that it maps homogeneous elements of the source module to homogeneous elements of the target module of the same parity) is automatically also right-linear (and vice-versa), but a similar statement for non-even morphisms is not true. 
One thus has two morphism modules between two graded bi-modules, the left-linear ones and the right-linear ones. 
These morphism modules are in a natural way graded bi-modules over $\CA$ and they are in a natural way isomorphic.

\begin{definition}{Convention}
As we wish to adhere to the convention that interchanging two homogeneous elements (of whatever nature) induces a minus sign whenever both are odd, we are led to denote evaluation of left-linear morphisms on the left instead of on the right as is usual for maps. More precisely, if $f:M\to N$ is a left-linear morphism between the graded bi-modules $M$ and $N$, we will denote the image of the element $m\in M$ by the map $f$ as $\contrf{m}f$. In that way, left-linearity gets the form
$$
\contrf{a\cdot m}f = a\cdot(\contrf{m}f)
\mapob,
$$
instead of
$$
f(am) = af(m)
\mapob,
$$
which would violate the convention in case $f$ is not even.
This notational convention is already used in ordinary differential geometry when one evaluates\slash contracts a $k$-form with a tangent vector (to yield a $k-1$-form).

\end{definition}

The (for us) important constructions on vector spaces can be carried out also for graded bi-modules over $\CA$, such as direct sums, direct products, morphism modules (left or right), tensor products (over the graded-commutative ring $\CA$), exterior algebras and free bi-modules on a set of homogeneous generators. 

\medskip

We next turn our attention to a special case of graded bi-modules over $\CA$: those that are of the form $V\otimes_\RR \CA$ with $V=V_0\oplus V_1$ a graded vector space over $\RR$, which we call \stresd{$\CA$-vector spaces}. 
Inside such an $\CA$-vector space we have the (real) subspace $V\otimes_\RR 1\cong V$ of elements of the form $v\otimes 1$ and the projection (body map) $\body_V:V\otimes_\RR \CA\to V\otimes_\RR 1\cong V$ defined by $\body_V(v\otimes a) = v\otimes \body(a)$.

If $M$ is any (left) module over $\CA$, we can define the subset $\nilpotent_M\subset M$ of those elements $m$ for which there exists $a\in \CA$, $a\neq0$ such that $a\cdot m=0$.
As the reals are included in $\CA$, any module $M$ over $\CA$ is in particular a real vector space and one can show that $\nilpotent_M$ is a vector subspace of $M$ (as a real vector space).
For an $\CA$-vector space $M=V\otimes_\RR \CA$ we have
$$
\nilpotent_M = \{\,v\otimes n\mid v\in V\,,\, n\in\nilpotent\subset \CA\,\}
\mapob.
$$
Another way to define an $\CA$-vector space thus is as those graded bi-modules for which there exists a supplement $V\subset M$ for $\nilpotent_M$ in the category of graded vector spaces over $\RR$.
The body map then becomes the projection onto this summand $M=V\oplus \nilpotent_M$. 
In this way the definition of an $\CA$-vector space obtains the same flavour as the conditions imposed on our graded commutative ring $\CA$.

\begin{definition}{Remark}
The map $V\mapsto V\otimes_\RR \CA$ from the category of graded vector spaces over $\RR$ to the category of graded bi-modules over $\CA$ is a functor (nearly an isomorphism of categories). In particular the constructions one can perform on these categories (direct sums, tensor products, exterior algebras) are preserved by this map. However, one has to make a choice what to do with morphisms, as in the category of graded vector spaces over $\RR$ there is no difference between left- and right-linear, whereas in the category of $\CA$-vector spaces there is.

\end{definition}

\subsection{Smooth functions}
\label{subsecsmoothfunc}
\vrule width0pt
\smallskip

In order to give a definition of smooth super functions from first principles, we first give an alternate description of smooth functions in the ordinary (non-super) case. We then just copy this alternate description to the super case to obtain our definition of super smooth functions.
The starting point is the formula
\begin{equation}\label{fminfequalsintoverder}
f(x) - f(y) = \sum_{i=1}^m (x_i-y_i)\cdot \left( \int_0^1 (\partial_if)\bigl(sx + (1-s)y\bigr)\,\extder s\right)
\mapob,
\end{equation}
valid for any function $f:O\subset \RR^m\to \RR^n$ of class $C^1$ defined on a convex set $O$. If we define the functions $g_i:O\times O\to \RR^n$ by
$$
g_i(x,y) = \int_0^1 (\partial_if)\bigl(sx + (1-s)y\bigr)\,\extder s
\mapob,
$$
then we have:
\begin{equation}\label{smoothnesstrick}
\forall x,y\in O : 
f(x)-f(y) = \sum_{i=1}^m (x_i-y_i)\cdot g_i(x,y)
\mapob.
\end{equation}
It then it is easy to show that, for any $k\in \NN$, $f$ is of class $C^{k+1}$ if and only if there exist functions $g_i$ of class $C^k$ such that \recalf{smoothnesstrick} is valid.
Moreover, the partial derivatives of $f$ are given by $(\partial_if)(x) = g_i(x,x)$. As is well known, these partial derivatives are unique, while the functions $g_i$ themselves are not (apart from the case $n=1$).

\begin{definition}{``Classical'' definitions}
A \stresd{smooth system} is an assigment $\smooths$ that associates to any open set $O\subset \RR^m$ and any target space $\RR^n$ a collection of \stress{continuous} functions $\smooths(O,\RR^n)\subset C^0(O,\RR^n)$ verifying the property
\begin{multline*}
\forall f\in \smooths(O,\RR^n)
\quad
\exists g_i\in \smooths(O^2,\RR^n)
\quad
\forall x,y\in O\subset \RR^m
\quad:\quad
\\
f(x) - f(y) = \sum_{i=1}^m (x_i-y_i)\cdot g_i(x,y)
\mapob.
\end{multline*}
A smooth system $\smooths_1$ is said to be \stresd{smaller than} a smooth system $\smooths_2$, denoted as $\smooths_1\le \smooths_2$, if we have
$$
\forall O,n: \smooths_1(O,\RR^n)\subset \smooths_2(O,\RR^n)
\mapob.
$$

\end{definition}


\begin{proclaim}{Theorem}
$C^\infty$ is the (not a) maximal (with respect to the order $\le$) smooth system.
\end{proclaim}

\begin{definition}{Remark}
The above definitions and statements should not be taken literally, as they are wrong as stated. Formula \recalf{fminfequalsintoverder} is valid for convex sets, but not in general for arbitrary (open) sets. \stress{If} one disposes of smooth partitions of unity, one can re-establish such a result, but that would create, in the case of $\RR$, a circular definition, and in the case of $\CC$ it would be impossible. The solution is to use covers at every stage, but that complicates the notation (not the idea), so we left it out in this summary.

\end{definition}

In order to mimick these definitions in the super case, we have to start with a topology.

\begin{definition}{Definition}
Let $E$ be any finite dimensional $\CA$-vector space, meaning that $\body E\subset E$ is a finite dimensional graded vector space over $\RR$ and $E=\body E \otimes_\RR \CA$. 
The \stresd{DeWitt topology} on $E$ is the coarsest topology on $E$ for which the projection $\body:E\to \body E$ is continuous (when $\body E$ is equipped with the usual euclidean topology).
In particular $\CA$ itself is a finite dimensional $\CA$-vector space with $\body \CA=\RR\oplus\{0\}$ (i.e., $(\body\CA)_0=\RR$ and $(\body\CA)_1=\{0\}$), and thus is equipped with the DeWitt topology.
It should be noted that the DeWitt topology is not separated\slash Hausdorff and that for any open set $U\subset E$ we have $\body U\subset U$ and $U=\body\mo(\body U)$.

\end{definition}

Now let $E$ be a finite dimensional $\CA$-vector space and let $e_1, \dots, e_p$, $f_1, \dots, f_q\in \body E$ be a homogeneous basis, i.e., $e_1, \dots, e_p$ is a basis of the vector space (over $\RR$) $(\body E)_0\subset E_0$ and $f_1, \dots, f_q$ a basis of $(\body E)_1\subset E_1$. It follows that any element $e\in E$ is described uniquely by $p+q$ elements $x_1, \dots, x_p, \xi_1, \dots, \xi_q\in\CA$ according to
$$
e = \sum_{i=1}^p x_ie_i + \sum_{j=1}^q \xi_jf_j
\mapob.
$$
Now if $e$ belongs to the \stress{even} part of $E$, then the $x_i$ belong to $\CA_0$ and the $\xi_j$ to $\CA_1$. In other words, an element of $E_0$ is described by $p$ even ``coordinates'' and $q$ odd ``coordinates.''
Moreover, a subset $U\subset E_0$ of the even part of $E$ is open (in the DeWitt topology) if and only if there exists an open set $O\subset \RR^p$ such that
\begin{equation}\label{defofopeninE0}
e = \sum_{i=1}^p x_ie_i + \sum_{j=1}^q \xi_jf_j \in U
\qquad\Longleftrightarrow\qquad
(\body x_1, \dots, \body x_p) \in O
\mapob.
\end{equation}
Another way to state this equivalence is the equality $U = (\body\caprestricted_{E_0})\mo(\body U)$ with $\body U = O\oplus \{0\}\subset (\body E)_0 \oplus (\body E)_1=\RR^p \oplus \RR^q$.

\begin{definition}[superdefinitions]{``Super'' definitions}
$\bullet$
A \stresd{super smooth system} is an assigment $\smooths$ that associates to any open set $U\subset E_0$ in the even part of any finite dimensional $\CA$-vector space $E$ and to any finite dimensional $\CA$-vector space $F$ (the target space) a collection of \stress{continuous} functions $\smooths(U,F)\subset C^0(U,F)$ verifying the two properties
$$
f(\body U) \subset \body F
$$
and
\begin{multline*}
\forall f\in \smooths(U,F)
\quad
\exists g_i\in \smooths(U^2,F)
\quad
\forall x,y\in U\subset E_0
\quad:\quad
\\
f(x) - f(y) = \sum_{i=1}^{p+q} (x_i-y_i)\cdot g_i(x,y)
\mapob.
\end{multline*}

$\bullet$
A super smooth system $\smooths_1$ is said to be \stresd{smaller than} a super smooth system $\smooths_2$, denoted as $\smooths_1\le \smooths_2$, if we have
$$
\forall U,F: \smooths_1(U,F)\subset \smooths_2(U,F)
\mapob.
$$

$\bullet$
$C^\infty$ is \stresd{the} maximal super smooth system with respect to the order $\le$.

\end{definition}


We now have a definition of super smooth functions that is a look-alike to a possible definition of smooth functions in the non-super case, but two questions remain: (1) why the additional condition $f(\body U)\subset \body F$ and (2) can we say anything interesting about these smooth functions?
One possible answer to (1) is that it allows us to give a positive answer to (2). 
Another possible answer to (1) is that we wish to stay as close as possible to non-super smooth functions and that this conditions is trivially satisfied in the non-super case.

\begin{proclaim}[Gextension]{Lemma}
Let $U\subset E_0$ be an open set in the even part of the $p\vert q$-dimensional $\CA$-vector space $E$, let $O=\body U\subset \RR^p$ be the corresponding open subset in $\RR^p$ (see \recalf{defofopeninE0}) and let $f:O\to \body F$ be an ordinary smooth function with values in the real vector space $\body F$ (where we thus ignore the grading).

With these data, we define the function $\mathbf{G}f:U\to F$ by the formula
$$
(\mathbf{G}f)(x_1, \dots, x_p, \xi_1, \dots, \xi_q)
=
\sum_{k=0}^\infty \frac1{k!} \bigl((D^kf)(r_1, \dots, r_p)\bigr)(\overbrace{n,\vrule width0pt height2.3ex\dots, n}^{k \text{ terms}})
\mapob,
$$
where $r_i=\body x_i\in \RR$, $n_i=x_i-r_i\in \nilpotent_0$, $n=(n_1, \dots, n_p)$, $(D^kf)(r_1, \dots, r_p)$ the $k$th order derivative of $f$ at $(r_1, \dots, r_p)\in O$ as $k$-linear symmetric map with values in $\body F$ and $\bigl((D^kf)(r_1, \dots, r_p)\bigr)({n,\dots, n})$ the formal evaluation of $(D^kf)(r_1, \dots, r_p)$ in the even nilpotent coefficients $n_i$.

This $\mathbf{G}f$ is a super smooth function.

\end{proclaim}

\begin{proclaim}[structureofsmoothfunctions]{Theorem}
Let $U\subset E_0$ be an open set in the even part of the $p\vert q$-dimensional $\CA$-vector space $E$, let $O=\body U\subset \RR^p$ be the corresponding open subset in $\RR^p$ and let $F$ be a finite dimensional $\CA$-vector space.
Then a function $f:U\to F$ belongs to $C^\infty(U,F)$ (i.e., is super smooth) if and only if there exist (ordinary) smooth functions $f_{I}:O\to \body F$ such that we have
$$
f(x_1, \dots, x_p, \xi_1, \dots, \xi_q)
=
\sum_{I\subset \{1, \dots, q\}} (\mathbf{G}f_I)(x_1, \dots, x_p)\cdot \xi^I
\mapob,
$$
where $\xi^I$ denotes the product
$$
I = \{j_1, \dots, j_k\} \text{ with } 1\le j_1<\cdots <j_k\le q
\qquad\Longrightarrow\qquad
\xi^I = \xi_{j_1} \cdots \xi_{j_k}
\mapob,
$$
with $\xi^\emptyset = 1$.

\end{proclaim}

\begin{proclaim}{Corollary}
Let $U\subset E_0$ be an open set in the even part of the $p\vert q$-dimensional $\CA$-vector space $E$, let $O=\body U\subset \RR^p$ be the corresponding open subset in $\RR^p$ and let $F$ be a finite dimensional $\CA$-vector space. Then we have the equality
$$
C^\infty(U,F)
=
C^\infty(O,\body F) \otimes \bigwedge \RR^q
\mapob.
$$

\end{proclaim}

Once we know the structure of super smooth functions, we can address the question of partial derivatives. For that we imitate the non-super case and want to define the partial derivative $\partial_if$ as
$$
(\partial_if)(x) = g_i\bigl(x,x)
\mapob,
$$
where the super smooth functions $g_i$ are given by the definition of super smoothness as
$$
f(x) - f(y) = \sum_{i=1}^{p+q} (x_i-y_i)\cdot g_i(x,y)
\mapob.
$$
But for this to be a coherent definition, one must show that the diagonal $g_i(x,x)$ is uniquely determined by $f$, even when the functions $g_i(x,y)$ are not unique. 
And it is here that one needs the fact that there exist $q+1$ odd elements in $\CA_1$ whose product is non-zero: if that condition is satisfied, the $\partial_if$ are well defined; if not, even the diagonal functions $g_i(x,x)$ will not be unique and thus our super smooth functions will not have derivatives. 
(Of course one could define the partial derivatives by hand in terms of the functions $f_I$ of \recalt{structureofsmoothfunctions}, but then the whole purpose of giving an intrinsic definition of super smoothness would be superfluous.)

\medskip

With our choice of $\CA$, the partial derivatives are well defined and behave exactly as expected. If the index $i$ is associated to an odd coordinate, $\partial_i$ is an odd derivation of $C^\infty(U,F)$, whereas it is an even derivation when associated to an even coordinate. 
One should note that the $\partial_i$ are \stress{right}-derivations (as opposed to left-derivations), meaning that they are right-linear, which is a consequence of the choice in \recalt{superdefinitions} to write the ``coefficients'' $(x_i-y_i)$ to the left of $g_i(x,y)$ instead of to the right.

Once we know what smooth functions and their derivatives are, ``all'' usual results of real analysis can be copied, in particular the inverse function theorem.

\subsection{\texorpdfstring{$\CA$}{A}-manifolds}
\label{subsecAmfds}
\vrule width0pt
\smallskip

Once we know what open subset are and what (local) diffeomorphisms are, we can copy (nearly) all standard differential geometric constructions, such as manifolds (using an atlas and local diffeomorphisms for chart-changing), fiber bundles, vector bundles, the tangent bundle of a manifold, its cotangent bundle etc{\ae}tera. 
There is however one important difference that concerns vector bundles. 
An $\CA$-manifold is modelled by charts that are open subsets in the \stress{even} part of a finite dimensional $\CA$-vector space $E$; as such it has an even and an odd dimension (the dimensions of $(\body E)_0$ and $(\body E)_1$). 
On the other hand, the typical fiber of a vector bundle is \stress{all} of a finite dimensional $\CA$-vector space $F$. 
There are multiple reasons for doing so, one of them being that the standard constructions such as tensor products and exterior powers can be performed on $\CA$-vector spaces, but not on their even parts (at least not in a way that gives satisfactory results). 
Now a full $\CA$-vector space $F$ can be seen in a natural way as the even part of another $\CA$-vector space: the direct sum of $F$ with its parity dual $\prod F$ (the even part of $\prod F$ is the odd part of $F$ and the odd part of $\prod F$ is the even part of $F$), so vector bundles are still $\CA$-manifolds.

Concerning these constructions, one can show the following results. 
\begin{itemize}
\item
The body map can be extended to $\CA$-manifolds and their smooth maps, the result being an ordinary (non-super) manifold and an ordinary smooth map (taking the body map is essentially mapping all nilpotent elements in $\CA$ to zero).

\item
If $M$ is an $\CA$-manifold of graded dimension $p\vert q$, then $\body M$ is an ordinary manifold of dimension $p$ and for any (sufficiently small) open set $O\subset \body M$ with $U=\body\mo O\subset M$ we have the equality
$$
C^\infty(U, \CA) = C^\infty(O)\otimes \bigwedge \RR^q
\mapob,
$$
providing the link with the sheaf theoretic/ringed spaces approach (in this way we create a ringed space on the ordinary manifold $\body M$, and conversely, every ringed space appears this way).

\item
The set $C^\infty(M)$ of all smooth maps $f:M\to \CA$ on an $\CA$-manifold $M$ is a graded $\RR$-algebra and the set $\Gamma(B)$ of all smooth sections of a vector bundle $B\to M$ is a graded bi-module over $C^\infty(M)$.

\item
The tangent map $Tf:TM\to TN$ associated to a smooth map $f:M\to N$ between $\CA$-manifolds is naturally \stress{left}-linear and even. In particular we have, for local coordinate systems $x$ in $M$ and $y$ in $N$ with $y=f(x)$, the formula
$$
\iota\Bigl({\sum_{i=1}^{p+q}X_i \cdot \fracp{}{x_i}\bigrestricted_x}\Bigr){Tf} = \sum_{i=1}^{p+q}\sum_{j=1}^{p'+q'} X_i \cdot (\partial_if_j)(x)\cdot \fracp{}{y_j}\bigrestricted_y
\mapob,
$$
where $p\vert q$ is the graded dimension of $M$ and $p'\vert q'$ that of $N$.

\item
Any smooth and \stress{even} vector field $X$ on an $\CA$-manifold $M$ can be integrated to produce a flow, i.e., a $1$-parameter group of local diffeomorphsims $\Phi_t$ with $t\in \CA_0$, which has the usual properties of a flow.
According to \recalt{structureofsmoothfunctions} it has a local expression $\Phi_t(x,\xi)\equiv \Phi(t,x,\xi) = \sum_I (\mathbf{G}\Phi_I)(t,x)\cdot \xi^I$ in terms of ordinary smooth functions of real variables and products of odd coordinates. 
These functions can be computed by induction on the number of elements in $I$ (details can be found in \cite[Ch.V\S4]{Tu04}); the equation for $\Phi_\emptyset$ is a ``standard'' first order differential equation for the flow of a non-super vector field (it is the flow equation for the ordinary smooth vector field $\body X$ on the ordinary manifold $\body M$) and the other terms are determined by first order inhomogeneous linear differential equations (and as such they do not restrict the domain of definition of the flow).

\item
A smooth \stress{odd} vector field $X$ can be integrated to a flow $\Phi_\tau$ with an odd time parameter $\tau\in \CA_1$ (but otherwise with the same properties of a flow and in particular the group law) if and only if the auto-commutator $[X,X]$ is zero.
In local coordinates the expression for the flow essentially boils down to the formula
$$
\Phi_\tau(x) = x+\tau\cdot X(x)
\mapob.
$$

\item
Lie's third theorem and its converse are true: to any $\CA$-Lie group $G$ (i.e., an $\CA$-manifold equipped with a group structure such that multiplication and inversion are smooth maps) is associated an $\CA$-Lie algebra $\Liealg g$ (isomorphic to the tangent space at the identity and isomorphic to the space of all left-invariant vector fields). 
Conversely, for any finite dimensional $\CA$-Lie algebra $\Liealg g$ there exists a (unique up to isomorphisms) simply connected $\CA$-Lie group $G$ whose $\CA$-Lie algebra is (isomorphic to) $\Liealg g$.

\item
The notion of an $\CA$-Lie group is equivalent to the notion of a \quote{Lie supergroup pair,} \ie, an \quote{ordinary} Lie group $H$, an $\CA$-Lie algebra $\Liealg g$ and an action of $H$ on $\Liealg g$ that satisfy the two conditions that the Lie algebra $\Liealg h$ of $H$ is isomorphic to $\body\Liealg g_0$ and such that the restriction to $\Liealg h\cong\body\Liealg g_0$ of the action of $H$ on $\Liealg g$ is the adjoint action.
Starting from an $\CA$-Lie group $G$, one obtains the Lie supergroup pair directly by taking $H=\body G$ and as the action of $H$ on $\Liealg g$ the restriction to $H\subset G$ of the adjoint action of $G$ on $\Liealg g$.
Conversely, if we have a Lie supergroup pair $(H,\Liealg g)$, we first construct the simply connected $\CA$-Lie group $\widetilde G$ whose $\CA$-Lie algebra is $\Liealg g$. It then follows that $\body \widetilde G$ is a covering group of $H$.
If we denote by $Z\subset \body \widetilde G\subset \widetilde G$ the kernel of this covering map, then the condition that the action of $H$ on $\Liealg g$ extends the adjoint action of $H$ on $\Liealg h\cong \body\Liealg g_0$ implies that $Z$ is central (and discrete) in $\widetilde G$. 
One then shows easily that $G=\widetilde G/Z$ is an $\CA$-Lie group whose associated Lie supergroup pair (via the previous construction) is the given Lie supergroup pair.

\end{itemize}

\begin{definition}{Remark}
If we replace (everywhere) in \S\ref{subsecsmoothfunc} the field of real numbers $\RR$ by the field of complex numbers $\CC$ (which means in particular that we replace $\CA$ by $\CA^\CC = \CA\otimes_\RR\CC = \CA \oplus i\CA$), we get complex smooth super functions. As ordinary smooth complex differentiable functions (from open sets in $\CC^p$ to complex vector spaces) are holomorphic, this means that for $U\subset E_o$ open in an $\CA^\CC$-vector space of dimension $p\vert q$ with $O= \body U\subset \CC^p$ and $F$ any $\CA^\CC$-vector space, we have the equality
$$
C^\infty(U,F) = \operatorname{Hol}(O,\body F) \otimes \bigwedge \CC^q
\mapob,
$$
where $\operatorname{Hol}(O,\body F)$ denotes the set of complex analytic\slash holomorphic functions on $O$ with values in $F$.
As for real analytic superfunctions (and despite the fact that I don't have an intrinsic definition in the same spirit as for smoothness), it is easy to define them as a subset of super smooth real functions, just by requiring the functions $f_I$ in \recalt{structureofsmoothfunctions} to be real analytic.

Applying this to manifold theory \S\ref{subsecAmfds}, we obtain the categories of real and complex analytic $\CA^{(\CC)}$-manifolds (split and non-split). 
As the results on integration of vector fields depends (only) upon the corresponding statement in ordinary non-super analysis, one obtains directly that an analytic vector field on a real or complex analytic $\CA^{(\CC)}$-manifold admits an analytic flow.
It follows that an $\CA$-Lie group admits a unique structure of a real analytic $\CA$-manifold, simply because in the non-super setting one can prove it by means of the flow of an analytic vector field (see \cite{DK00}).

\end{definition}

\subsection{Foliations}
\vrule width0pt
\smallskip
\label{sectionappendixfoliations}

Before we attack the notion of a foliation of an $\CA$-manifold, we first have to discuss a negative side effect of the condition imposed on super smooth functions: that it maps the body of the source to the body of the target. 
On the positive side, this condition allows us to prove the structure of super smooth functions \recalt{structureofsmoothfunctions} and the fact that the body of an $\CA$-manifold is a subset of the $\CA$-manifold.
However, on the negative side we have the fact that fixing some coordinates in a super smooth function does not necessarily yield a super smooth function!

To explain the situation, we consider first a function $f:\RR^2\to \RR$, we take a fixed value $x\in \RR$ and we define the function $g:\RR\to \RR$ by
$$
g(y) = f(x,y)
\mapob.
$$
If $f$ is smooth, then $g$ is smooth. 
Let us next consider a function $f:(\CA_0)^2\to \CA$ and a fixed value $x\in \CA_0$, with which we define the function $g=\CA_0\to \CA$ again by $g(y) = f(x,y)$.
But now there is no guarantee that $g$ is smooth if $f$ is smooth!
The problem stems from the condition that for a smooth function we required the condition $g(\body U)\subset \body F$, which gives here the condition
\begin{equation}\label{restrictionbeingreal}
y\in \RR=\body \CA_0
\quad\Rightarrow\quad
g(y)=f(x,y)\in \RR=\body\CA
\mapob.
\end{equation}
But, since $f$ is smooth, we only know the property
$$
x,y\in \RR=\body\CA_0 
\quad\Rightarrow\quad
f(x,y)\in \RR=\body\CA
\mapob.
$$
It follows that if $x$ belongs to $\RR=\body\CA_0$, \stress{then} we can be sure that $g$ is super smooth. 
For any other $x$ it is highly unlikely that $g$ satisfies \recalf{restrictionbeingreal}, in which case it is not super smooth.
A very simple example (which can be interpreted as addition in the $\CA$-Lie group $\CA_1$ of dimension $0\vert1$) is the map $f:(\CA_1)^2\to \CA_1$ given by
$$
f(\xi,\eta) = \xi+\eta
\mapob.
$$
This map is smooth, but for any fixed non-zero element $\xi\in \CA_1$, the resulting map $g:\CA_1\to \CA_1$ given by
$$
g(\eta) = \xi+\eta
$$
is not smooth, as it maps $\{0\}=\body\CA_1$ to $\{\xi\}$, which is not the body of $\CA_1$.

More generally, fixing some coordinates to \stress{real} values in a super smooth function yields a super smooth function (of less variables), but when we restrict to non-real values (arbitrary values in $\CA_0$ or $\CA_1$), it is highly unlikely that the resulting function is super smooth.
For \stress{even} coordinates this is not a real problem, as a super smooth function of even coordinates is an extension of an ordinary smooth function of the same number of real coordinates (see \recalt{Gextension}).
On the other hand, for odd coordinates it is essential: if a super smooth function depends upon an odd coordinate, it is \stress{not} possible to fix this coordinate to an arbitrary value in $\CA_1$ (except $0$) and getting a super smooth function in the remaining coordinates.
This should not come as a surprise for readers familiar with the sheaf-theoretic\slash ringed spaces approach to supermanifolds. 
There a smooth (odd) coordinate function is a single element in the ring of local smooth superfunctions, it is not a function in the ordinary sense. 
Hence evaluating such a ``function'' in a point does not make sense.
On the other hand, we have shown that the local rings of smooth superfunctions are the same in the $\CA$-manifold approach and in the sheaf-theoretic\slash ringed spaces approach, so we should not be able to create more smooth superfunctions by restricting a superfunction to some (arbitrary) value.

One of the consequences of this fact is the following: the tangent map $Tf:TM\to TN$ associated to a smooth map $f:M\to N$ between $\CA$-manifolds is again smooth. However, the restriction to a specific tangent space
$$
Tf\caprestricted_{T_mM}:T_mM \to T_{f(m)}N
$$
will in general be smooth only when $m$ belongs to the body $\body M$ (which is equivalent to saying that $m$ has real coordinates in any and then all local coordinate systems). 
On the other hand, it retains all its algebraic properties, in particular it still is (left-) linear.
And even when it is not (allowed to be called) smooth, we still can speak about its partial derivatives! It suffices to take the partial derivatives \stress{before} fixing $m\in M$, and then to fix $m$.
Readers familiar with the sheaf-theoretic\slash ringed spaces approach to supermanifolds will recognize this in a different form. 
For there indeed $Tf$ is a smooth map, but it can not be applied to all individual tangent spaces, one has to consider the full map $Tf$ on the whole of $TM$ in order to be able to apply it to an arbitrary tangent vector $X\in T_mM$, seen as an arbitrary point in $TM$.

Another instance where this phenomenon becomes important is when one wants to compute the exponential map from a Lie superalgebra $\Liealg g$ to its Lie group $G$. 
As the $\CA$-Lie algebra $\Liealg g$ of an $\CA$-Lie group $G$ of dimension $p\vert q$ ($p$ even coordinates and $q$ odd coordinates) has dimension $p\vert q$, i.e., $p$ even generators and $q$ odd generators, it is ``immediate'' that $G$ is modelled on the even part $\Liealg g_0$. 
The exponential map thus should be a map $\exp:\Liealg g_0\to G$, for which one wants to integrate all even left-invariant vector fields associated to all even Lie algebra elements.
However, most of these will not be smooth and in particular some will be odd smooth vector fields multiplied by an odd parameter.
The ``solution'' to this problem is to consider \stress{all} these even vector fields simultaneously as a single even smooth vector field $Z$ on the bigger $\CA$-manifold $\Liealg g_0 \times G$ given by
$$
Z(X,g) = X_g
\mapob,
$$
where $X_g$ is the value of the left-invariant vector field on $G$ associated to $X\in \Liealg g_0\cong (T_{id}G)_0$ at the point $g\in G$.
As an even smooth vector field, $Z$ has a flow $\Phi_t$ and the exponential map then is defined by the equation
$$
\bigl(X,\exp(X)\bigr) = \Phi_1(X,id)\in \Liealg g \times G
\mapob.
$$
The map $\exp:\Liealg g_0\to G$ is smooth because $\Phi$ is and because $1$ and $id$ have real coordinates.

\begin{definition}{Definitions}
$\bullet$
A \stresd{subbundle of rank $k$} of the tangent bundle $TM$ of an $\CA$-manifold is a subset $D\subset TM$ satisfying the condition
\begin{itemize}[\quad]
\item
for all $m_o\in M$ there exists an open neighbourhood $U\subset M$ of $m_o$ and $k$ homogeneous smooth vector fields $X_1, \dots X_k$ on $U$ such that for all $m\in U$  the subset $D_m=D\cap T_mM$ is the subspace of $T_mM$ generated by the $k$ independent (a condition on the $X_i$) elements $X_1(m), \dots, X_k(m)$.

\end{itemize}

$\bullet$
A subbundle $D$ of rank $k$ is said to be \stresd{involutive} if for any two smooth vector fields $X$ and $Y$ on $M$ such that $X(m),Y(m)\in D_m$ for all $m\in M$ we have the property that $[X,Y](m)\in D_m$ (again for all $m\in M$). An involutive subbundle is also called \stresd{a foliation}.

$\bullet$
An \stresd{integral manifold of a subbundle $D$} is an $\CA$-manifold $N$ together with an immersion $i:N\to M$ (i.e., a smooth injective map for which $Ti$ is injective at all points, which makes sense even for points outside the body of $N$ as it is an algebraic condition) such that $(Ti)(T_nN) \subset D_{i(n)}$ for all $n\in N$.

\end{definition}

\begin{proclaim}[Frobenius]{Frobenius' theorem}
Let $D$ be a subbundle of rank $k$ of $TM$. Then $D$ is involutive if and only if for all $m_o\in M$ there exists a local system of coordinates $x_1, \dots, x_{p+q}$ around $m_o$ (even and odd together) such that $D_m$ is generated by $\partial_1\caprestricted_m, \dots, \partial_k\caprestricted_m$ (for all $m$ in the coordinate neighbourhood).

\end{proclaim}

\begin{proclaim}{Theorem}
Let $D\subset TM$ be a foliation of rank $k$ and let $m_o\in \body M$ be a point with real coordinates. 
Then there exists an integral manifold $i:N\to M$ for $D$ passing through $m_o$ (i.e., $m_o\in i(N)$) which is maximal, i.e., $(Ti)(T_nN) = D_{i(n)}$ for all $n\in N$ and  if $i':N'\to M$ is another integral manifold for $D$ passing through $m_o$, then $i'(N')\subset i(N)$. $N$ is usually called a \stresd{leaf} of the foliation $D$.

\end{proclaim}

We now come to yet another negative side effect. In non-super geometry, one shows that a subbundle is involutive if and only if through \stress{every} point of $M$ passes an integral manifold of maximal dimension, whereas for $\CA$-manifolds, we can only ascertain that for points with real coordinates (i.e., in the body of $M$).
To explain why, let us consider the very easy case of $M=(\CA_0)^2$ with two even coordinates $x$ and $y$, and the subbundle $D$ spanned by $\partial_x$.
This coordinate system thus automatically satisfies Frobenius' theorem \recalt{Frobenius} and it would seem obvious that the maps $i_y:\CA_0\to M=(\CA_0)^2$ defined as
$$
i_y(x) = (x,y)
$$
are integral manifolds (and leaves). 
Unfortunately, when $y$ does not belong to the body of $\CA_0$, the map $i_y$ does not belong to the class of smooth maps!
The obvious candidates for the leaves of this foliation thus can not all be called immersed submanifolds.
One can, as is done in the main part of this paper, define the leaves as topological submanifolds, and then through every point passes a topological leaf.
However, most of these leaves do not have the status of an immersed submanifold.

\section*{Acknowledgements}

I am indebted to the anonymous referee for her\slash his critical and constructive remarks, which helped enormously to improve the readability, correctness and scope of this paper.

This work was supported in part by the Labex CEMPI  (ANR-11-LABX-0007-01).

\bibliographystyle{amsalpha}

\bibliography{BiblioGMT}

\end{document}